\pdfoutput=1
\RequirePackage{ifpdf}
\ifpdf 
\documentclass[pdftex]{sigma}
\else
\documentclass{sigma}
\fi

\numberwithin{equation}{section}

\newtheorem{Theorem}{Theorem}[section]
\newtheorem*{Theorem*}{Theorem}
\newtheorem{Corollary}[Theorem]{Corollary}
\newtheorem{Lemma}[Theorem]{Lemma}
\newtheorem{Proposition}[Theorem]{Proposition}
\newtheorem{Conjecture}[Theorem]{Conjecture}
\newtheorem*{Claim}{Claim}

\theoremstyle{definition}
\newtheorem{Definition}[Theorem]{Definition}

\newtheorem{Example}[Theorem]{Example}
\newtheorem{Remark}[Theorem]{Remark}

\def\Ad{\operatorname{Ad}}
\def\tr{\operatorname{tr}}
\def\End{\operatorname{End}}

\def\Stab{\operatorname{Stab}}
\def\Mat{\operatorname{Mat}}

\def\diag{\operatorname{diag}}
\def\Id{\operatorname{Id}}
\def\id{\operatorname{id}}
\def\rdet{\operatorname{rdet}}
\def\cdet{\operatorname{cdet}}
\def\symdet{\operatorname{symdet}}
\def\sp{\operatorname{span}}
\def\sgn{\operatorname{sgn}}
\def\wt{\operatorname{wt}}
\def\DKir{\mathbf{D}}
\def\Func{\operatorname{Func}}

\def\Spec{\operatorname{Spec}}
\def\stirling{\genfrac\{\}{0pt}{}}

\usepackage{dsfont}
\usepackage{mathrsfs}
\usepackage{mathtools}
\usepackage[all,cmtip]{xy}
\usepackage{stmaryrd}

\begin{document}

\allowdisplaybreaks

\newcommand{\arXivNumber}{2501.04605}

\renewcommand{\PaperNumber}{024}

\FirstPageHeading

\ShortArticleName{Big Algebra in Type $A$ for the Coordinate Ring of the Matrix Space}

\ArticleName{Big Algebra in Type $\boldsymbol{A}$ for the Coordinate Ring \\ of the Matrix Space}

\Author{Nhok Tkhai Shon NGO}

\AuthorNameForHeading{N.T.S.~Ngo}

\Address{Institute of Science and Technology Austria (ISTA),\\ Am Campus 1, 3400 Klosterneuburg, Austria}
\Email{\mail{nhoktkhaishon.ngo@ista.ac.at}}

\ArticleDates{Received April 30, 2025, in final form February 22, 2026; Published online March 14, 2026}

\Abstract{In this paper, we consider the big algebra recently introduced by Hausel for the $\mathrm{GL}_n$-action on the coordinate ring of the matrix space $\mathrm{Mat}(n,r)$. In particular, we~obtain explicit formulas for the big algebra generators in terms of differential operators with polynomial coefficients. We show that big algebras in type $A$ are commutative and relate them to the Bethe subalgebra in the Yangian $\operatorname{Y}(\mathfrak{gl}_{n})$. We apply these results to big algebras of symmetric powers of the standard representation of $\mathrm{GL}_n$.}


\Classification{17B10; 16T25}

\section{Introduction}
Recently, Hausel introduced in \cite{Hausel_24} the notion of the \emph{big algebra} associated to an irreducible representation $V(\lambda)$ of a complex semi-simple Lie group $G$. In fact, the corresponding big algebra $\mathscr{B}(\lambda)$ can be seen as a certain maximal commutative subalgebra of the algebra of all $G$-equivariant polynomial maps from $\mathfrak{g}=\operatorname{Lie}(G)$ to $\End V(\lambda)$.

From the geometric point of view, big algebras were motivated by the study of the geometry of affine Schubert varieties. Namely, big algebras are isomorphic to the equivariant intersection cohomology of affine Schubert varieties. In particular, this fact can be used to define a~ring structure on equivariant intersection cohomology for such varieties, which is not possible in general.
From the algebraic side, big algebras resemble objects that are well-studied in mathematical physics and representation theory. For example, Gaudin spin chain models and Mischenko--Fomenko integrable systems, see \cite{Rybnikov}. Moreover, big algebras give a novel approach to representations of complex semi-simple Lie groups. For instance, Hausel in \cite{Hausel_24} suggests how one can read off from the spectrum of the big algebra a lot of representation-theoretic data, such as the weight diagram, the crystal structure and Lusztig's $q$-weight multiplicity polynomial.\looseness=1

In this paper, we work with big algebras from a purely algebraic point of view. We study the big algebras associated to polynomial representations of the general linear group $G=\mathrm{GL}_{n}$.
We~obtain certain explicit formulas for the generators of the big algebra in this case. Using these, we give an independent proof of the commutativity of big algebras by relating them to \emph{Bethe subalgebras} of the \emph{Yangian} $\operatorname{Y}(\mathfrak{gl}_{n})$. As an application of our results, we use them to study the big algebras associated to symmetric powers of the vector representation of $\mathrm{GL}_{n}$. Before we state our results more precisely, we first review the definition of big algebras.\looseness=1

\subsection{Big algebras}
Let $G$ be a connected complex reductive group with Lie algebra $\mathrm{Lie}(G)=\mathfrak{g}$ (note that in~\cite{Hausel_24}~$\mathfrak{g}$~is semi-simple, but most constructions can be extended to the reductive case as well). Fix a~triangular decomposition $\mathfrak{g}=\mathfrak{n}_{-}\oplus\mathfrak{h}\oplus\mathfrak{n}_{+}$, where $\mathfrak{h}$ is a Cartan subalgebra and $\mathfrak{n}_{\pm}$ are nilpotent subalgebras spanned by positive and negative root spaces. In particular, we fix a~choice of a~set of positive roots of $\mathfrak{g}$ and this defines the set of dominant weights.
For any dominant integral weight~$\lambda$, let $V(\lambda)$ be the irreducible representation of the highest weight $\lambda$ and denote by~${\pi_{\lambda}\colon\mathfrak{g}\to \End V(\lambda)}$ be the corresponding Lie algebra homomorphism.

Although in \cite{Hausel_24} big algebras are defined only for irreducible representations, it will be useful for us to define them in a more general setting.
Let $\pi\colon\mathfrak{g}\to\End V$ be a representation of~$\mathfrak{g}$ which is isomorphic to a direct sum of finite-dimensional representations. Denote by $S(\mathfrak{g}^*)$ the symmetric algebra of $\mathfrak{g}^*$, which we identify with the coordinate ring $\mathbb{C}[\mathfrak{g}]$ of polynomial functions on $\mathfrak{g}$. Then, the \emph{Kirillov algebra} \cite{Kirillov_00} of the representation $(\pi,V)$ is defined\footnote{Note that Kirillov in \cite{Kirillov_00} studied the algebra $(S(\mathfrak{g})\otimes\End V)^{G}$ which can be identified with $\mathscr{C}(V)$ since $\mathfrak{g}\simeq\mathfrak{g}^*$ as $\mathfrak{g}$-modules. We decided to use $S(\mathfrak{g}^*)$ instead of $S(\mathfrak{g})$ because the former is more natural for geometric applications. For instance, for the $G$-equivariant cohomology of a point we have $H_{G}(\mathrm{pt})\simeq\mathbb{C}[\mathfrak{g}]^{G}\simeq S(\mathfrak{g}^*)^{G}$.} as $\mathscr{C}(V)=(S(\mathfrak{g}^{*})\otimes\End V)^{G}$ (here we use the adjoint action of $G$ on $S(\mathfrak{g}^{*})$). Alternatively, we can view the elements of the Kirillov algebras as polynomial maps $F\colon\mathfrak{g}\to\End V$ which satisfy the following equivariance condition:
\begin{equation*}
F(\Ad(g)(X))=\pi(g)F(X)\pi(g)^{-1}, \qquad X\in\mathfrak{g},\ g\in G.
\end{equation*}
Note that $\mathscr{C}(V)$ is an \smash{$S(\mathfrak{g}^{*})^{G}$}-algebra since \smash{$S(\mathfrak{g}^{*})^{G}$} embeds into $\mathscr{C}(V)$ as the subalgebra of scalar operators, i.e., operators of the form $P\cdot\Id_{V}$ with \smash{$P\in S(\mathfrak{g}^{*})^{G}$}.

In general, Kirillov algebras need not be commutative.
However, one can construct a large commutative subalgebra inside $\mathscr{C}(V)$ and this is the \emph{big algebra} $\mathscr{B}(V)$ which we define next. The subalgebra $\mathscr{B}(V)$ can be defined by means of a certain ``differential-like'' operator $\DKir=\DKir_{V}$ acting on $\mathscr{C}(V)$ which we call the \emph{Kirillov--Wei operator} \cite[Section~1.2]{Kirillov_00} and \cite[equations~(42) and~(57)]{Wei}.

Fix a non-degenerate invariant symmetric bilinear form on $\mathfrak{g}$ (e.g., the Killing form in the case of simple $\mathfrak{g}$) and let \smash{$\{X_{i}\}_{i=1}^{\dim\mathfrak{g}}$} and \smash{$\{X^{i}\}_{i=1}^{\dim\mathfrak{g}}$} be two bases of $\mathfrak{g}$ dual with respect to this form. Then, for any $F\in\mathscr{C}(V)$ define
\begin{equation*}
(\DKir F)(X)=\sum_{i=1}^{\dim\mathfrak{g}}\frac{\partial F}{\partial X_{i}}(X)\cdot\pi\bigl(X^{i}\bigr),
\end{equation*}
where \smash{$\frac{\partial F}{\partial X_i}$} denotes the directional derivative of $F\colon\mathfrak{g}\to\End V$ along $X_i\in\mathfrak{g}$. It turns out that in types $A$, $B$, $C$, $D$ and $G$ there exist homogeneous generators $c_1,\dots,c_r$ of the ring \smash{$S(\mathfrak{g}^{*})^{G}$} such that the subalgebra generated by elements $\DKir^{p}(c_{q})$, is commutative (see \cite[Section~2.1]{Hausel_24} and Definition \ref{dfn:gl_n_big_alg} for $\mathfrak{g}=\mathfrak{gl}_n$).
The resulting commutative subalgebra of $\mathscr{C}(V)$ is called the \emph{big algebra} of $V$ and is denoted by $\mathscr{B}(V)$. Namely, we set
\begin{equation*}
\mathscr{B}(V)= \langle \DKir^{p}(c_q)\mid 1\le q\le r,\, 0\le p\le \deg(c_q) \rangle\subset\mathscr{C}(V).
\end{equation*}
For any dominant integral weight $\lambda$, we define $\mathscr{B}(\lambda)\coloneqq\mathscr{B}(V(\lambda))$ and $\mathscr{C}(\lambda)\coloneqq\mathscr{C}(V(\lambda))$.

More generally, big algebras can be defined for any simple Lie algebra $\mathfrak{g}$ in a different way,
using the notion of the \emph{Feigin--Frenkel center} (see \cite[Section~2]{Hausel_24} and also \cite[Section~4]{Rybnikov}). Let~$\hat{\mathfrak{g}}$ be the affine Kac--Moody algebra associated to $\mathfrak{g}$. As a Lie algebra, $\hat{\mathfrak{g}}$ is a central extension of the loop algebra $\mathfrak{g}\bigl[t,t^{-1}\bigr]$. Then, the Feigin--Frenkel center $\mathfrak{z}(\hat{\mathfrak{g}})$ can be regarded as a certain commutative subalgebra of the universal enveloping algebra $U(\hat{\mathfrak{g}}_{-})$, where $\hat{\mathfrak{g}}_{-}=t^{-1}\mathfrak{g}\bigl[t^{-1}\bigr]$ is the negative part of $\hat{\mathfrak{g}}$. The big algebra $\mathscr{B}(V)$ now can be defined as the image of $\mathfrak{z}(\hat{\mathfrak{g}})$ under certain homomorphism $\varrho_{z}\colon U(\hat{\mathfrak{g}}_{-})\to S(\mathfrak{g})\otimes \End V\simeq S(\mathfrak{g}^{*})\otimes\End V$ which depends on a non-zero complex parameter $z$. The homomorphism $\varrho_{z}$ acts on $\hat{\mathfrak{g}}_{-}$ as follows:
\begin{equation*}
\varrho_{z}(x\otimes t^{-k})=z^{-k} \cdot 1\otimes \pi(x)+\delta_{k,1}\cdot x\otimes 1,\qquad x\in\mathfrak{g},\ k=1,2,\dots\,.
\end{equation*}
One can show that the image $\varrho_{z}(\mathfrak{z}(\hat{\mathfrak{g}}))$ does not depend on the choice of $z\in\mathbb{C}^{\times}$, see \cite[Section~4, Corollary~3]{Rybnikov}.
This approach, however, does not give a simple way to compute big algebras since the formulas for the generators of $\mathfrak{z}(\hat{\mathfrak{g}}_{-})$ are rather complicated (see \cite[Chapters~7 and~8]{Molev_18}). Besides that, the commutativity of big algebras then relies on deep results of Feigin and Frenkel~\cite{Feigin_Frenkel} on $\mathfrak{z}(\hat{\mathfrak{g}}_{-})$. One of the motivations for our work was to obtain a more elementary proof of the commutativity of big algebras (in type~$A$) using the definition involving the Kirillov--Wei operator.

\subsection{Main results}

In this paper, we consider big algebras for $\mathfrak{g}=\mathfrak{gl}_{n}$. Note that Hausel considered big algebras only for semi-simple Lie algebras, but the construction can be applied in our case as well. We~fix a~positive integer $r$ and consider the coordinate ring $\mathcal{P}(n,r)=\mathbb{C}[\Mat(n,r)]$ of the affine space $\Mat(n,r)$ of complex $n\times r$ matrices.
The natural action of $\mathrm{GL}_{n}\times\mathrm{GL}_{r}$ on $\Mat(n,r)$ endows~$\mathcal{P}(n,r)$ with the structure of $(\mathrm{GL}_{n}\times\mathrm{GL}_{r})$-module.
Moreover, one can explicitly describe the decomposition of $\mathcal{P}(n,r)$ into the direct sum of irreducible $(\mathrm{GL}_n\times\mathrm{GL}_r)$-representations. Namely, by the \emph{Howe duality} \cite[Section~2.1.2]{Howe_95}, we have
\begin{equation*}
\mathbb{C}[\Mat(n,r)]=\mathcal{P}(n,r)\simeq\bigoplus_{\lambda}V_{\mathrm{GL}_n}(\lambda)\otimes V_{\mathrm{GL}_r}(\lambda),
\end{equation*}
where $\lambda$ runs over all partitions with length $\ell(\lambda)\leq\min\{r,n\}$, where $\lambda$ is regarded as a dominant weight for $\mathrm{GL}_n$ and $\mathrm{GL}_r$, respectively. Note that if $r=n$, then $\mathcal{P}(n,r)=\mathcal{P}(n,n)$ contains all polynomial irreducible representations of $\mathrm{GL}_n$.

In particular, we can regard $\mathcal{P}(n,r)$ as a $\mathfrak{gl}_n$-module and can consider the big algebra $\mathscr{B}(\mathcal{P}(n,r))$. It is not difficult to verify that the elements of $\mathscr{B}(\mathcal{P}(n,r))$ can be viewed as certain polynomial differential operators on $\Mat(n,r)$.

We find explicit formulas for the generators of $\mathscr{B}(\mathcal{P}(n,r))$ using direct calculations and the Kirillov--Wei operator $\DKir$ (Theorem~\ref{big_generator_thm} and Corollary~\ref{cor:F_operators}). Using these formulas and the Capelli identities, we identify the generators of $\mathscr{B}(\mathcal{P}(n,r))$ with homomorphic images of some elements of the Bethe subalgebra of the Yangian $\operatorname{Y}(\mathfrak{gl}_{n})$. As a consequence, we obtain that the big algebra $\mathscr{B}(\mathcal{P}(n,r))$ is commutative (Theorem~\ref{thm:big_alg_commut}). It follows from the big algebra construction that for any $V(\lambda)$ appearing in the decomposition of $\mathcal{P}(n,r)$ there is a surjective algebra homomorphism $\mathscr{B}(\mathcal{P}(n,r))\twoheadrightarrow\mathscr{B}(\lambda)$. Therefore, the discussion above implies the following fact.

\begin{Theorem}[Corollary~\ref{cor:big_alg_commut_irrep}]
The big algebra of any polynomial finite-dimensional irreducible representation of $\mathrm{GL}_{n}$ is commutative.
\end{Theorem}

The approach outlined above has an advantage of considering big algebras of all polynomial representations of $\mathfrak{gl}_n$ simultaneously. Besides that, our proof of commutativity relies only on direct calculations and thus, avoids the use of more complicated constructions such as the Feigin--Frenkel center.

Using the formulas for generators of big algebra, we can give a more explicit description of the big algebra $\mathscr{B}(m\varpi_1)$ for any positive integer $m$ (Proposition~\ref{prop:big_alg_sym_power_descr} and Theorem~\ref{thm:big_sym_power_gen_rel}). Here~$\varpi_1$ is the first fundamental weight of $\mathfrak{g}=\mathfrak{gl}_{n}$ which corresponds to the $n$-dimensional vector representation $V(\varpi_1)$. The representation $V(m\varpi_1)$ is isomorphic to the $m$-th symmetric power~$S^m(V(\varpi_1))$.

The representations $V(m\varpi_1)$, $m=1,2,\dots$, of $\mathfrak{gl}_{n}$ are particular examples of \emph{weight multiplicity free representations}. Recall that a representation $V(\lambda)$ of $\mathfrak{g}$ is called weight multiplicity free, if all summands in its weight decomposition
\begin{equation*}
V(\lambda)=\bigoplus_{\mu}V_{\mu}(\lambda), \qquad \text{where}~V_{\mu}(\lambda)=\{v\in V(\lambda) \mid \pi_{\lambda}(h)v=\mu(h)v~\text{for all}~h\in\mathfrak{h}\},
\end{equation*}
are at most one-dimensional. In the weight multiplicity free case, the Kirillov algebra was studied by several authors, see, e.g., \cite{Kirillov_00,Kirillov_01,Panyushev,Rozhkovskaya}. In particular, it is known that $\mathscr{C}(\lambda)$ is commutative (see also Section~\ref{subsect:kirillov_wmf}). Using the result conjectured by Hausel (see\footnote{In \cite{Hausel_24} this claim was stated as Theorem~1.1.3. However, we were informed by Hausel that, unfortunately, the original argument contained a gap, see Remark~\ref{rem:conj_status} for a discussion.} Conjecture~\ref{conj:kirillov_alg_center}), one can then deduce that in the weight multiplicity free case the big algebra $\mathscr{B}(\lambda)$ coincides with the Kirillov algebra $\mathscr{C}(\lambda)$. In particular, the conjecture of Hausel implies that $\mathscr{B}(m\varpi_1)=\mathscr{C}(m\varpi_1)$ for all $m$.

The description of $\mathscr{B}(m\varpi_1)$ from Proposition~\ref{prop:big_alg_sym_power_descr} also allows us to prove purely algebraically the following fact, first established geometrically by Hausel \cite[equation~(4.3)]{Hausel_24} using the technique of \emph{fixed point schemes} from \cite{Hausel_Rychlewicz}.

\begin{Theorem}[Theorem~\ref{thm:big_alg_sym_pow_iso}]
There exists an isomorphism of $S(\mathfrak{gl}_{n}^{*})^{\mathfrak{gl}_{n}}$-algebras
$\mathscr{B}(m\varpi_1)\simeq S^m(\mathscr{B}(\varpi_1))$. Here $S^m(\mathscr{B}(\varpi_1))$ denotes the symmetric part of the $m$-th tensor power of $S(\mathfrak{gl}_{n}^{*})^{\mathfrak{gl}_{n}}$-algebra $\mathscr{B}(\varpi_1)$, i.e., the subalgebra $(\mathscr{B}(\varpi_1)^{\otimes m})^{\mathfrak{S}_m}\subset\mathscr{B}(\varpi_1)^{\otimes m}$ \big(the tensor products are taken over $S(\mathfrak{gl}_{n}^{*})^{\mathfrak{gl}_{n}}$\big).
\end{Theorem}

Under the assumption that $\mathscr{B}(m\varpi_1)=\mathscr{C}(m\varpi_1)$, the statement above is essentially equivalent to a special case\footnote{Namely, the case $(\mathrm{A}_n,m\varphi_1)$ in Panyushev's notation.} of a conjecture due to Panyushev, see \cite[Conjecture 6.2\,(2)]{Panyushev}.
The conjecture asserts that there exists a variety $X_{\lambda}$ whose $G$-equivariant cohomology is isomorphic to $\mathscr{C}(\lambda)$. Then, Theorem~\ref{thm:big_alg_sym_pow_iso} concerns the case $\lambda=m\varpi_1$. For $m=1$, the weight $\lambda=\varpi_1$ is minuscule and one takes $X_{\varpi_1}=G/P_{\varpi_1}\simeq\mathbb{P}^{n-1}$, and indeed \smash{$H^*_{\mathrm{GL}_n}\bigl(\mathbb{P}^{n-1}\bigr)\simeq\mathscr{C}(\varpi_1)$} (see \cite{Panyushev}). When $m$ is arbitrary, Panyushev suggests to consider \smash{$X_{m\varpi_1}=\bigl(\mathbb{P}^{n-1}\bigr)^{m}/\mathfrak{S}_m$}. The theorem above confirms the conjecture in this case since
\begin{align*}
H^*_{\mathrm{GL}_n}(X_{m\varpi_1})& \simeq \bigl(H^*_{\mathrm{GL}_n}\bigl(\mathbb{P}^{n-1}\bigr)^{\otimes m}\bigr)^{\mathfrak{S_m}}\\
& \simeq S^m(\mathscr{C}(\varpi_1)) =S^m(\mathscr{B}(\varpi_1))
 \simeq\mathscr{B}(m\varpi_1)=\mathscr{C}(m\varpi_1).
\end{align*}

Finally, we observe that the big algebra $\mathscr{B}(m\varpi_1)$ is invariant under the action of the Kirillov--Wei operator $\DKir$. In general this does not seem to be the case, but in our situation this would also be a consequence of equality $\mathscr{B}(m\varpi_1)=\mathscr{C}(m\varpi_1)$. We were able to derive a formula for $\DKir$ on $\mathscr{B}(m\varpi_1)$ in terms of the description from Proposition~\ref{prop:big_alg_sym_power_descr}, see Section~\ref{subsect:kirillov_wei_sym_power}. The formula for~$\DKir$ turns out to be related to some constructions from symmetric function theory and the ring of diagonal invariants. Although we do not yet have a conceptual explanation of these phenomena, we believe that these observations might be of independent interest.

\subsection{Related works and generalizations} The algebras of the form $\mathscr{C}(V)=(S(\mathfrak{g}^*)\otimes\End V)^{G}$ and $\mathscr{Q}(V)=(U(\mathfrak{g})\otimes\End V)^{G}$ were called \emph{classical and quantum family algebras} by Kirillov in \cite{Kirillov_00,Kirillov_01}. However, these objects were also considered before in earlier works, in particular, by Kostant \cite{Kostant_63,Kostant_75}, Bracken and Green \cite{Bracken_Green,Green}, Gould \cite{Gould} and others. For instance, Kostant in \cite{Kostant_75} uses the name \emph{strongly commuting ring} for algebras $\mathscr{Q}(V)$. Several papers including \cite{Bracken_Green,Gould,Green} are devoted to the study of \emph{characteristic identities} for special elements in $\mathscr{Q}(V)$. In view of this, a more recent work by Hausel \cite{Hausel_25} adopts the term \emph{Kostant algebra} for $\mathscr{Q}(V)$.

We believe that our approach to proving the commutativity of big algebras in type $A$ can be extended to other classical types. In types $B$, $C$ and $D$, one also have many polynomial representations and one might use the properties of the evaluation homomorphism and Bethe subalgebras for twisted Yangians (see~\cite{Molev_07} and~\cite{Nazarov_Olshanski}).

Using the work of \cite{Rybnikov}, one can construct in a similar way a large commutative subalgebra~$\mathscr{G}(V)$ in $\mathscr{Q}(V)=(U(\mathfrak{g})\otimes \End V)^{G}$ (see \cite[Section~2]{Hausel_24}). A more explicit construction of this algebra involves a variant of the Kirillov--Wei operator $\DKir$ on $\mathscr{Q}(V)$. We speculate that a direct proof of the commutativity of $\mathscr{G}(V)$ in this case would require some generalizations of Capelli identities, e.g., from \cite{Mukhin_Tarasov_Varchenko}. Note that a part of the algebra $\mathscr{G}(V)$ can be described very explicitly, see~\mbox{\cite[Theorem~1.2]{Hausel_25}}. Namely, the algebra $\mathscr{G}(V)$ contains the center of $\mathscr{Q}(V)$ and in \cite{Hausel_25} the author obtains a description of the center of $\mathscr{Q}(V(\lambda))$ as a certain quotient of $Z(\mathfrak{g})\otimes Z(\mathfrak{g})$, where $Z(\mathfrak{g})$ is the center of $U(\mathfrak{g})$. Similarly to $\mathscr{B}(V)$, the algebra $\mathscr{G}(V)$ also has a geometric interpretation via the affine Schubert varieties (see \cite[Theorem~3.1]{Hausel_24}).

Finally, let us mention that certain analogues of algebras $\mathscr{C}(V)$ and $\mathscr{Q}(V)$ and characteristic identities for their elements were also studied for quantum groups $U_q(\mathfrak{g})$ (e.g., in \cite{Gould_Zhang_Bracken,Gurevich_Pyatov_Saponov, Zhang_Gould_Bracken}). Besides that, there exists a version of Capelli identities for $U_q(\mathfrak{gl}_n)$ where the algebra of differential operators is replaced by the so-called \emph{reflection equation algebra} (see \cite{Gurevich_Saponov_Zaitsev, Jing_Liu_Molev,Zaitsev} for more details). It seems plausible that some of our constructions can be generalized to these settings as well.

\subsection{Contents}

Now let us briefly outline the contents of this paper.

In Section~\ref{s2}, we fix the notation and recall the necessary facts about the representation theory of $\mathfrak{gl}_n$ and its action on $\mathcal{P}(n,r)=\mathbb{C}[\Mat(n,r)]$.

In Section~\ref{s3}, we recall the notions of the Kirillov algebra and the big algebra. Then, we state the explicit formulas for the generators of big algebra of the coordinate ring of $\Mat(n,r)$ (see Theorem~\ref{big_generator_thm} and Corollary~\ref{cor:F_operators}).

Section~\ref{s4} is rather technical and is devoted to proofs of Theorem~\ref{big_generator_thm} and Corollary~\ref{cor:F_operators}.

Sections \ref{sect:capelli_identities}--\ref{s7} contain the proof of the commutativity of big algebras. In Section~\ref{sect:capelli_identities}, we recall the Capelli identity and its variants which we use in the proofs. In Section~\ref{s6}, we review the construction of a certain commutative subalgebra of $\operatorname{Y}(\mathfrak{gl}_n)$, called \emph{Bethe subalgebra} following Molev \cite[Section~1.14]{Molev_07}. In Section~\ref{s7}, we prove the commutativity of the big algebra (in type~$A$) using the explicit formulas obtained in Section~\ref{s3} (Corollary~\ref{cor:F_operators}) and the results from Sections~\ref{sect:capelli_identities} and~\ref{s6}.

In Section~\ref{sect:sym_pow_vect_rep}, we use the formulas obtained in Section~\ref{s3} to prove several results about the big algebras $\mathscr{B}(m\varpi_1)$ including Proposition~\ref{prop:big_alg_sym_power_descr}, which gives a description of $\mathscr{B}(m\varpi_1)$ in terms of certain functions on the weight diagram of $V(m\varpi_1)$, and Theorem~\ref{thm:big_alg_sym_pow_iso} on the isomorphism between $\mathscr{B}(m\varpi_1)$ and $S^m(\mathscr{B}(\varpi_1))$. We obtain a presentation of $\mathscr{B}(m\varpi_1)$ in terms of generators and relations (see Theorem~\ref{thm:big_sym_power_gen_rel} and Corollary~\ref{cor:big_sym_power_gen_rel}) and compare it to the previous computations of Hausel--Rychlewicz \cite{Hausel_Rychlewicz}, Rozhkovskaya
\cite{Rozhkovskaya} and Hausel \cite{Hausel_24} (see Section~\ref{subsect:examples}). Besides that, we also prove a different formula for the Kirillov--Wei operator $\DKir$ on $\mathscr{B}(m\varpi_1)$ (see Proposition~\ref{prop:kirillov_wei_operator_formula}).

Finally, in Appendix~\ref{appendix}, we prove Lemma~\ref{lemma:vanishing_lemma} which was used in the proof of Theorem~\ref{thm:big_alg_sym_pow_iso}.

\section{Notation and preliminaries}\label{s2}
Most of the proofs in this paper involve many direct calculations. To simplify the formulas, we introduce the following notation.

\subsection{Operations with tuples} For every positive integer $m$, we denote $[m]\coloneqq\{1,\dots,m\}$ and let $\mathfrak{S}_{m}$ be the symmetric group of~$[m]$. For any integer $k$ such that $0\le k\le m$, we define $\binom{[m]}{k}$ to be the set of all $k$-element subsets of~$[m]$ and $[m]^{\underline{k}}$ to be the set of all $k$-tuples which consist of $k$ distinct elements of $m$. Clearly,
\begin{equation*}
\#\binom{[m]}{k}=\binom{m}{k},\qquad \#[m]^{\underline{k}}=m^{\underline{k}},
\end{equation*}
where $m^{\underline{k}}$ is the so-called \emph{falling factorial},
\begin{equation*}
m^{\underline{k}}=m(m-1)\cdots(m-k+1)=k!\cdot\binom{m}{k}.
\end{equation*}
We also denote by $[m]^k=[m]^{[k]}$ the set of all $k$-tuples with entries in $[m]$. It is occasionally convenient to view a $k$-tuple $I=(i_1,\dots,i_k)$ as a function on the set $[k]=\{1,\dots,k\}$, namely, we set $I(l)=i_l$ for $l\in[k]$.

\subsubsection{Action of the symmetric group}
For any $k$-tuple $I=(i_1,\dots,i_k)$ and any permutation $\pi\in\mathfrak{S}_k$, define
\begin{equation}\label{sym_gr_action_tuple}
\pi I=\bigl(i_{\pi^{-1}(1)},\dots,i_{\pi^{-1}(k)}\bigr).
\end{equation}
Regarding $I$ as a function on $[k]$, we can alternatively write $\pi I=I\circ\pi^{-1}$.

\begin{Remark}\label{subsets_tuples_rem}
Using this group action, we can identify $\binom{[m]}{k}$ with quotient $\mathfrak{S}_k\backslash[m]^{\underline{k}}$. In particular, we will often regard a $k$-element subset as a $k$-tuple with arbitrarily chosen ordering (and in~these situations, the choice of ordering will be irrelevant).
\end{Remark}

\subsubsection{Sign functions for tuples}
For any $I,J\in[n]^{\underline{k}}$, define the generalized sign function as follows:
\begin{equation*}
\sgn\binom{I}{J}=
\begin{cases}
\sgn(\tau) &\text{if as sets}~I=J,
\\
0 &\text{otherwise},
\end{cases}
\end{equation*}
where in the first case $\tau$ is the unique permutation in $\mathfrak{S}_{k}$ that maps the $k$-tuple $J$ to $I$. In other words, $\tau$ satisfies $I=\tau\cdot J$.
\begin{Remark}
Clearly, this agrees with the usual sign of permutation in the case when both $I$ and $J$ are permutations of $(1,2,\dots,n)$. For instance, if $I=(1,\dots,n)$ and $J=(\tau(1),\dots,\tau(n))$ for $\tau\in\mathfrak{S}_n$, then
\begin{equation*}
\sgn\binom{I}{J}=\sgn
\begin{pmatrix}
1 & 2 & \cdots & n
\\
\tau(1) & \tau(2) & \cdots & \tau(n)
\end{pmatrix}
=\sgn(\tau).
\end{equation*}
\end{Remark}
Now assume that we have tuples $I_1,\dots,I_k$ and $J_1,\dots,J_k$ such that $I_l,J_l\in[n]^{\underline{p_l}}$ for each $l=1,\dots,k$ and some positive integers $p_1,\dots,p_k$. We assume additionally that $I_1,\dots,I_k$ are disjoint as sets and similarly for $J_1,\dots,J_k$. Denote also $p=p_1+\dots+p_k$. Then, we define
\begin{equation*}
\sgn
\begin{pmatrix}
I_1 & \cdots & I_k
\\
J_1 & \cdots & J_k
\end{pmatrix}
=\sgn\binom{I}{J},
\end{equation*}
where $I$ and $J$ are $p$-tuples obtained by concatenating $I_1,\dots,I_k$ and $J_1,\dots,J_k$, respectively. For example, if $k=2$, then
\begin{gather*}
I(s)=
\begin{cases}
I_1(s), &s\in\{1,\dots,p_1\},
\\
I_2(s-p_1), &s\in\{p_1+1,\dots,p_1+p_2\},
\end{cases}
\\
J(s)=
\begin{cases}
J_1(s), &s\in\{1,\dots,p_1\},
\\
J_2(s-p_1), &s\in\{p_1+1,\dots,p_1+p_2\}.
\end{cases}
\end{gather*}

In our calculations we also use another variant of the signature function. For any tuples $I_1,J_1\in[m]^{\underline{p}}$ and $I_2,J_2\in[m]^{\underline{q}}$ with $p\ge q$, define
\begin{equation*}
\varepsilon(I_1,J_1,I_2,J_2)=
\begin{cases}
\sgn(\tau_1\tau_2) &\text{if as sets}~I_1\setminus I_2=J_1\setminus J_2\in\binom{[m]}{p-q},
\\
0 &\text{otherwise}.
\end{cases}
\end{equation*}
Here, in the first case, $\tau_1$ and $\tau_2$ are elements of $\mathfrak{S}_p$ that satisfy the following equalities:
\begin{equation*}
\tau_1 I_1|_{[q]}=I_2,\qquad \tau_2 J_1|_{[q]}=J_2, \qquad\text{and}\qquad \tau_1 I_1|_{\{p+1,\dots,q\}}=\tau_2 J_1|_{\{p+1,\dots,q\}}.
\end{equation*}
Note that a pair $(\tau_1,\tau_2)$ is not defined uniquely in general. However, for any other such pair $(\tau_1',\tau_2')$, there exists an element $\sigma\in\mathfrak{S}_p$ which fixes each element of $[q]$ and such that $(\tau_1',\tau_2')=(\sigma\tau_1,\sigma\tau_2)$. In particular, $\sgn(\tau_1\tau_2)=\sgn(\tau_1'\tau_2')$ and, consequently, $\varepsilon(I_1,J_1,I_2,J_2)$ is well defined.

Observe that both generalized sign functions are skew-symmetric in the sense of the following lemma.
\begin{Lemma}
For any tuples $I,J\in[m]^{\underline{p}}$ and any permutations $\sigma,\tau\in\mathfrak{S}_p$, one has
\begin{equation*}
\sgn\binom{\sigma I}{\tau J}=\sgn(\sigma\tau)\sgn\binom{I}{J}.
\end{equation*}
Similarly, for any $I_1,J_1\in[m]^{\underline{p}}$, $I_2,J_2\in[m]^{\underline{q}}$ with $p\ge q$ and any permutations $\sigma_1,\tau_1\in\mathfrak{S}_p$, $\sigma_2,\tau_2\in\mathfrak{S}_{q}$, one has
\begin{equation*}
\varepsilon(\sigma_1I_1,\tau_1J_1,\sigma_2I_2,\tau_2J_2)=\sgn(\sigma_1\tau_1)\sgn(\sigma_2\tau_2)\cdot\varepsilon(I_1,J_1,I_2,J_2).
\end{equation*}
\end{Lemma}

\subsection{Non-commutative matrices} We often work with matrices whose entries are elements of non-commutative (associative) algebras. In particular, many computations involve non-commutative versions of determinants.

Let \smash{$A=[a_{ij}]_{i,j=1}^{N}$} be an $N\times N$ matrix with entries in a certain non-commutative algebra. We define for matrix $A$ its
\begin{itemize}\itemsep=0pt
\item
\emph{row determinant}:
\begin{equation*}
\rdet(A)=\sum_{\sigma\in\mathfrak{S}_{N}}\sgn(\sigma)a_{1,\sigma(1)}a_{2,\sigma(2)}\cdots a_{N,\sigma(N)},
\end{equation*}

\item
\emph{column determinant}:
\begin{equation*}
\cdet(A)=\sum_{\sigma\in\mathfrak{S}_{N}}\sgn(\sigma)a_{\sigma(1),1}a_{\sigma(2),2}\cdots a_{\sigma(N),N},
\end{equation*}

\item
\emph{symmetrized determinant}:
\begin{equation*}
\symdet(A)=\sum_{\sigma,\tau\in\mathfrak{S}_{N}}\sgn(\sigma\tau)a_{\sigma(1),\tau(1)}a_{\sigma(2),\tau(2)}\cdots a_{\sigma(N),\tau(N)}.
\end{equation*}
\end{itemize}
Observe that if the entries of $A$ do commute, then the row and the column determinants coincide with the usual one while for the symmetrized version one has $\symdet(A)=N!\cdot\det(A)$. Finally, note that the row (column) determinant is skew-symmetric with respect to columns (rows) while the symmetrized determinant is skew-symmetric with respect to both rows and columns.

For any $N\times N$ matrix $M$ and any tuples $I=(i_1,\dots,i_k)\in[N]^{k}$ and $J=(j_1,\dots,j_l)\in[N]^{l}$, we denote by $M_{IJ}$ the following $k\times l$ matrix:
\begin{equation*}
M_{IJ}=[M_{i_{\alpha},j_{\beta}}]_{\alpha\in[k],\beta\in[l]}=
\begin{bmatrix}
M_{i_{1}j_{1}} & \cdots & M_{i_{1}j_{l}}
\\
\vdots & \ddots & \vdots
\\
M_{i_{k}j_{1}} & \cdots & M_{i_{k}j_{l}}
\end{bmatrix}.
\end{equation*}
In the case, when the entries of $I$ and $J$ are strictly increasing, $M_{IJ}$ is a submatrix of $M$.

\subsection[The Lie algebra gl\_n]{The Lie algebra $\boldsymbol{\mathfrak{gl}_n}$}
Let $\mathfrak{g}=\mathfrak{gl}_{n}$ be the general linear Lie algebra of complex $n\times n$ matrices. We denote by $\{E_{ij}\}_{i,j=1}^{n}$ the standard basis of $\mathfrak{gl}_n$ consisting of matrix units. They satisfy the following commutation relations:
\begin{equation*}
[E_{ij},E_{kl}]=\delta_{jk}E_{il}-\delta_{li}E_{kj}.
\end{equation*}
Denote also by $\{y_{ij}\}_{i,j=1}^{n}$ the corresponding coordinates on $\mathfrak{gl}_n$.
In particular, we can view the algebra $\mathbb{C}[\mathfrak{gl}_n]=S(\mathfrak{gl}_n^*)$ as the polynomial ring $\mathbb{C}[y_{ij},1\leq i,j\leq n]$.

Let $\mathfrak{h}=\sp\{E_{ii}\}_{1\le i\le n}$ be the Cartan subalgebra consisting of diagonal matrices. Let $\mathfrak{n}_{+}=\sp\{E_{ij}\}_{1\le i<j\le n}$ and $\mathfrak{n}_{-}=\sp\{E_{ij}\}_{1\le j<i\le n}$ be the nilpotent subalgebras consisting of upper-triangular and lower-triangular matrices, respectively. Thus, we obtain a triangular decomposition $\mathfrak{g}=\mathfrak{n}_{-}\oplus\mathfrak{h}\oplus\mathfrak{n}_{+}$.

\subsection{The coordinate ring of the matrix space}
Consider the coordinate ring $\mathcal{P}(n,r)=\mathbb{C}[\Mat(n,r)]$ of the space $\Mat(n,r)$ of complex $n\times r$ matrices. Denote by $\{x_{ij}\mid 1\le i\le n,\,1\le j\le r\}$ the standard coordinates on $\Mat(n,r)$. Then, $\mathcal{P}(n,r)$ is the polynomial ring $\mathbb{C}[x_{ij},\,1\le i\le n,\, 1\le j\le r]$ in $rn$ variables. We denote by $\partial_{ij}$ the partial derivative with respect to the variable $x_{ij}$.

\subsubsection[The action of gl\_n]{The action of $\boldsymbol{\mathfrak{gl}_n}$}
The matrix space $\Mat(n,r)$ possesses the following $\mathrm{GL}_n$-action:
\begin{equation*}
(g,A)\mapsto \bigl(g^{-1}\bigr)^{t}\cdot A,\qquad g\in\mathrm{GL}_n,\ A\in\Mat(n,r).
\end{equation*}
This action induces a $\mathrm{GL}_n$-action on the coordinate ring $\mathcal{P}(n,r)=\mathbb{C}[\Mat(n,r)]$.
Denote by~\smash{$\widetilde{L}(g)$} the corresponding linear operator in $\End\mathcal{P}(n,r)$.
Differentiating this action along one-parameter subgroups in $\mathrm{GL}_n$ yields the infinitesimal $\mathfrak{gl}_n$-action on $\mathcal{P}(n,r)$, which we denote by~$L$. A direct calculation shows that
\begin{equation}\label{eq:L_representation}
L(E_{ij})=\sum_{\alpha=1}^{r}x_{i \alpha}\partial_{j \alpha},\qquad 1\le i,j\le n.
\end{equation}

\subsubsection[The algebra PD(n,r) of differential operators on P(n,r)]{The algebra $\boldsymbol{\mathcal{PD}(n,r)}$ of differential operators on $\boldsymbol{\mathcal{P}(n,r)}$}
Since $L$ is a representation of $\mathfrak{gl}_n$ it also gives rise to a representation of the universal enveloping algebra $U(\mathfrak{gl}_n)$ on $\mathcal{P}(n,r)$. In particular, the formula for $L(E_{ij})$ above implies that the elements of $U(\mathfrak{gl}_n)$ act on $\mathcal{P}(n,r)$ as differential operators with polynomial coefficients.

Let $\mathcal{PD}(n,r)$ be the (non-commutative) algebra of differential operators on $\mathcal{P}(n,r)$ with polynomial coefficients. In other words, $\mathcal{PD}(n,r)$ is the Weyl algebra generated by $\{x_{i\alpha},\partial_{i\alpha} \mid 1\le i\le n,\, 1\le\alpha\le r\}$ and relations of the form
\begin{equation*}
[\partial_{i\alpha},x_{j\beta}]=\delta_{ij}\delta_{\alpha\beta},\qquad 1\le i,j\le n,\ 1\le\alpha,\beta\le r.
\end{equation*}
We will need the following well-known fact in Section~\ref{sect:capelli_identities}.
\begin{Proposition}\label{prop:L_nn_injectivity}
If $r=n$, then the map $L\colon U(\mathfrak{gl}_n)\to\mathcal{PD}(n,r)$ is injective.
\end{Proposition}
\begin{Remark}
In fact, this follows from the observation that $U(\mathfrak{gl}_n)$ is isomorphic to the algebra of all left-invariant differential operators on $\mathrm{GL}_n$. In general, the image of $L$ coincides with the subalgebra of all differential operators in $\mathcal{PD}(n,r)$ which commute with the $\mathrm{GL}_r$-action on~$\mathbb{C}[\Mat(n,r)]$. See also \cite[Section~1]{Howe_Umeda}.
\end{Remark}

Denote by $E$, $X$ and $D$ the matrices, whose $(i,j)$-th entry equals $E_{ij}$, $x_{ij}$ and $\partial_{ij}$, respectively.
In particular, we have a formal identity $L(E)=X\cdot D^{t}$.

\section[Kirillov algebra, big algebra and medium algebra on P(n,r)]{Kirillov algebra, big algebra and medium algebra on $\boldsymbol{\mathcal{P}(n,r)}$}\label{s3}

\subsection[Invariant polynomials on gl\_n]{Invariant polynomials on $\boldsymbol{\mathfrak{gl}_n}$}\label{subsect:inv_polys}

Define the elements $c_1,\dots,c_n\in S(\mathfrak{gl}_n^*)^{\mathfrak{gl}_n}$ as the coefficients of the characteristic polynomial of~$Y$:
\begin{equation*}
\det(Y-z\cdot \Id_n)=\sum_{k=0}^{n}(-1)^{n-k}c_k(Y)\cdot z^{n-k}.
\end{equation*}
Denote by $y_{ij}$ the coordinates on $\mathfrak{gl}_n$ which correspond to standard matrix units $E_{ij}\in\mathfrak{gl}_n$. Then, $Y=[y_{ij}]_{i,j=1}^{n}$ and
\begin{equation}\label{eq:c_k_formula}
c_k(Y)=\sum_{I\in\binom{[n]}{k}}\det Y_{II}.
\end{equation}
It is well known that the elements $c_1,\dots,c_n$ are free generators of the ring $S(\mathfrak{gl}_n^*)^{\mathfrak{gl}_n}$ of $\mathrm{GL}_n$-invariants of $S(\mathfrak{gl}_n^*)$.

\subsection{Construction of the Kirillov algebra}
Let $\widetilde{\pi}\colon\mathrm{GL}_n\to V$ be a finite-dimensional representation of $\mathrm{GL}_n$ and let $\pi\colon\mathfrak{gl}_n\to\End V$ be the associated representation of the Lie algebra $\mathfrak{gl}_n$.
\begin{Definition}
The \emph{Kirillov algebra}\footnote{Also known as the \emph{classical family algebra}, see Kirillov's original papers \cite{Kirillov_00,Kirillov_01}.} of $V$
is defined as the algebra
\[
\mathscr{C}(V)=(S(\mathfrak{gl}_n^*)\otimes\End V)^{\mathrm{GL}_n}.
\]
\end{Definition}
In other words, the Kirillov algebra is the algebra of $\mathrm{GL}_n$-equivariant polynomial maps from~$\mathfrak{gl}_n$ to $\End V$, i.e., any $F\in\mathscr{C}(V)$ satisfies the following equivariance condition:
\begin{equation}\label{eq:equivar_condition}
F(\Ad(g)(Y))=\widetilde{\pi}(g)F(Y)\widetilde{\pi}(g)^{-1},\qquad g\in\mathrm{GL}_n,\ Y\in\mathfrak{gl}_n,
\end{equation}
Note that $\mathscr{C}(V)$ is an algebra over the ring $S(\mathfrak{gl}_n^*)^{\mathfrak{gl}_n}$. The elements of $S(\mathfrak{gl}_n^*)^{\mathfrak{gl}_n}$ are realized inside $\mathscr{C}(V)$ as scalar operators.

\begin{Remark}
If $V=\bigoplus_{\alpha} V_{\alpha}$ is an infinite direct sum of finite-dimensional representations $V_{\alpha}$, then we replace the algebra $\End V$ in the definition $\mathscr{C}(V)$ with the direct sum $\bigoplus_{\alpha,\beta}V_{\alpha}^*\otimes V_{\beta}$.
\end{Remark}

Define the \emph{Kirillov--Wei operator} $\DKir=\DKir_V$ which acts on $\mathscr{C}(V)$ as follows: for any $F\in\mathscr{C}(V)$, set
\begin{equation}\label{eq:kirillov_wei_op_def}
(\DKir F)(Y)=\sum_{i,j=1}^{n}\frac{\partial F}{\partial y_{ji}}(Y)\cdot \pi(E_{ij}).
\end{equation}
It follows from the definition that for any positive integer $p$, we have
\begin{equation*}
(\DKir^pF)(Y)=\sum_{i_1,\dots,i_p=1}^{n}\sum_{j_1,\dots,j_p=1}^{n}\frac{\partial^p F}{\partial y_{j_1i_1}\cdots\partial y_{j_pi_p}}(Y)\cdot \pi(E_{i_1j_1}\cdots E_{i_pj_p}).
\end{equation*}

In \cite[Lemma 1]{Kirillov_00} and \cite[Section~1.4]{Kirillov_01}, Kirillov hinted at the following fact.

\begin{Proposition}[{\cite[Proposition 5.2]{Wei}}]
The operator $\DKir$ maps $\mathscr{C}(V)$ to itself.
\end{Proposition}
\begin{proof}
This is shown via a direct verification of \eqref{eq:equivar_condition} for $F=\DKir G$, $G\in\mathscr{C}(V)$.
\end{proof}

\begin{Remark}
One can generalize the construction of the Kirillov algebra, the operator $\DKir$ and the results of this subsection for any semi-simple Lie algebra $\mathfrak{g}$. The corresponding operator $\DKir$ can be described using a non-degenerate invariant symmetric bilinear form on $\mathfrak{g}$ (e.g., the Killing form in the case of simple $\mathfrak{g}$).
\end{Remark}

The explicit formulas for specific elements of $\mathscr{C}(V)$ are often quite complicated. However, these formulas can be simplified if one restricts elements of Kirillov algebra to Cartan subalgebra~$\mathfrak{h}$. Moreover, these restrictions uniquely determine the elements of $\mathscr{C}(V)$.

\begin{Proposition}[{\cite[Section~1.2]{Kirillov_01} and \cite[Section~2]{Panyushev}}]\label{prop:restr_cartan}
The restriction map $\Phi=\Phi_{V}\colon F\mapsto F|_{\mathfrak{h}}$, $F\in\mathscr{C}(V)$ is an injective homomorphism of algebras over $S(\mathfrak{gl}_{n}^{*})^{\mathrm{GL}_{n}}\simeq S(\mathfrak{h}^{*})^{W}$. Moreover, the image of $\Phi$ is contained inside the algebra
\begin{equation*}
(S(\mathfrak{h}^{*})\otimes\End_{\mathfrak{h}}V)^{W}=\biggl(\bigoplus_{\mu\in\wt(V)}S(\mathfrak{h}^{*})\otimes\End V_{\mu}\biggr)^{W}.
\end{equation*}
Here $W\simeq\mathfrak{S}_{n}$ denotes the Weyl group of $\mathrm{GL}_{n}$, $\End_{\mathfrak{h}}V$ is the space of all $\mathfrak{h}$-equivariant linear operators on $V$, and $\wt(V)$ is the set of weights of $V$.
\end{Proposition}

\begin{proof}
We use the argument outlined in \cite[Theorem~2]{Kirillov_00}. It follows from the equality \eqref{eq:equivar_condition} that the restriction $F|_{\mathfrak{h}}$ completely determines the map $F\colon\mathfrak{gl}_n\to\End V$ on the set of all regular semi-simple elements of $\mathfrak{gl}_n$ (those are conjugate to elements of $\mathfrak{h}$). Since the latter set is Zariski dense in $\mathfrak{gl}_n$, the injectivity of $\Phi$ follows.

The second statement follows from the fact that the Weyl group $W$ is the normalizer of $\mathfrak{h}$ in~$\mathrm{GL}_{n}$ and the equivariance condition \eqref{eq:equivar_condition} applied for $g\in W$.
\end{proof}

The following statement shows that in order to study the Kirillov algebra of $V$ one can study the Kirillov algebra of a ``larger'' representation.
\begin{Proposition}\label{prop:kir_direct_sum_prop}
Let $V=\bigoplus_{\alpha} V_{\alpha}$ be a direct sum of $\mathrm{GL}_n$-modules. Then, the natural algebra homomorphism $\prod_{\alpha}\mathscr{C}(V_{\alpha})\to\mathscr{C}(V)$ induced by the embeddings $\iota_{\alpha}\colon\End(V_{\alpha})\hookrightarrow\End(V)$ is injective.
\end{Proposition}

\begin{proof}
The map in question sends a family of maps $F_{\alpha}\colon\mathfrak{gl}_n\to\End V_{\alpha}$ to the map $\sum_{\alpha}\iota_{\alpha}\circ F_{\alpha}\colon\mathfrak{gl}_n\to\mathscr{C}(V)$. Note that this expression is well defined since for every $Y\in\mathfrak{gl}_n$ the summand $(\iota_{\alpha}\circ F_{\alpha})(Y)$ belongs to $\iota_{\alpha}(\End V_{\alpha})$. In particular, this map is indeed an algebra homomorphism. The injectivity now follows from the injectivity of the maps $\iota_\alpha$.
\end{proof}

\subsection{Big and medium algebras}
First we recall the definitions of the big and the medium algebras in type~$A$.

\begin{Definition}\label{dfn:gl_n_big_alg}
Let $V$ be a direct sum of finite-dimensional $\mathrm{GL}_{n}$-modules. The \emph{big algebra}~$\mathscr{B}(V)$ is the subalgebra of $\mathscr{C}(V)$ generated by the elements
\begin{equation}\label{eq:M_pq_operators}
M_{p,q}=\DKir^{q}(c_{p+q}),\qquad 0\le p, q\le n,\ p+q\le n.
\end{equation}
The \emph{medium algebra} $\mathscr{M}(V)$ of $V$ is a subalgebra of $\mathscr{B}(V)$ defined as follows:
\begin{equation}\label{eq:medium_algebra_def}
\mathscr{M}(V)=\bigl\langle F,\DKir(F)\mid F\in S(\mathfrak{gl}_n^*)^{\mathfrak{gl}_n}\bigr\rangle\subset\mathscr{B}(V).
\end{equation}
\end{Definition}
\begin{Remark}
For a direct sum $V=\bigoplus_{\alpha} V_{\alpha}$ of finite-dimensional representation the elements of $\mathscr{B}(V)$ and $\mathscr{M}(V)$ belong to $S(\mathfrak{g}^*)\otimes \bigl(\bigoplus_{\alpha}\End V_{\alpha}\bigr)$ (this follows from the definition of $\DKir$, see \eqref{eq:kirillov_wei_op_def}). Therefore, the algebra structure on $\mathscr{B}(V)$ and $\mathscr{M}(V)$ is well defined even if $V$ is infinite-dimensional.
\end{Remark}

\begin{Remark}
A similar approach for defining big algebras can be used for complex semi-simple Lie groups of types $B$, $C$, $D$ and $G$, see \cite[Sections~1 and~2]{Hausel_24}.
\end{Remark}

Now let us give a more explicit formula for the generators of $\mathscr{B}(V)$. Recall that $ \pi\colon\mathfrak{gl}_n \to \End V $ is the representation of $\mathfrak{gl}_n$ which comes from the $\mathrm{GL}_n$-action on $V$.

\begin{Proposition}\label{symmetrized_m_pq_prop}
For any $Y\in\mathfrak{gl}_n$, we have
\begin{equation*}
M_{p,q}(Y)=\DKir^{q}(c_{p+q})(Y)=\sum_{\substack{I_1,J_1\in\binom{[n]}{p}\\I_2,J_2\in\binom{[n]}{q}\\I_1\sqcup I_2=J_1\sqcup J_2}}\sgn\binom{I_1~I_2}{J_1~J_2}\det Y_{I_1J_1}\cdot\symdet \pi(E)_{J_2I_2}.
\end{equation*}
Here $\pi(E)$ stands for the $n\times n$ matrix whose $(i,j)$-entry equals $\pi(E_{ij})$.
\end{Proposition}

\begin{proof}
This follows by a direct calculation using \eqref{eq:kirillov_wei_op_def} and \eqref{eq:c_k_formula}.
\end{proof}

The next proposition relates the big (resp.\ medium) algebras of direct sums with big (resp. medium) algebras of summands.

\begin{Proposition}\label{prop:big_alg_direct_sum_prop}
\quad
\begin{enumerate}
\item[$(i)$] The image of the algebra homomorphism $\prod_{\alpha}\mathscr{C}(V_{\alpha})\to\mathscr{C}(V)$ from Proposition~$\ref{prop:kir_direct_sum_prop}$ contains the big algebra $\mathscr{B}(V)$.

\item[$(ii)$] For each $\alpha$, there exist surjective homomorphisms $\mathscr{B}(V)\to\mathscr{B}(V_{\alpha})$ and $\mathscr{M}(V)\to\mathscr{M}(V_{\alpha})$. These homomorphisms are induced by the natural algebra homomorphisms $\prod_{\alpha'}\End V_{\alpha'}\to\End V_{\alpha}$.
\end{enumerate}
\end{Proposition}

\begin{proof}
The definition of the Kirillov--Wei operator implies the subalgebra $\prod_{\alpha}\mathscr{C}(V_{\alpha})$ of $\mathscr{C}(V)$ is invariant under the action of $\DKir_{V}$. Moreover, $\DKir_{V}$ acts on $\prod_{\alpha}\mathscr{C}(V_{\alpha})$ in the same way as the operator $\prod_{\alpha}\DKir_{V_{\alpha}}$. The first part now follows from the fact that $c_1,\dots,c_n\in S(\mathfrak{gl}_n^*)^{\mathfrak{gl}_n}$ are realized inside $\mathscr{C}(V)$ as scalar operators and the identity \smash{$\DKir^{p}_{V}(c_{k}\otimes\id_{V})=\prod_{\alpha}\DKir^{p}_{V_{\alpha}}(c_{k}\otimes\id_{V_{\alpha}})$} for all $p\geq 1$. The required homomorphisms $\mathscr{B}(V)\to\mathscr{B}(V_{\alpha})$ and $\mathscr{M}(V)\to\mathscr{M}(V_{\alpha})$ are induced by the projection homomorphisms $\prod_{\alpha'}\mathscr{C}(V_{\alpha'})\to\mathscr{C}(V_{\alpha})$.
\end{proof}

\subsection[Kirillov algebra for P(n,r)]{Kirillov algebra for $\boldsymbol{\mathcal{P}(n,r)}$}

There are natural actions of the groups $\mathrm{GL}_n$ and $\mathrm{GL}_r$ on the matrix space $\Mat(n,r)$ that induce actions of these groups on $\mathcal{P}(n,r)$.
Namely, for any $A\in\Mat(n,r)$, $g\in\mathrm{GL}_n$ and $h\in\mathrm{GL}_r$, we~have
\begin{equation*}
\widetilde{L}(g)(P)(A)=P\bigl(g^{t}\cdot A\bigr),\qquad \widetilde{R}(h)(P)(A)=P(A\cdot h).
\end{equation*}
Note that these two actions commute and give rise to an action of $\mathrm{GL}_n\times\mathrm{GL}_r$ on $\mathcal{P}(n,r)$.
Then, \emph{Howe duality} (see \cite[Section~2.1.2]{Howe_95}) implies that the $(\mathrm{GL}_n\times\mathrm{GL}_r)$-module $\mathcal{P}(n,r)=\mathbb{C}[\Mat(n,r)]$ decomposes as follows:
\begin{equation}\label{eq:p_gl_n_decomp}
\mathcal{P}(n,r)\simeq\bigoplus_{\lambda\colon \ell(\lambda)\le r}V_{\mathrm{GL}_n}(\lambda)\otimes V_{\mathrm{GL}_r}(\lambda),
\end{equation}
where the summation is over all partitions $\lambda=(\lambda_1,\dots,\lambda_n)$ of length at most $r$ (i.e., with~${\lambda_{i}=0}$ for $i>r$). Here $V_{\mathrm{GL}_n}(\lambda)$ and $V_{\mathrm{GL}_r}(\lambda)$ denote the irreducible representations of highest weight~$\lambda$ of $\mathrm{GL}_n$ and $\mathrm{GL}_r$, respectively.
In particular, we conclude that $\mathcal{P}(n,r)$ contains as subrepresentations all irreducible representations of $\mathfrak{gl}_n$ that have highest weight $\lambda$ with $\ell(\lambda)\leq r$.

\subsection[Generators of the big algebra of P(n,r)]{Generators of the big algebra of $\boldsymbol{\mathcal{P}(n,r)}$}

From now on, we consider the case $V=\mathcal{P}(n,r)$ and the corresponding Kirillov and big algebras. Proposition~\ref{prop:big_alg_direct_sum_prop} implies that many properties of big algebras of irreducible representations can be read off from $\mathscr{B}(\mathcal{P}(n,r))$.

Observe that the big algebra $\mathscr{B}(\mathcal{P}(n,r))$ is contained in $S(\mathfrak{gl}_n^*)\otimes\mathcal{PD}(n,r)$. In other words, the elements of $\mathscr{B}(\mathcal{P}(n,r))$ are certain differential operators on $\mathcal{P}(n,r)$ whose coefficients are polynomials in variables $x_{i\alpha}$ and $y_{jk}$. One of the main results of this paper is an explicit formula (the so-called \emph{normal form}) for the operators~$M_{p,q}$.

\begin{Theorem}\label{big_generator_thm}
The normal form of the operator $M_{p,q}(Y)$ is as follows:
\begin{align}\nonumber
M_{p,q}(Y)=
{}&\sum_{\ell=0}^{q}(-1)^{q-\ell}(q-\ell)!\,\ell!\,\stirling{q}{\ell}\sum_{\substack{I_1,J_1\in\binom{[n]}{p}\\I_2,J_2\in\binom{[n]}{q}\\I_1\sqcup I_2=J_1\sqcup J_2}}\sgn\binom{I_1~I_2}{J_1~J_2}\det Y_{I_1J_1}
\\
&{}{\times}\,
\sum_{R\in\binom{[r]}{\ell}}\sum_{V,W\in\binom{[q]}{\ell}}\varepsilon(J_2,I_2,J_2(V),I_2(W))\det\bigl(X_{J_2(V),R}\bigr)\det\bigl(D_{I_2(W),R}\bigr).\label{big_generator_formula}
\end{align}
Here, for tuples $I$ and $J$ we denote by $X_{IJ}$, $Y_{IJ}$ and $D_{IJ}$ the corresponding submatrices of \smash{$X=[x_{ij}]_{i,j=1}^{n,r}$}, \smash{$D=[\partial_{ij}]_{i,j=1}^{n,r}$} and \smash{$Y=[y_{ij}]_{i,j=1}^{n,n}$}, respectively.
\end{Theorem}
\begin{Remark}
$\stirling{k}{\ell}$ stands for the \emph{Stirling number of the second kind}, i.e., the number of ways to split a $k$-element set into $\ell$ non-empty subsets.
\end{Remark}

\begin{Corollary}\label{cor:F_operators}
The big algebra $\mathscr{B}(\mathcal{P}(n,r))$ is generated by the operators $\{F_{p,q}\mid p,q\ge 0,\allowbreak {p+q\le n}\}$, where
\begin{equation}\label{eq:F_operators_formula}
F_{p,q}(Y)=\sum_{\substack{I_1,J_1\in\binom{[n]}{p}\\I_2,J_2\in\binom{[n]}{q}\\I_1\sqcup I_2=J_1\sqcup J_2}}\sgn\binom{I_1~I_2}{J_1~J_2}\det Y_{I_1J_1}\sum_{R\in\binom{[r]}{q}}\det(X_{J_2,R})\det(D_{I_2,R}).
\end{equation}
In particular, the operators $\{M_{p,q}\mid p,q\ge 0,\, p+q\le n\}$ and $\{F_{p,q}\mid p,q\ge 0,\, p+q\le n\}$ are related to each other in the following way:
\begin{equation}\label{m_f_operators_transition}
M_{p,q}=\sum_{\ell=0}^{q}(-1)^{q-\ell}(q-\ell)!\,\ell!\,\stirling{q}{\ell}\binom{n-p-\ell}{q-\ell}F_{p,\ell}.
\end{equation}
\end{Corollary}
We prove these formulas for $M_{p,q}$ and $F_{p,q}$ in the next section.
The following fact is an~immediate consequence of the explicit formulas from Corollary~\ref{cor:F_operators}.

\begin{Corollary}\label{cor:F_generators_r}
The big algebra $\mathscr{B}(\mathcal{P}(n,r))$ is generated by operators $\{F_{p,q}\mid 1\le p\le n, \allowbreak {0\le q\le r},\, {p+q\le n}\}$.
\end{Corollary}
\begin{proof}
Indeed, it follows from \eqref{eq:F_operators_formula} that all operators $F_{p,q}$ with $q>r$ are zero.
\end{proof}

\subsection{Restriction to the Cartan subalgebra} The formulas \eqref{big_generator_formula} and \eqref{eq:F_operators_formula} simplify a lot if we restrict $Y\in\mathfrak{gl}_n$ to Cartan subalgebra $\mathfrak{h}$, i.e., to diagonal matrices.
We expect that, at least for small values of $r$, these formulas can be used for finding common eigenfunctions in $\mathbb{C}[\Mat(n,r)]$ of operators $\{F_{p,q}\}$ and the corresponding eigenvalues (for $r=1$ this is done in Section~\ref{sect:sym_pow_vect_rep}). This would be helpful for obtaining presentations of $\mathscr{B}(\lambda)$ in terms of generators and relations.
\begin{Proposition}\label{f_cartan_restr_prop}
For $Y=\diag(z_1,\dots,z_n)\in\mathfrak{h}$, we have
\begin{equation*}
F_{p,q}(Y)=\sum_{\substack{I\in\binom{[n]}{p}\\J\in\binom{[n]}{q}\\I\cap J=\varnothing}}\biggl(\prod_{i\in I}z_i\biggr)\biggl(\sum_{R\in\binom{[r]}{q}}\det(X_{JR})\det(D_{JR})\biggr).
\end{equation*}
\end{Proposition}
\begin{proof}
This is a consequence of Corollary~\ref{cor:F_operators} and the observation that for $Y\in\mathfrak{h}$ the determinant $\det Y_{I_1J_1}$ vanishes unless $I_1=J_1$.
\end{proof}

We use this formula in the proof of commutativity of the big algebra (see Theorem~\ref{thm:big_alg_commut}). More precisely, we match these normal form expressions with images of certain elements of the Bethe subalgebra using the Capelli identities (see Section~\ref{s7}).

\section{Proofs of Theorem~\ref{big_generator_thm} and Corollary~\ref{cor:F_operators}}\label{s4}

In this section, we prove Theorem~\ref{big_generator_thm} and Corollary~\ref{cor:F_operators}, i.e., formulas \eqref{big_generator_formula} and \eqref{eq:F_operators_formula}. The proof is purely computational and reduces to some identities in the Weyl algebra.

\subsection{Computation of the symmetrized determinants} In view of Proposition~\ref{symmetrized_m_pq_prop}, to prove Theorem~\ref{big_generator_thm}, we first need to obtain a formula for $\symdet L(E)_{IJ}$. Recall that for $k$-tuples $I=(i_1,\dots,i_k)$ and $J=(j_1,\dots,j_k)$ in $[n]^{\underline{k}}$ the symmetrized determinant $\symdet L(E_{IJ})$ is defined as (see \eqref{eq:L_representation})
\begin{equation*}
\begin{aligned}
\symdet L(E_{IJ})
&=\sum_{\sigma,\tau\in\mathfrak{S}_k}\sgn(\sigma\tau)L\bigl(E_{i_{\sigma(1)}i_{\tau(1)}}\bigr)\cdots L\bigl(E_{i_{\sigma(k)}i_{\tau(k)}}\bigr)
\\
&=\sum_{a_1,\dots,a_k=1}^{n}\sum_{\sigma,\tau\in\mathfrak{S}_{k}}\sgn(\sigma\tau)\cdot x_{i_{\sigma(1)}a_1}\partial_{j_{\tau(1)}a_1}\cdots x_{i_{\sigma(k)}a_k}\partial_{j_{\tau(k)}a_k}.
\end{aligned}
\end{equation*}
Our aim is to obtain a reduced expression (the normal form) for $\symdet L(E_{IJ})$ in the algebra~$\mathcal{PD}(n,r)$ of differential operators on $\mathfrak{gl}_n$ with polynomial coefficients.
The computation is divided into several steps.

We start with the following auxiliary identity.

\begin{Lemma}\label{sym_det_ident_lemma}
Let $N$ be a positive integer. For any permutations $\sigma,\tau\in\mathfrak{S}_{N}$ and any $s\in\{1,\dots,N\}$, define
\begin{equation*}
\delta_{s}(\sigma,\tau)=
\begin{cases}
0 &\text{if}\ \sigma(s)\le\tau(s),
\\
1 &\text{if}\ \sigma(s)>\tau(s).
\end{cases}
\end{equation*}
Then, for any $\alpha_1,\dots,\alpha_{N}\in\mathbb{C}$ the following identity holds:
\begin{gather}
\sum_{\sigma,\tau\in\mathfrak{S}_{N}}\sgn(\sigma\tau)\cdot(\alpha_1+\delta_1(\sigma,\tau))\cdots(\alpha_{N}+\delta_{N}(\sigma,\tau))\nonumber\\
\qquad{}=(-1)^{N-1}(N-1)!(\alpha_1+\dots+\alpha_{N}).\label{sym_det_alpha}
\end{gather}
\end{Lemma}

\begin{proof} The identity clearly holds for $N=1$, so let us assume that $N\ge 2$.
Observe that the left-hand side can be expressed as
\begin{equation*}
\sum_{\ell=0}^{N}C_{\ell}\cdot e_{\ell}(\alpha_1,\dots,\alpha_N),
\end{equation*}
where $e_{\ell}(\alpha_1,\dots,\alpha_N)$ is the $\ell$-th elementary symmetric polynomial in $\alpha_1,\dots,\alpha_N$ and $C_0,\dots,C_{\ell}$ are certain real numbers. Therefore, to check the identity \eqref{sym_det_alpha} it suffices to calculate the coefficients $C_{\ell}$ for every $0\le\ell\le N$. Note that $C_{\ell}$ is equal to the coefficient of the term $\alpha_1\dots\alpha_{\ell}$
on the left-hand side. Then, it is not difficult to see that
\begin{equation*}
C_{\ell}=\sum_{\sigma,\tau\in\mathfrak{S}_N}\sgn(\sigma\tau)\prod_{s=\ell+1}^{N}\delta_{s}(\sigma,\tau).
\end{equation*}
In other words, $C_{\ell}$ is the sum of $\sgn(\sigma\tau)$, where $(\sigma,\tau)$ runs over all pairs of permutations in $\mathfrak{S}_N$ such that
\begin{equation*}
\sigma(s)>\tau(s)\qquad\text{for all}~\ell+1\le s\le N.
\end{equation*}
Let $\Gamma_{\ell}\subset\mathfrak{S}_{N}\times\mathfrak{S}_{N}$ be the set of all such pairs.
Now let us consider several cases:
\begin{itemize}\itemsep=0pt
\item
$\ell=0$. In this case, we have
\begin{equation*}
C_{0}=\sum_{\sigma,\tau\in\mathfrak{S}_{N}}\sgn(\sigma\tau)\cdot\delta_1(\sigma,\tau)\cdots\delta_N(\sigma,\tau)=0
\end{equation*}
because for any $\sigma,\tau\in\mathfrak{S}_{N}$ at least one of $\delta_{s}(\sigma,\tau)$ is zero \big(for example, one can take $s=\sigma^{-1}(1)$\big).

\item
$\ell=1$. We claim that in this case the set $\Gamma_{\ell}=\Gamma_{1}$ contains exactly $(N-1)!$ elements. Indeed, by definition a pair $(\sigma,\tau)\in\mathfrak{S}_{N}\times\mathfrak{S}_{N}$ belongs to $\Gamma_{1}$ if $\sigma(s)>\tau(s)$ for every $s\in\{2,\dots,N\}$. One checks that this holds if and only if these permutations satisfy $\sigma(1)=1$, $\tau(1)=N$ and $\sigma(s)=\tau(s)+1$ for all $s=2,\dots,N$. In particular, $|\Gamma_{1}|=(N-1)!$ and for any $(\sigma,\tau)\in\Gamma_{1}$ the permutation $\sigma\tau^{-1}$ is the cycle $(1~2~\dots~N)$. Hence,
\begin{equation*}
C_{1}=\sum_{(\sigma,\tau)\in\Gamma_{1}}\sgn(\sigma\tau)=|\Gamma_{1}|\cdot(-1)^{N-1}=(-1)^{N-1}\cdot(N-1)!.
\end{equation*}

\item
$2\le\ell\le N$. Note that for any permutations $\sigma',\tau'\in\mathfrak{S}_{N}$ which fix each of $\ell+1,\dots,N$ the pair $(\sigma,\tau)$ belongs to $\Gamma_{\ell}$ if and only if $(\sigma\sigma',\tau\tau')\in\Gamma_{\ell}$. This and the equality
\begin{equation*}
\sum_{\sigma',\tau'\in\mathfrak{S}_{\ell}}\sgn(\sigma'\tau')=\biggl(\sum_{\sigma'\in\mathfrak{S}_{\ell}}\sgn\sigma'\biggr)\biggl(\sum_{\tau'\in\mathfrak{S}_{\ell}}\sgn\tau'\biggr)=0
\end{equation*}
imply that $C_{\ell}=0$ for all $2\le\ell\le N$.

\end{itemize}
Combining everything, we obtain
\begin{equation*}
C_{\ell}=
\begin{cases}
0,&~\ell\neq 1,
\\
(-1)^{N-1}(N-1)!,&~\ell=1,
\end{cases}
\end{equation*}
which is equivalent to \eqref{sym_det_alpha}.
\end{proof}

\begin{Corollary}\label{sym_det_ident_cor}
Let $N$ be a positive integer and let $\ell\in\{0,1,\dots,N\}$. Then, in the notation of Lemma $\ref{sym_det_ident_lemma}$ we have
\begin{gather*}
\sum_{\sigma,\tau\in\mathfrak{S}_{N}}\sgn(\sigma\tau)\prod_{j=1}^{\ell}(\alpha_{j}+\delta_{j}(\sigma,\tau))\\
\qquad{}=
\begin{cases}
0,&\ell\in\{0,1,\dots,N-2\},
\\
(-1)^{N-1}(N-1)!,&\ell=N-1,
\\
(-1)^{N-1}(N-1)!\cdot(\alpha_1+\dots+\alpha_N),&\ell=N.
\end{cases}
\end{gather*}
\end{Corollary}
\begin{proof}
This follows from \eqref{sym_det_alpha} by differentiating both sides with respect to $\alpha_{\ell+1},\dots,\alpha_N$.
\end{proof}

Next, we prove the following identity in the Weyl algebra which essentially computes the symmetrized determinant $\symdet L(E)_{IJ}$ in the simplest case $r=1$.

\begin{Lemma}\label{sym_det_weyl_lemma}
Consider the Weyl algebra generated by variables $u_1,\dots,u_n$ and the corresponding partial derivatives $\partial_1,\dots,\partial_n$. Let $I=(i_1,\dots,i_k)$ and $J=(j_1,\dots,j_k)$ be two $k$-tuples in $[n]^{\underline{k}}$. Assume that $I$ and $J$ have $\ell$ common elements, $0\le\ell\le k$. Denote
\begin{equation*}
\Psi(I,J)=\symdet\bigl(\bigl[u_{i_{\alpha}}\partial_{j_{\beta}}\bigr]_{\alpha,\beta=1}^{k}\bigr)=\sum_{\sigma,\tau\in\mathfrak{S}_k}\sgn(\sigma\tau)u_{i_{\sigma(1)}}\partial_{j_{\tau(1)}}\cdots u_{i_{\sigma(k)}}\partial_{j_{\tau(k)}}.
\end{equation*}
Then, $\Psi(I,J)$, as an element of the Weyl algebra, can be simplified as follows:
\begin{itemize}\itemsep=0pt
\item
if $\ell=k$ and $\sigma\in\mathfrak{S}_k$ is such that $i_{l}=j_{\sigma(l)}$ for all $l\in\{1,\dots,k\}$, then
\[
\Psi(I,J)=(-1)^{k-1}(k-1)!\cdot\sgn(\sigma)\sum_{i\in I}u_i\partial_{i},
\]

\item
if $\ell=k-1$ and $\pi,\sigma\in\mathfrak{S}_k$ are such that $i_{\pi(l)}=j_{\sigma(l)}$ for all $l\in\{1,\dots,k-1\}$ and $i_{\pi(k)}\neq j_{\sigma(k)}$, then
\[
\Psi(I,J)=(-1)^{k-1}(k-1)!\cdot\sgn(\pi\sigma)u_{i_{\pi(k)}}\partial_{j_{\sigma(k)}},
\]

\item
if $\ell\le k-2$, then $\Psi(I,J)=0$.
\end{itemize}
Moreover, viewing $I$ and $J$ as functions on $[k]=\{1,\dots,k\}$, one can express $\Psi(I,J)$ as follows:
\begin{equation*}
\Psi(I,J)=(-1)^{k-1}(k-1)!\sum_{\pi\in\mathfrak{S}_k}\sum_{s=1}^{k}\mathds{1}\bigl(I|_{[k]\setminus\{s\}}=\pi^{-1}J|_{[k]\setminus\{s\}}\bigr)\cdot\sgn(\pi)u_{i_s}\partial_{j_{\pi(s)}}.
\end{equation*}
Here, for an assertion $A$ we put $\mathds{1}(A)=1$, if $A$ holds, and $\mathds{1}(A)=0$ otherwise.
\end{Lemma}

\begin{proof}
We start with the following observation: for any permutation $\pi\in\mathfrak{S}_k$, we have
\[
\Psi(\pi I,\sigma J)=\sgn(\pi\sigma)\cdot\Psi(I,J),
\]
where $\pi I=\bigl(i_{\pi^{-1}(1)},\dots,i_{\pi^{-1}(k)}\bigr)$ and $\sigma J=\bigl(j_{\sigma^{-1}(1)},\dots,j_{\sigma^{-1}(k)}\bigr)$. Therefore, it suffices to prove the statement in the case when $i_1=j_1$, \dots, $i_{\ell}=j_{\ell}$ and
\begin{equation*}
\{i_1,\dots,i_k\}\cap\{j_1,\dots,j_k\}=\{i_1,\dots,i_{\ell}\}.
\end{equation*}

Denote $K=(i_1,\dots,i_{\ell})=(j_1,\dots,j_{\ell})
$. To simplify the notation, let us also assume that $K=(1,\dots,\ell)$. Now note that elements $u_{i_{\ell+1}},\dots,u_{i_k}$ and $\partial_{j_{\ell+1}},\dots,\partial_{j_k}$ commute with each other and also with $u_i$ and $\partial_{i}$ for $i\in K$. Therefore, we can rewrite $\Psi(I,J)$ as
\begin{equation}\label{psi_i_j_reduced}
\Psi(I,J)=u_{i_{\ell+1}}\cdots u_{i_{k}}\partial_{j_{\ell+1}}\cdots\partial_{j_k}\cdot\Phi(\ell;K),
\end{equation}
where $\Phi(\ell;K)$ is a certain element in the Weyl algebra generated by $u_i$, $\partial_i$ with $i\in K$. Namely, $\Phi(\ell;K)$ is obtained from the expression for $\Phi(I,J)$ by removing all $u_i$ and $\partial_i$ with $i\notin K$. One finds that the action of $\Phi(\ell;K)$ on a monomial $u_{1}^{\alpha_1}\cdots u_{\ell}^{\alpha_{\ell}}$, where $\alpha_1,\dots,\alpha_{\ell}\in\mathbb{N}_0$, is given by
\begin{gather*}
\Phi(\ell;K)\bigl(u_{1}^{\alpha_1}\cdots u_{\ell}^{\alpha_{\ell}}\bigr)\\
\qquad=\biggl(\sum_{\sigma,\tau\in\mathfrak{S}_{k}}\sgn(\sigma\tau)\bigl(\alpha_{1}+\delta_1\bigl(\sigma^{-1}, \tau^{-1}\bigr)\bigr)\cdots\bigl(\alpha_{\ell}+\delta_{\ell}\bigl(\sigma^{-1},\tau^{-1}\bigr)\bigr)\biggr)\cdot u_{1}^{\alpha_1}\cdots u_{\ell}^{\alpha_{\ell}},
\end{gather*}
where the $\delta_i$'s are defined as in Lemma~\ref{sym_det_ident_lemma}. Indeed, this follows from the identities
\begin{gather*}
(u_i\partial_i)\bigl(u_{1}^{\alpha_1}\cdots u_{\ell}^{\alpha_{\ell}}\bigr)=\alpha_i\cdot u_{1}^{\alpha_1}\cdots u_{\ell}^{\alpha_{\ell}},\\
(\partial_i u_i)\bigl(u_{1}^{\alpha_1}\cdots u_{\ell}^{\alpha_{\ell}}\bigr)=(\alpha_i+1)\cdot u_{1}^{\alpha_1}\cdots u_{\ell}^{\alpha_{\ell}},\qquad 1\le i\le\ell.
\end{gather*}
Now observe that
\begin{equation*}
\sum_{\sigma,\tau\in\mathfrak{S}_{k}}\sgn(\sigma\tau)\prod_{j=1}^{\ell}\bigl(\alpha_{j}+\delta_{j}\bigl(\sigma^{-1},\tau^{-1}\bigr)\bigr)=\sum_{\sigma,\tau\in\mathfrak{S}_k}\sgn(\sigma\tau)\prod_{j=1}^{\ell}(\alpha_{j}+\delta_{j}(\sigma,\tau)),
\end{equation*}
and hence Corollary~\ref{sym_det_ident_cor} implies that
\begin{itemize}\itemsep=0pt
\item
if $0\le\ell\le k-2$, then $\Phi(\ell;K)\bigl(u_{1}^{\alpha_1}\cdots u_{\ell}^{\alpha_{\ell}}\bigr)=0$;

\item
if $\ell=k-1$, then $\Phi(\ell;K)\bigl(u_{1}^{\alpha_1}\cdots u_{\ell}^{\alpha_{\ell}}\bigr)=(-1)^{k-1}(k-1)!\cdot u_{1}^{\alpha_1}\cdots u_{\ell}^{\alpha_{\ell}}$;

\item
if $\ell=k$, then $\Phi(\ell;K)\bigl(u_{1}^{\alpha_1}\cdots u_{\ell}^{\alpha_{\ell}}\bigr)=(-1)^{k-1}(k-1)!\cdot(\alpha_1+\dots+\alpha_k) u_{1}^{\alpha_1}\cdots u_{\ell}^{\alpha_{\ell}}$.
\end{itemize}
Since the elements of Weyl algebra are uniquely defined by their action on the corresponding polynomial ring $\mathbb{C}[u_1,\dots,u_n]$, we obtain
\begin{equation*}
\Phi(\ell;K)=
\begin{cases}
0, & 0\le\ell\le k-2,
\\
(-1)^{k-1}(k-1)!, & \ell=k-1,
\\
(-1)^{k-1}(k-1)!\sum_{i\in K}u_i\partial_i, & \ell=k.
\end{cases}
\end{equation*}
Combining this with \eqref{psi_i_j_reduced} concludes the proof.
\end{proof}

\begin{Lemma}\label{lemma:delta_a_IJ}
Fix $a_1,\dots,a_k\in\{1,\dots,r\}$. For any two $k$-tuples $I=(i_1,\dots,i_k)$ and $J=(j_1,\dots,j_k)$ of distinct elements of $\{1,\dots,n\}$, denote
\begin{equation*}
\Delta(I,J)=\sum_{\sigma,\tau\in\mathfrak{S}_k}\sgn(\sigma\tau)\cdot x_{i_{\sigma(1)}a_1}\partial_{j_{\tau(1)}a_1}\cdots x_{i_{\sigma(k)}a_k}\partial_{j_{\tau(k)}a_k}.
\end{equation*}
Then, we have the following identity:
\begin{gather}\label{delta_a_IJ_formula}
\Delta(I,J)=(-1)^{k-\ell}(k-\ell)!\sum_{V,W\in\binom{[k]}{\ell}}\varepsilon(I,J,I(V),J(W)) \det\bigl(X_{I(V),R}\bigr)\cdot\det\bigl(D_{J(W),R}\bigr),
\end{gather}
where $R=\{a_1,\dots,a_k\}$ and $\ell=|R|$.
\end{Lemma}
\begin{Remark}
Here $I(V)$ and $J(W)$ are the subtuples of $I$ and $J$ that correspond to $V$ and $W$, respectively. The $k$-element subsets $V$ and $W$ are viewed as $k$-tuples with the ordering chosen in an arbitrary way (see Remark~\ref{subsets_tuples_rem}). Note that the summand does not depend on a choice of the ordering of $V$ and $W$.
\end{Remark}
\begin{Remark}
Observe that the formula for $\Psi(I,J)$ from Lemma~\ref{sym_det_weyl_lemma} is a particular case of \eqref{delta_a_IJ_formula} when $\ell=|R|=1$. Indeed, for any $I,J\in[n]^{\underline{k}}$ and any $l,s\in[k]$, we have
\begin{equation*}
\sum_{\substack{\pi\in\mathfrak{S}_k\\\pi(s)=l}}
\mathds{1}\bigl(I|_{[k]\setminus\{s\}}=\pi^{-1}J|_{[k]\setminus\{s\}}\bigr)\cdot\sgn(\pi)=\varepsilon(I,J,i_{s},j_{l}),
\end{equation*}
where in $\varepsilon(I,J,i_{s},j_{l})$ the elements $i_{s}$ and $j_{l}$ are viewed as $1$-tuples.
\end{Remark}
\begin{Remark}
In the notation of Section~\ref{sect:capelli_identities}, $\Delta(I,J)$ is a particular summand of the symmetrized row determinant $\rdet[L(E_{i_{\alpha}j_{\beta}})]_{\alpha,\beta=1}^{k}$.
When $I=J$, the formula \eqref{delta_a_IJ_formula} can be simplified as follows:
\begin{equation*}
\begin{aligned}
\Delta(I,I)=(-1)^{k-\ell}(k-\ell)!\sum_{V\in\binom{[k]}{\ell}}\det\bigl(X_{I(V),R}\bigr)\cdot\det\bigl(D_{I(V),R}\bigr).
\end{aligned}
\end{equation*}
Thus, we obtain an expression from the Cauchy--Binet type identity, see Proposition~\ref{prop:cauchy_binet_capelli} from Section~\ref{sect:capelli_identities}. Later, in the proof of Theorem~\ref{thm:big_alg_commut} we use these observations to match the formulas for big algebra generators with homomorphic images of elements of Bethe subalgebra.
\end{Remark}
\begin{proof}[Proof of Lemma~\ref{lemma:delta_a_IJ}]
Define the decomposition \smash{$\{1,\dots,k\}=\bigsqcup_{l=1}^{r}K_l$} via
\begin{equation*}
K_l=\{j\in\{1,\dots,k\}\mid a_{j}=l\}.
\end{equation*}
Let $m_l$ be the cardinality of $K_l$ for each $l$. Note that $R=\{l\in[r]\mid m_l>0\}$ and $\ell=|R|$. Denote by $\mathfrak{S}(K_l)$ the group of permutations of the set $K_l$ viewed as a subgroup of $\mathfrak{S}_k$. In particular, $\mathfrak{S}(K_l)\simeq\mathfrak{S}_{m_l}$ for all $l$. We also denote by $\mathfrak{S}(K)$ the subgroup of $\mathfrak{S}_k$ which stabilizes each of the subsets~$R_l$, i.e., $\mathfrak{S}(K)=\mathfrak{S}(K_1)\times\dots\times\mathfrak{S}(K_r)$.

Towards the end of the proof we will regard any sequence $\{c_j\}_{j\in K_l}$ indexed by elements of~$K_l$ as an $m_l$-tuple by considering the elements of $K_l$ in ascending order. Since we are working in a non-commutative algebra, let us make a convention that products over $K_l$ are considered in ascending order as well.

For any permutation $\sigma\in\mathfrak{S}_k$ and any $l\in[r]$ define the $m_l$-tuples $\sigma^{-1}I|_{K_l}$ and $\sigma^{-1}J|_{K_l}$ as follows (cf.\ \eqref{sym_gr_action_tuple}):
\begin{equation*}
\sigma^{-1}I|_{K_l}=\bigl\{i_{\sigma(s)}\bigr\}_{s\in K_l},\qquad \sigma^{-1}J|_{K_l}=\bigl\{j_{\sigma(s)}\bigr\}_{s\in K_l}.
\end{equation*}
For any two $p$-tuples $U=(u_1,\dots,u_p)$ and $V=(v_1,\dots,v_p)$ in $[k]^{\underline{p}}$, define
\begin{equation*}
\Theta_l(U,V)=\prod_{i=1}^{p}x_{u_i l}\partial_{v_i l}=x_{u_1 l}\partial_{v_1 l}\cdots x_{u_p l}\partial_{v_p l}.
\end{equation*}
We also denote (cf.\ Lemma~\ref{sym_det_weyl_lemma})
\begin{equation*}
\Psi_l(U,V)=
\begin{cases}
1, & l\notin R,
\\
\symdet\bigl(\bigl[x_{u_{\alpha}l}\partial_{v_{\beta}l}\bigr]_{\alpha,\beta=1}^{p}\bigr), & l\in R.
\end{cases}
\end{equation*}
Observe that for any $\pi_1,\pi_2\in\mathfrak{S}(K)$, we have
\begin{equation*}
\Psi_l\bigl((\sigma\pi_1)^{-1}I|_{K_l},(\tau\pi_2)^{-1}J|_{K_l}\bigr)=\sgn(\pi_1|_{K_l})\sgn(\pi_2|_{K_l})\cdot
\Psi_l\bigl(\sigma^{-1}I|_{K_l},\tau^{-1}J|_{K_l}\bigr),
\end{equation*}
where $\pi_i|_{K_l}$ is the restriction of permutation $\pi_i\in\mathfrak{S}(K)$ to the subset $K_l$ (recall that ${\pi_i(K_l)\!=\!K_l}$).
Finally, denote by $\mathfrak{S}_k/\mathfrak{S}(K)$ the set of all left cosets of $\mathfrak{S}(K)$ in $\mathfrak{S}_k$.

With all this notation, we can now proceed to the proof of \eqref{delta_a_IJ_formula}.
Since $x_{ia}$ and $\partial_{jb}$ commute whenever $a\neq b$, we can rewrite $\Delta(I,J)$ as follows:
\begin{equation*}
\begin{aligned}
\Delta(I,J)
&=\sum_{\sigma,\tau\in\mathfrak{S}_k}\sgn(\sigma\tau)\cdot x_{i_{\sigma(1)}a_1}\partial_{j_{\tau(1)}a_1}\cdots x_{i_{\sigma(k)}a_k}\partial_{j_{\tau(k)}a_k}
\\
&=\sum_{\sigma,\tau\in\mathfrak{S}_k}\sgn(\sigma\tau)\cdot\prod_{l=1}^{r}\Theta_l\bigl(\sigma^{-1}I|_{K_l},\tau^{-1}J|_{K_l}\bigr)
\\
&=\sum_{[\sigma],[\tau]\in\mathfrak{S}_k/\mathfrak{S}(K)}\sgn(\sigma\tau) \biggl(\sum_{\pi_1,\pi_2\in\mathfrak{S}(K)}\sgn(\pi_1\pi_2)\prod_{l=1}^{r}\Theta_l\bigl((\sigma\pi_1)^{-1}I|_{K_l},(\tau\pi_2)^{-1}J|_{K_l}\bigr)\biggr).
\end{aligned}
\end{equation*}
Note that for any $\sigma\in\mathfrak{S}_k$ the tuples $\sigma^{-1}I|_{K_l}$ and $\sigma^{-1}J|_{K_l}$ depend only on how $\sigma$ acts on $K_l$. Recall also that $\mathfrak{S}(K)=\mathfrak{S}(K_1)\times\dots\times\mathfrak{S}(K_r)$. Using this for any $\sigma,\tau\in\mathfrak{S}_k$, we obtain
\begin{gather*}
\sum_{\pi_1,\pi_2\in\mathfrak{S}(K)}\prod_{l=1}^{r}\bigl(\sgn(\pi_1\pi_2|_{K_l})\Theta_l\bigl((\sigma\pi_1)^{-1}I|_{K_l},(\tau\pi_2)^{-1}J|_{K_l}\bigr)\bigr)
\\
\qquad = \prod_{l=1}^{r}\symdet\bigl(\bigl[x_{i_{\sigma(s)}l}\partial_{j_{\tau(t)}l}\bigr]_{s,t\in K_l}\bigr)=\prod_{l=1}^{r}\Psi_l\bigl(\sigma^{-1}I|_{K_l},\tau^{-1}J|_{K_l}\bigr).
\end{gather*}
Therefore, we get the following expression for $\Delta(I,J)$:
\begin{equation*}
\Delta(I,J)=\sum_{[\sigma],[\tau]\in\mathfrak{S}_k/\mathfrak{S}(K)}\sgn(\sigma\tau)\prod_{l=1}^{r}\Psi_l\bigl(\sigma^{-1}I|_{K_l},\tau^{-1}J|_{K_l}\bigr).
\end{equation*}
\begin{Remark}
The summation here runs over all the left cosets of $\mathfrak{S}(K)$ in $\mathfrak{S}_k$. Note that the expression $\sgn(\sigma\tau)\prod_{l=1}^{r}\Psi_l\bigl(\sigma^{-1}I|_{K_l},\tau^{-1}J|_{K_l}\bigr)$ does not depend on the choice of representatives of $[\sigma],[\tau]\in\mathfrak{S}_k/\mathfrak{S}(K)$.
\end{Remark}
Now we apply Lemma~\ref{sym_det_weyl_lemma} in order to compute for every $l\in R$ the symmetrized determinant $\Psi_l\bigl(\sigma^{-1}I|_{K_l},\tau^{-1}J|_{K_l}\bigr)$. Firstly, observe that $\Psi_l(I_l(\sigma),J_l(\tau))=0$ unless there exist a permutation $\pi_l\in\mathfrak{S}(K_l)$ and an element $s_l\in K_l$ such that $j_{\tau(\pi(s))}=i_{\sigma(s)}$ for all $s\in K_l\setminus\{s_l\}$. Moreover, in~this case we have
\begin{gather*}
\Psi_l\bigl(\sigma^{-1}I|_{K_l},\tau^{-1}J|_{K_l}\bigr)=(-1)^{m_l-1}(m_l-1)!
\\
\qquad\times\sum_{\pi_l\in\mathfrak{S}(K_l)}\sum_{s_l\in K_l}\mathds{1}\bigl(\sigma^{-1} I|_{K_l\setminus\{s_l\}}=(\tau\pi_l)^{-1} J|_{K_l\setminus\{s_l\}}\bigr)\cdot\sgn(\pi_l) x_{(\sigma^{-1}I)(s_l),l}\partial_{(\tau^{-1}J)(\pi_l(s_l)),l}.
\end{gather*}
Combining everything, we get
\begin{gather*}
\Delta(I,J)
=\sum_{[\sigma],[\tau]\in\mathfrak{S}_k/\mathfrak{S}(K)}\sgn(\sigma\tau)\prod_{l=1}^{r}\Psi_l\bigl(\sigma^{-1}I|_{K_l},(\tau\pi_l)^{-1}J|_{K_l}\bigr)
\\
=\biggl(\prod_{l\in R}(-1)^{m_l-1}(m_l-1)!\biggr)\sum_{[\sigma],[\tau]\in\mathfrak{S}_k/\mathfrak{S}(K)}
\sgn(\sigma\tau)
\\
\quad{}\times\sum_{\substack{s_l\in K_l\\\pi_l\in\mathfrak{S}(K_l)\\\text{for each}~l\in R}}\prod_{l\in R}\bigl(\mathds{1}\bigl(\sigma^{-1} I|_{K_l\setminus\{s_l\}}=(\tau\pi_l)^{-1} J|_{K_l\setminus\{s_l\}}\bigr)\sgn(\pi_l) x_{(\sigma^{-1}I)(s_l),l)}\partial_{(\tau^{-1}J)(\pi_l(s_l)),l}\bigr).
\end{gather*}
Now observe that $\mathfrak{S}(K)=\prod_{l\in R}\mathfrak{S}(K_l)$. Hence, when $\tau$ runs over $\mathfrak{S}_k/\mathfrak{S}(K)$ and $\pi_l$ runs over~$\mathfrak{S}(K_l)$ for each $l\in R$, the product $\tau\cdot\prod_{l\in R}\pi_l$ runs over $\mathfrak{S}_k$. Therefore, we can rewrite the formula above using the summation over $\tau\in\mathfrak{S}_k$ as follows:
\begin{gather*}
\Delta(I,J)=\biggl(\prod_{l\in R}(-1)^{m_l-1}(m_l-1)!\biggr)\sum_{[\sigma]\in\mathfrak{S}_k/\mathfrak{S}(K)}\sgn(\sigma)
\\
\quad\times\sum_{\substack{s_l\in K_l\\ \text{for each}~l\in R}}\sum_{\tau\in\mathfrak{S}_k}\mathds{1}\bigl(\sigma^{-1} I|_{[k]\setminus\{s_l\colon l\in R\}}=\tau^{-1} J|_{[k]\setminus\{s_l\colon l\in R\}}\bigr)\sgn(\tau)\prod_{l\in R}\bigl(x_{I(\sigma(s_l)),l}\partial_{J(\tau(s_l)),l}\bigr).
\end{gather*}
Replacing $\tau$ with $\tau\sigma$ and introducing the summation over $\sigma\in\mathfrak{S}_k$ instead of the summation over $[\sigma]\in\mathfrak{S}_k/\mathfrak{S}(K)$ yields
\begin{align*}
&\Delta(I,J)=(-1)^{k-\ell}\biggl(\prod_{l\in R}m_l^{-1}\!\biggr)\!\sum_{\sigma\in\mathfrak{S}_k}\sum_{\substack{s_l\in K_l\\\text{for each}~l\in R}}\sum_{\tau\in\mathfrak{S}_k}\!\mathds{1}\bigl(I|_{[k]\setminus\{\sigma(s_l)\colon l\in R\}}=\tau^{-1} J|_{[k]\setminus\{\sigma(s_l)\colon l\in R\}}\bigr)
\\
&\hphantom{\Delta(I,J)=(-1)^{k-\ell}\biggl(\prod_{l\in R}m_l^{-1}\biggr)\sum_{\sigma\in\mathfrak{S}_k}\sum_{\substack{s_l\in K_l\\\text{for each}~l\in R}}\sum_{\tau\in\mathfrak{S}_k}}
\times\sgn(\tau)\prod_{l\in R}\bigl(x_{I(\sigma(s_l)),l}\partial_{J(\tau(\sigma(s_l))),l}\bigr).
\end{align*}
For a given $\ell$-tuple $\{s_l\}_{l\in R}$, as $\sigma$ runs over $\mathfrak{S}_k$ the tuple $\{\sigma(s_l)\}_{l\in R}$ runs over all $\ell$-tuples in $[k]^{\underline{\ell}}$ and each of them occurs exactly $(k-\ell)!$ times. Since there are precisely $\prod_{l\in R}m_l$ tuples $\{s_l\}_{l\in R}$ in $\prod_{l\in R}R_l$, we have
\begin{align*}
&\Delta(I,J)=(-1)^{k-\ell}(k-\ell)!\sum_{U\in[k]^{\underline{R}}}\sum_{\tau\in\mathfrak{S}_k}\mathds{1}\bigl(I|_{[k]\setminus\{u_l\colon l\in R\}}=\tau^{-1} J|_{[k]\setminus\{u_l\colon l\in R\}}\bigr)
\\
&\hphantom{\Delta(I,J)=(-1)^{k-\ell}(k-\ell)!\sum_{U\in[k]^{\underline{R}}}\sum_{\tau\in\mathfrak{S}_k}}{}
\times\sgn(\tau)\prod_{l\in R}\bigl(x_{I(u_l),l}\partial_{J(\tau(u_l)),l}\bigr).
\end{align*}
Here by $[k]^{\underline{R}}$ we denote the set of all $\ell$-tuples $U=\{u_l\}_{l\in R}$ of distinct elements in $[k]$. The next step is to replace the summation over $\ell$-tuples $U\in[k]^{\underline{R}}$ by the summation over $\ell$-element subsets of $R$. Namely, we collect the terms which correspond to the same $\ell$-element subset of~$R$ and after some algebraic manipulations obtain the determinants in $x_{ij}$ and $\partial_{ij}$. Let $\mathfrak{S}(U)$ be the permutation group of the $\ell$-element set $\{u_l\mid l\in R\}$. Then, we can rewrite the last sum as~follows:
\begin{gather*}
\Delta(I,J)=\frac{(-1)^{k-\ell}(k-\ell)!}{\ell!}\sum_{U\in[k]^{\underline{R}}}\sum_{\pi\in\mathfrak{S}(U)}\sum_{\tau\in\mathfrak{S}_k}\mathds{1}\bigl(I|_{[k]\setminus\{u_l\colon l\in R\}}=\tau^{-1} J|_{[k]\setminus\{u_l\colon l\in R\}}\bigr)
\\
\hphantom{\Delta(I,J)=\frac{(-1)^{k-\ell}(k-\ell)!}{\ell!}\sum_{U\in[k]^{\underline{R}}}\sum_{\pi\in\mathfrak{S}(U)}\sum_{\tau\in\mathfrak{S}_k}}{}
\times\sgn(\tau)\prod_{l\in R}\bigl(x_{I(\pi(u_l)),l}\partial_{J(\tau(\pi(u_l))),l}\bigr).
\end{gather*}
Clearly, for any $\pi\in\mathfrak{S}(U)$, we have
\begin{equation}\label{indicator_IJ}
\mathds{1}\bigl(I|_{[k]\setminus\{u_l\colon l\in R\}}=\tau^{-1} J|_{[k]\setminus\{u_l\colon l\in R\}}\bigr)=\mathds{1}\bigl(I|_{[k]\setminus\{u_l\colon l\in R\}}=(\tau\pi)^{-1} J|_{[k]\setminus\{u_l\colon l\in R\}}\bigr).
\end{equation}
Therefore, substituting $\tau\mapsto\tau\pi^{-1}$, we obtain
\begin{align*}
\Delta(I,J)
={}&\frac{(-1)^{k-\ell}(k-\ell)!}{\ell!}\sum_{U\in[k]^{\underline{R}}}\sum_{\pi\in\mathfrak{S}(U)}\sum_{\tau\in\mathfrak{S}_k}\mathds{1}\bigl(I|_{[k]\setminus\{u_l\colon l\in R\}}=\tau^{-1} J|_{[k]\setminus\{u_l\colon l\in R\}}\bigr)
\\
&\hphantom{\frac{(-1)^{k-\ell}(k-\ell)!}{\ell!}\sum_{U\in[k]^{\underline{R}}}\sum_{\pi\in\mathfrak{S}(U)}\sum_{\tau\in\mathfrak{S}_k}}{}
\times\sgn(\tau\pi)\prod_{l\in R}x_{I(\pi(u_l)),l}\partial_{J(\tau(u_l)),l}
\\
={}&\frac{(-1)^{k-\ell}(k-\ell)!}{\ell!}\sum_{U\in[k]^{\underline{R}}}\sum_{\tau\in\mathfrak{S}_k}\mathds{1}\bigl(I|_{[k]\setminus\{u_l\colon l\in R\}}=\tau^{-1} J|_{[k]\setminus\{u_l\colon l\in R\}}\bigr)\sgn(\tau)
\\
&\hphantom{\frac{(-1)^{k-\ell}(k-\ell)!}{\ell!}\sum_{U\in[k]^{\underline{R}}}\sum_{\tau\in\mathfrak{S}_k}}{}
\times\det\bigl(\bigl[x_{I(u_\alpha),\beta}\bigr]_{\alpha,\beta\in R}\bigr)\prod_{l\in R}\partial_{J(\tau(u_l)),l}.
\end{align*}
Similarly, we have
\begin{gather*}
\sum_{U\in[k]^{\underline{R}}}\det\bigl(\bigl[x_{I(u_\alpha),\beta}\bigr]_{\alpha,\beta\in R}\bigr)\prod_{l\in R}\partial_{J(\tau(u_l)),l}\\
\qquad{}=\sum_{V\in\binom{[k]}{\ell}}\det\bigl(\bigl[x_{I(v),l}\bigr]_{v\in V,l\in R}\bigr)\cdot\det\bigl(\bigl[\partial_{J(\tau(v)),l}\bigr]_{v\in V,l\in R}\bigr).
\end{gather*}
Plugging this into the previous formula and taking into account \eqref{indicator_IJ} gives
\begin{align}
\Delta(I,J)={}&
\frac{(-1)^{k-\ell}(k-\ell)!}{\ell!}\sum_{V\in\binom{[k]}{\ell}}\sum_{\tau\in\mathfrak{S}_k}\mathds{1}\bigl(I|_{[k]\setminus V}=J\circ\tau|_{[k]\setminus V}\bigr)
\nonumber\\
&{}{\times}\,\sgn(\tau)\det\bigl(\bigl[x_{I(v),l}\bigr]_{v\in V,l\in R}\bigr)\cdot\det\bigl(\bigl[\partial_{J(\tau(v)),l}\bigr]_{v\in V,l\in R}\bigr).\label{fin1_expr_sym_det}
\end{align}
Now observe that for any given subset $V\subset[k]$ with $|V|=\ell$ there exist either $\ell!$ permutations $\tau\in\mathfrak{S}_k$ such that
\begin{equation}\label{tau_cond}
I|_{[k]\setminus V}=J\circ\tau|_{[k]\setminus V},
\end{equation}
or none. Indeed, suppose such a permutation $\tau_0$ exists. Then, any other permutation $\tau$ that satisfies \eqref{tau_cond} is of the form $\tau=\tau_0\pi$, where $\pi\in\mathfrak{S}(V)$. In particular,
\begin{equation*}
\begin{split}
\sgn(\tau)\det\bigl(\bigl[\partial_{J(\tau(v)),l}\bigr]_{v\in V,l\in R}\bigr)
&=\sgn(\tau_0\pi)\sgn(\pi)\det\bigl(\bigl[\partial_{J(\tau_0(v)),l}\bigr]_{v\in V,l\in R}\bigr)
\\
&=\sgn(\tau_0)\det\bigl(\bigl[\partial_{J(\tau_0(v)),l}\bigr]_{v\in V,l\in R}\bigr).
\end{split}
\end{equation*}
Hence, the summation over $\tau$ in \eqref{fin1_expr_sym_det} in fact contains $\ell!$ equal summands.
Note that $\tau$ satisfying~\eqref{tau_cond} exists if and only if all elements of $I|_{[k]\setminus V}$ are contained in $J$. This allows us to write the final formula for $\Delta(I,J)$:
\begin{equation}\label{fin2_expr_sym_det}
\Delta(I,J)=(-1)^{k-\ell}(k-\ell)!\sum_{V,W\in\binom{[k]}{\ell}}\varepsilon(I,J,I(V),J(W))\det(X_{I(V),R})\cdot\det(D_{J(W),R}).
\end{equation}
Indeed, to get \eqref{fin2_expr_sym_det} from~\eqref{fin1_expr_sym_det}, we just note that \smash{$\mathds{1}\bigl(I|_{[k]\setminus V}=J\circ\tau|_{[k]\setminus V}\bigr)\sgn(\tau)$} is equal to $\varepsilon(I,J,I(V),J(W))$ if $W=\tau(V)$ and $\tau$ satisfies~\eqref{tau_cond}, and is zero otherwise.
\end{proof}

\begin{Lemma}\label{main_sym_det_lem}
For any $I,J\in[n]^{\underline{k}}$, we have
\begin{gather*}
\symdet L(E)_{IJ}
\\
=\sum_{\ell=0}^{k}(-1)^{k-\ell}(k-\ell)!\ell!\stirling{k}{\ell}\! \sum_{R\in\binom{[r]}{\ell}}\sum_{V,W\in\binom{[k]}{\ell}}\varepsilon(I,J,I(V),J(W))\det\bigl(X_{I(V),R}\bigr)\det\bigl(D_{J(W),R}\bigr).
\end{gather*}
Here by \smash{$\stirling{k}{\ell}$} we denote the Stirling number of the second kind, i.e., the number of ways to split a $k$-element set into $\ell$ non-empty subsets.
\end{Lemma}
\begin{proof}
We start with the identity
\begin{equation*}
\symdet(L(E))_{IJ}=\sum_{a_1,\dots,a_k=1}^{n}\sum_{\sigma,\tau\in\mathfrak{S}_{k}}\sgn(\sigma\tau)\cdot x_{i_{\sigma(1)}a_1}\partial_{j_{\tau(1)}a_1}\cdots x_{i_{\sigma(k)}a_k}\partial_{j_{\tau(k)}a_k}.
\end{equation*}
For each choice of $a_1,\dots,a_k\in\{1,\dots,r\}$ we rewrite the products using Lemma~\ref{lemma:delta_a_IJ}. It remains to notice that for any given \smash{$R\in\binom{r}{\ell}$} the number of sequences $(a_1,\dots,a_k)\in[r]^{k}$ such that $\{a_1,\dots,a_k\}=R$ equals \smash{$\stirling{k}{\ell}\cdot \ell!$}
\end{proof}

\subsection{Proof of the main formulas} Now we are ready to prove Theorem~\ref{big_generator_thm} and Corollary~\ref{cor:F_operators}.
\begin{proof}[Proof of Theorem~\ref{big_generator_thm}]
By Proposition~\ref{symmetrized_m_pq_prop}, we have
\begin{equation*}
M_{p,q}(Y)=\sum_{\substack{I_1,J_1\in\binom{[n]}{p}\\I_2,J_2\in\binom{[n]}{q}\\I_1\sqcup I_2=J_1\sqcup J_2}}\sgn\binom{I_1~I_2}{J_1~J_2}\det Y_{I_1J_1}\cdot\symdet L(E)_{J_2I_2}.
\end{equation*}
Applying Lemma~\ref{main_sym_det_lem} to $\symdet L(E)_{J_2I_2}$ gives the required identity \eqref{big_generator_formula}.
\end{proof}

\begin{proof}[Proof of Corollary~\ref{cor:F_operators}]
Indeed, the identity \eqref{m_f_operators_transition} is a direct consequence of \eqref{big_generator_formula} and the fact that the expression \smash{$\sgn\binom{I_1~I_2}{J_1~J_2}\varepsilon(J_2,I_2,J_2(V),I_2(W))$}, if non-zero, equals \smash{$\sgn\binom{I_1~I_2(W)}{J_1~J_2(V)}$}. The~coefficient \smash{$\binom{n-p-\ell}{q-\ell}$} appears as the number of ways to choose a $(q-\ell)$-element subset $I_2\setminus I_2(W)=J_2\setminus J_2(V)$ in the complement of the $(p+\ell)$-element subset $I_1\sqcup I_2(W)=J_1\sqcup J_2(V)$ of~$[n]$.

Then, it follows from \eqref{m_f_operators_transition} that the difference $M_{p,q}-q!\cdot F_{p,q}$ is a linear combination of~$F_{p,\ell}$ with~$\ell<q$. Therefore, the sets $\{M_{p,q}\}_{p+q\le n}$ and $\{F_{p,q}\}_{p+q\le n}$ generate the same algebra inside~${\mathscr{C}(\mathcal{P}(n,r))}$.
\end{proof}

\section{Capelli identities and their variations}\label{sect:capelli_identities}

In this section, we give an account of Capelli identities and their variants which we use in the proof of Theorem~\ref{thm:big_alg_commut}. Since some of these generalizations are not widely known, we also include the proofs for the sake of completeness.

\subsection{Classical Capelli identity} Assume that $r=n$, i.e., that $\mathcal{P}(n,n)$ is a polynomial ring $\mathbb{C}[\Mat(n,n)]$ in $n^2$ variables $x_{ij}$, $1\le i,j\le n$. Recall that $\mathcal{PD}(n,n)$ is the algebra of differential operators on $\mathcal{P}(n,n)$ with polynomial coefficients. Define an element $\Pi\in U(\mathfrak{gl}_n)$ as
\begin{equation*}
\Pi=\rdet(E_{ij}+(i-1)\delta_{ij})_{i,j=1}^{n}=\rdet
\begin{bmatrix}
E_{11}+0 & E_{12} & \cdots & E_{1k}
\\
E_{21} & E_{22}+1 & \cdots & E_{2k}
\\
\vdots & \vdots & \ddots & \vdots
\\
E_{k1} & E_{k2} & \cdots & E_{kk}+k-1
\end{bmatrix}.
\end{equation*}

The next statement was observed by Alfredo Capelli \cite{Capelli} and was used by Hermann Weyl in his treatment of invariant theory for $\mathrm{GL}_n$ \cite[Section~II.4]{Weyl}.

\begin{Proposition}[{Capelli identity, \cite{Capelli}}]\label{capelli_prop} The image of $\Pi$ in $\mathcal{PD}(n,n)$ equals
\begin{equation*}
L(\Pi)=\det(X)\det(D),
\end{equation*}
where $L$ is defined as in \eqref{eq:L_representation}.
In particular, the expression on the left is a $\mathrm{GL}_n$-invariant differential operator.
\end{Proposition}

\subsection{Capelli identities for rectangular matrices} Now let $r$ be an arbitrary positive integer.
It turns out that one can generalize the Capelli identity for all minors of the matrix $E$ (see \cite{CSS, Howe_Umeda, Umeda} and also \cite{Howe_89,Okounkov}).

Firstly, for arbitrary $k$-tuples $I=(i_1,\dots,i_k)$ and $J=(j_1,\dots,j_k)$ in $[n]^{k}$, introduce the following element $\Pi_{IJ}\in U(\mathfrak{gl}_n)$:
\begin{equation*}
\Pi_{IJ}=\rdet [E_{i_{\alpha}j_{\beta}}+(\alpha-1)\delta_{i_{\alpha}j_{\beta}}]_{\alpha,\beta=1}^{k}.
\end{equation*}
For example, if $I=J=(i_1,\dots,i_k)\in[n]^{\underline{k}}$, then
\begin{equation*}
\Pi_{II}=\rdet
\begin{bmatrix}
E_{i_1i_1}+0 & E_{i_1i_2} & \cdots & E_{i_1i_k}
\\
E_{i_2i_1} & E_{i_2i_2}+1 & \cdots & E_{i_2i_k}
\\
\vdots & \vdots & \ddots & \vdots
\\
E_{i_ki_1} & E_{i_ki_2} & \cdots & E_{i_ki_k}+k-1
\end{bmatrix},
\end{equation*}
which coincides with $\Pi$ above when $I=(1,2,\dots,n)$.

\begin{Proposition}[{\cite[Proposition 4 and Theorem~5]{Umeda}}]\label{lower_capelli_prop}
The following identity holds in $U(\mathfrak{gl}_n)[z]$:%
\begin{equation}\label{capelli_det_z_exp}
\rdet(E_{ij}+(i-1-z)\delta_{ij})_{i,j=1}^{n}=\sum_{k=0}^{n}(-1)^kz^{\underline{k}}\cdot C_{n-k},
\end{equation}
where $z^{\underline{k}}=z(z-1)\cdots (z-k+1)$ and
$C_k=\sum_{I\in\binom{[n]}{k}}\Pi_{II}$.
Moreover, the elements $C_k$ belong to the center $Z(\mathfrak{gl}_n)$ of the universal enveloping algebra $U(\mathfrak{gl}_n)$ and their images in $\mathcal{PD}(n,r)$ are as follows:
\begin{equation}\label{lower_capelli_ident}
L(C_k)=\sum_{I\in\binom{[n]}{k}}\sum_{K\in\binom{[r]}{k}}\det(X_{IK})\det(D_{IK})
\end{equation}
In particular, if $r=n$, then plugging in $z=0$ yields the classical Capelli identity $($Proposi\-tion~{\rm \ref{capelli_prop})}.
\end{Proposition}
\begin{Remark}
In formula \eqref{lower_capelli_ident}, we regard \smash{$K\in\binom{[r]}{k}$} as an element of \smash{$\mathfrak{S}_k\backslash[r]^{\underline{k}}$}, i.e., as $k$-tuple with arbitrarily chosen ordering of elements. Note that the term $\det(X_{IK})\det(D_{IK})$ does not depend on this choice, so the expression on the right is well defined.
\end{Remark}
\begin{Remark}
The elements $C_k$ are often called \emph{Capelli generators}. Moreover, $Z(\mathfrak{gl}_n)$ is a free polynomial ring generated by $C_1,\dots,C_n$. The classical counterparts of the Capelli generators are the coefficients $c_k$ (see \eqref{eq:c_k_formula}) of the characteristic polynomial which are the generators of~$S(\mathfrak{gl}_n^*)^{\mathfrak{gl}_n}$.
\end{Remark}

Instead of proving Proposition~\ref{lower_capelli_prop} directly, we prove more general versions of equalities \eqref{capelli_det_z_exp} and \eqref{lower_capelli_ident}, see Propositions \ref{capelli_minor_z_exp_prop} and \ref{prop:cauchy_binet_capelli} below. These generalizations are then verified by a~straightforward inductive argument.

\subsection{Cauchy--Binet type identity}
We start with a formula for $L(\Pi_{IJ})$ for an arbitrary minor $\Pi_{IJ}$. It generalizes the classical Capelli identity and can be viewed as a non-commutative analogue of the Cauchy--Binet identity from linear algebra (see also \cite[Corollary 1.3\,(b)]{CSS}).

\begin{Proposition}[{Cauchy--Binet type identity, \cite[Theorem~2]{Umeda}}]\label{prop:cauchy_binet_capelli}
For any $I,J\in[n]^{k}$, the image of $\Pi_{IJ}$ in $\mathcal{PD}(n,r)$ under the map $L$ equals
\begin{equation}\label{capelli_minor_ident}
L(\Pi_{IJ})=\sum_{K\in\binom{[r]}{k}}\det(X_{IK})\det(D_{JK}).
\end{equation}
\end{Proposition}
\begin{Remark}\label{remark:skew_sym_capelli_minor}
Note that the right-hand side of \eqref{capelli_minor_ident} is skew-symmetric in the entries of $I$ and~$J$ and hence the same holds for $L(\Pi_{IJ})$. From the definition of $\Pi_{IJ}$, it is clear that $\Pi_{IJ}$ is skew-symmetric in $J$ (as row determinant), but the skew-symmetry in $I$ is not immediate from this definition. Formula \eqref{capelli_minor_ident} in particular implies that $L(\Pi_{IJ})$ is zero whenever $I$ or $J$ contains equal elements.
\end{Remark}
\begin{proof}[Proof of Proposition~\ref{prop:cauchy_binet_capelli}]
Induct on $k$. For $k=1$, the identity follows immediately from the definition of $L(E_{ij})$, see \eqref{eq:L_representation}. Now assume that $k>1$. Consider arbitrary $k$-tuples $I=(i_1,\dots,i_k)$ and $J=(j_1,\dots,j_k)$. Denote $I'=(i_1,\dots,i_{k-1})$ and \smash{$J^{(l)}=\bigl(j_1,\dots,\widehat{j_{l}},\dots,j_{k}\bigr)$} for every $l\in[k]$. Expanding $\Pi_{IJ}$ along the $k$-th row yields
\begin{equation*}
\Pi_{IJ}=\sum_{l=1}^{k}(-1)^{k-l}\Pi_{I'J^{(l)}}\cdot (E_{i_k j_{l}}+(k-1)\delta_{i_k j_{l}}),
\end{equation*}
and hence by the inductive hypothesis
\begin{equation*}
\begin{aligned}
L(\Pi_{IJ})
&=\sum_{l=1}^{k}(-1)^{k-l}L\big(\Pi_{I'J^{(l)}}\big)\cdot \biggl(\sum_{\alpha\in[r]}x_{i_k\alpha}\partial_{j_{l}\alpha}+(k-1)\delta_{i_k j_{l}}\biggr)
\\
&=\sum_{l=1}^{k}(-1)^{k-l}\sum_{K'\in\binom{[r]}{k-1}}\det(X_{I'K'})\det\bigl(D_{J^{(l)}K'}\bigr)\cdot \biggl(\sum_{\alpha\in[r]}x_{i_k\alpha}\partial_{j_{l}\alpha}+(k-1)\delta_{i_k j_{l}}\biggr).
\end{aligned}
\end{equation*}
Using the identity \smash{$\det\bigl(D_{J^{(l)}K'}\bigr)x_{i_k\alpha}=x_{i_k\alpha}\det\bigl(D_{J^{(l)}K'}\bigr) +\bigl[\det\bigl(D_{J^{(l)}K'}\bigr),x_{i_k\alpha}\bigr]$}, we can rewrite the last sum as follows:
\begin{align}
L(\Pi_{IJ})
={}&\sum_{l=1}^{k}\sum_{\alpha\in[r]}\sum_{K'\in\binom{[r]}{k-1}}(-1)^{k-l}\det(X_{I'K'})x_{i_k\alpha}\cdot\det\bigl(D_{J^{(l)}K'}\bigr)\partial_{j_{l}\alpha} \nonumber
\\
&{}{+}\,(k-1)\sum_{l=1}^{k}\sum_{K'\in\binom{[r]}{k-1}}(-1)^{k-l}\det(X_{I'K'})\cdot\det\bigl(D_{J^{(l)}K'}\bigr)\delta_{i_k j_{l}} \nonumber
\\
&{}{+}\,\sum_{l=1}^{k}\sum_{\alpha\in[r]}\sum_{K'\in\binom{[r]}{k-1}}(-1)^{k-l}\det(X_{I'K'})\cdot\bigl[\det\bigl(D_{J^{(l)}K'}\bigr),x_{i_k\alpha}\bigr]\partial_{j_{l}\alpha}.\label{L_Pi_IJ}
\end{align}
Denote the three summands on the right-hand side of \eqref{L_Pi_IJ} by \smash{$S^{(1)}_{IJ}$}, \smash{$S^{(2)}_{IJ}$} and \smash{$S^{(3)}_{IJ}$}, respectively.
Observe that \smash{$S^{(1)}_{IJ}$} equals
\begin{align*}
S^{(1)}_{IJ}&{}=\sum_{l=1}^{k}\sum_{\alpha\in[r]}\sum_{K'\in\binom{[r]}{k-1}}(-1)^{k-l}\det(X_{I'K'})x_{i_k\alpha}\cdot\det\bigl(D_{J^{(l)}K'}\bigr)\partial_{j_{l}\alpha}
\\
&{}=\sum_{\alpha\in[r]}\sum_{K'\in\binom{[r]}{k-1}}\det(X_{I'K'})x_{i_k\alpha}\cdot\det\bigl(D_{JK'_{(\alpha)}}\bigr)=\sum_{K\in\binom{[r]}{k}}\det(X_{IK})\cdot\det(D_{JK}).
\end{align*}
Here by \smash{$K'_{(\alpha)}$} we denote a $k$-tuple in $[n]^{\underline{k}}$ whose first $k-1$ entries coincide with those of $K'$ and whose $k$-th entry equals $\alpha$. (Recall that $K'$ can be viewed as a $(k-1)$-tuple.) The last equality follows from the cofactor expansion of $\det(X_{IK})$ along the last row. Indeed, for $K=(\beta_1,\dots,\beta_k)$ we have \smash{$\det(X_{IK})=\sum_{l=1}^{l}(-1)^{k-l}\det\bigl(X_{I'K^{(l)}}\bigr)x_{i_k\beta_l}$}, where \smash{$K^{(l)}=\bigl(\beta_1,\dots,\widehat{\beta_l},\dots,\beta_k\bigr)\in[r]^{\underline{k-1}}$}. If we denote \smash{$\widetilde{K}_{(l)}=(\beta_1,\dots,\beta_{l-1},\beta_{l+1},\dots,\beta_{k},\beta_{l})\in[r]^{\underline{k}}$} for $l=1,\dots,k$, then we have the equality \smash{$\det(D_{JK})=(-1)^{k-l}\det\bigl(D_{J\widetilde{K}_{(l)}}\bigr)$} and the formula for \smash{$S^{(1)}_{IJ}$} above follows.

Thus, it remains to show that the sum \smash{$S^{(2)}_{IJ}+S^{(3)}_{IJ}$} is zero. We consider three cases depending on the number $q$ of occurrences of $i_k$ in the $k$-tuple $J$.

\textit{Case 1:} $q=0$. Then, for any $l\in[k]$ and $\alpha\in[r]$, we have $\delta_{i_kj_{l}}=0$ and ${\bigl[\det\bigl(D_{J^{(l)}K'}\bigr),x_{i_k\alpha}\bigr]\!=\!0}$. Hence, both \smash{$S^{(2)}_{IJ}$} and \smash{$S^{(3)}_{IJ}$} are zero which concludes the proof.

In the remaining two cases, we use the following lemma.
\begin{Lemma}\label{commut_d_det}
For any $k$-tuples $I=(i_1,\dots,i_k)$ and $J=(j_1,\dots,j_k)$ in $[n]^{k}$ and any $i,j\in[n]$, we have
\begin{equation*}
\sum_{\alpha\in[r]}[\det(D_{IJ}),x_{i\alpha}]\partial_{j\alpha}=
\begin{cases}
0 &\text{if}~I\notin[n]^{\underline{k}}~\text{or}~i\notin I,
\\
\det(D_{KJ}) &\text{if}~i_{p}=i~\text{for}~p\in[k],
\end{cases}
\end{equation*}
where in the second case we put $K=(i_1,\dots,i_{p-1},j,i_{p+1},\dots,i_{k})$.
\end{Lemma}
\begin{proof}
If $I$ contains equal entries or $i\notin I$, then it clear that $[\det(D_{IJ}),x_{i\alpha}]=0$ for all $\alpha\in[r]$. Otherwise, if $p\in[k]$ is such that $i_p=i$, then $[\det(D_{IJ}),x_{i\alpha}]$ is zero if $\alpha\notin J$ and equals to the cofactor of the element $\partial_{j\alpha}$ of the matrix $D_{KJ}$ otherwise. Thus, $\sum_{\alpha\in[r]}[\det(D_{IJ}),x_{i\alpha}]\partial_{j\alpha}=\det(D_{KJ})$, as claimed.
\end{proof}

\textit{Case 2:} $q=1$. Let $p\in[k]$ be such that $i_k=j_{p}$. Observe that both sides of \eqref{capelli_minor_ident} are skew-symmetric with respect to~$J$ (see also Remark~\ref{remark:skew_sym_capelli_minor}). Thus, we may assume without loss of generality that $p=k$. In this case, we can rewrite \smash{$S^{(2)}_{IJ}$} as follows:
\begin{align*}
\begin{split}
S^{(2)}_{IJ}
&=(k-1)\sum_{l=1}^{k}\sum_{K'\in\binom{[r]}{k-1}}(-1)^{k-l}\det(X_{I'K'})\cdot\det\bigl(D_{J^{(l)}K'}\bigr)\delta_{i_k j_{l}}
\\
&=(k-1)\sum_{K'\in\binom{[r]}{k-1}}\det(X_{I'K'})\cdot\det\bigl(D_{J^{(k)}K'}\bigr).
\end{split}
\end{align*}
Applying Lemma~\ref{commut_d_det} gives the following expression for \smash{$S^{(3)}_{IJ}$}:
\begin{equation*}
\begin{split}
S^{(3)}_{IJ}
&=\sum_{l=1}^{k}\sum_{\alpha\in[r]}\sum_{K'\in\binom{[r]}{k-1}}(-1)^{k-l}\det(X_{I'K'})\cdot\bigl[\det\bigl(D_{J^{(l)}K'}\bigr),x_{i_k\alpha}\bigr]\partial_{j_{l}\alpha}
\\
&=\sum_{l\in[k-1]}\sum_{K'\in\binom{[r]}{k-1}}(-1)^{k-l}\det(X_{I'K'})\cdot(-1)^{k-l-1}\det\bigl(D_{J^{(k)}K'}\bigr)
\\
&=-(k-1)\sum_{K'\in\binom{[r]}{k-1}}\det(X_{I'K'})\cdot\det\bigl(D_{J^{(k)}K'}\bigr).
\end{split}
\end{equation*}
Therefore, \smash{$S^{(2)}_{IJ}$} and \smash{$S^{(3)}_{IJ}$} add up to zero, as claimed.

\textit{Case 3:} $q>1$. Let us show that in this case both \smash{$S^{(2)}_{IJ}$} and \smash{$S^{(3)}_{IJ}$} vanish. If $q\geq 3$, then each $(k-1)$-tuple $J^{(l)}$, $l=1,\dots,k$, contains at least two entries that are equal to $i_k$. Hence, all terms in sums \smash{$S^{(2)}_{IJ}$} and \smash{$S^{(3)}_{IJ}$} are zero since $\det(D_{J^{l}K'})=0$.

Now assume that $q=2$. Similar to the second case, we may assume that $j_{k-1}=j_{k}=i_k$ and $j_l\neq i_k$ for $l<k-1$. Then, for \smash{$S^{(2)}_{IJ}$}, we have
\begin{equation*}
\begin{split}
S^{(2)}_{IJ}
&=(k-1)\sum_{l=1}^{k}\sum_{K'\in\binom{[r]}{k-1}}(-1)^{k-l}\det(X_{I'K'})\cdot\det\bigl(D_{J^{(l)}K'}\bigr)\delta_{i_k j_{l}}
\\
&=(k-1)\sum_{K'\in\binom{[r]}{k-1}}\bigl(-\det(X_{I'K'})\cdot\det\bigl(D_{J^{(k-1)}K'}\bigr)+\det(X_{I'K'})\cdot\det\bigl(D_{J^{(k)}K'}\bigr)\bigr)=0,
\end{split}
\end{equation*}
because \smash{$\det\bigl(D_{J^{(l)}K'}\bigr)=0$} for $l<k-1$ while $J^{(k-1)}=J^{(k)}$.
Lemma~\ref{commut_d_det} implies that \smash{$S^{(3)}_{IJ}$} equals
\begin{align*}
S^{(3)}_{IJ}
={}&\sum_{l=1}^{k}\sum_{\alpha\in[r]}\sum_{K'\in\binom{[r]}{k-1}}(-1)^{k-l}\det(X_{I'K'})\cdot\bigl[\det\bigl(D_{J^{(l)}K'}\bigr),x_{i_k\alpha}\bigr]\partial_{j_{l}\alpha}
\\
={}&\sum_{K'\in\binom{[r]}{k-1}}\det(X_{I'K'})\\
&{}{\times}\,\sum_{\alpha\in[r]}\bigl(-\bigl[\det\bigl(D_{J^{(k-1)}K'}\bigr),x_{i_k\alpha}\bigr]\partial_{j_{k-1}\alpha} +\bigl[\det\bigl(D_{J^{(k)}K'}\bigr),x_{i_k\alpha}\bigr]\partial_{j_{k}\alpha}\bigr)=0,
\end{align*}
since $J^{(k-1)}=J^{(k)}$ and $j_{k-1}=j_{k}=i_k$.

Therefore, if $q>1$, we have \smash{$S^{(2)}_{IJ}=S^{(3)}_{IJ}=0$} which concludes the proof of the proposition.
\end{proof}

\begin{Corollary}
For any $I,J\in[n]^{k}$ and any $\sigma,\tau\in\mathfrak{S}_k$, we have $\Pi_{\sigma I,\tau J}=\sgn(\sigma\tau)\cdot \Pi_{IJ}$.
\end{Corollary}
\begin{proof}
By the previous proposition, we have the identity $L(\Pi_{\sigma I,\tau J})=\sgn(\sigma\tau)\cdot L(\Pi_{IJ})$. By~Pro\-position~\ref{prop:L_nn_injectivity}, for $r=n$ the map $L$ faithfully maps $U(\mathfrak{gl}_n)$ to $\mathcal{PD}(n,r)$, hence the result follows.
\end{proof}

\begin{Corollary}\label{cdet_rdet_cor}
For any $k$-tuples $I=(i_1,\dots,i_k)$ and $J=(j_1,\dots,j_k)$ in $[n]^{k}$, we have
\begin{equation*}
\Pi_{IJ}=\rdet\bigl[E_{i_{\alpha}j_{\beta}}+(\alpha-1)\delta_{i_{\alpha}j_{\beta}}\bigr]_{\alpha,\beta=1}^{k} =\cdet\bigl[E_{i_{\alpha}j_{\beta}}+(k-\beta)\delta_{i_{\alpha}j_{\beta}}\bigr]_{\alpha,\beta=1}^{k}.
\end{equation*}
\end{Corollary}
\begin{proof}
We prove that the images of $\Pi_{IJ}$ and \smash{$\cdet\bigl[E_{i_{\alpha}j_{\beta}}+(k-\beta)\delta_{i_{\alpha}j_{\beta}}\bigr]_{\alpha,\beta=1}^{k}$} under the map $L$ coincide. This is sufficient since by Proposition~\ref{prop:L_nn_injectivity}, for $r=n$ the map $L\colon U(\mathfrak{gl}_n)\to\mathcal{PD}(n,r)$ is injective.

In view of formula \eqref{capelli_minor_ident}, the equality \smash{$L(\Pi_{IJ})=\cdet\bigl[L\bigl(E_{i_{\alpha}j_{\beta}}\bigr)+(k-\beta)\delta_{i_{\alpha}j_{\beta}}\bigr]_{\alpha,\beta=1}^{k}$} is equivalent~to
\begin{equation*}
\cdet\bigl[L\bigl(E_{i_{\alpha}j_{\beta}}\bigr)+(k-\beta)\delta_{i_{\alpha}j_{\beta}}\bigr]_{\alpha,\beta=1}^{k}=\sum_{K\in\binom{[r]}{k}}\det(X_{IK})\det(D_{JK}).
\end{equation*}
This identity can be proved in essentially the same way as Proposition~\ref{prop:cauchy_binet_capelli}.
\end{proof}

\begin{Remark}
One can prove this identity directly in $U(\mathfrak{gl}_n)$ using the ternary relations and the Yang--Baxter equation (see Lemma~\ref{capelli_minor_char_poly_lem} in the next section).
\end{Remark}

\subsection{Characteristic polynomial of the Capelli matrix} Now we state a generalization of the formula \eqref{capelli_det_z_exp} for $k\times k$ minors.
Define (see also \cite[equation~(11.1.11)]{Howe_Umeda})
\begin{equation*}
C_k(z)=\sum_{I\in\binom{[n]}{k}}\rdet\bigl(E_{i_{\alpha}i_{\beta}}+(\alpha-1-z)\delta_{ij}\bigr)_{\alpha,\beta\in [k]}.
\end{equation*}
Observe that $C_k(0)=C_k$ and $C_n(z)$ is the left-hand side of \eqref{capelli_det_z_exp}.

\begin{Proposition}[{\cite[equation~(11.1.13)]{Howe_Umeda} and \cite[Section~2]{Umeda}}]\label{capelli_minor_z_exp_prop}
For any $0\le k\le n$, we have
\begin{equation}\label{capelli_minor_z_exp}
C_k(z)=\sum_{m=0}^{k}(-1)^m\binom{n-k+m}{m}z^{\underline{m}}\cdot C_{k-m}.
\end{equation}
In particular, for $k=n$ this identity is equivalent to \eqref{capelli_det_z_exp}.
\end{Proposition}

\begin{proof} We follow the proof from \cite[Section~11.1]{Howe_Umeda}.
Induct on $k$. The cases $k=0$ and $k=1$ are clear. Denote the right hand-side of \eqref{capelli_minor_z_exp} by $B_k(z)$. Let $\Delta$ be the difference operator defined as~${(\Delta f)(z)=f(z+1)-f(z)}$. Since both $B_k(z)$ and $C_k(z)$ are $U(\mathfrak{gl}_n)$-valued polynomials in $z$ and~${B_k(0)=C_k(0)=C_k}$, it suffices to check that $(\Delta B_k)(z)=(\Delta C_k)(z)$. We have
\begin{align*}
(\Delta B_k)(z)&{}=\sum_{m=0}^{k}(-1)^m\binom{n-k+m}{m}\Delta(z^{\underline{m}})\cdot C_{k-m}\\
&{}=\sum_{m=1}^{k}(-1)^m m\binom{n-k+m}{m}z^{\underline{m-1}}\cdot C_{k-m}.
\end{align*}
Note that $m\binom{n-k+m}{m}=(n-k+1)\binom{n-k+m}{m-1}$, and hence we obtain
\begin{align*}
(\Delta B_k)(z)&{}=(n-k+1)\sum_{m=0}^{k-1}(-1)^{m+1}\binom{n-k+m+1}{m}z^{\underline{m}}\cdot C_{k-1-m}\\
&{}=-(n-k+1)B_{k-1}(z).
\end{align*}
Now let us compute $(\Delta C_k)(z)$. For any $k$-element subset $I=\{i_1,\dots,i_k\}$ of $[n]$, where $i_1<\dots<i_k$, we set
\begin{gather*}
E_{II}^{\natural}
=
\begin{bmatrix}
E_{i_1i_1}+0 & E_{i_1i_2} & \cdots & E_{i_1i_k}
\\
E_{i_2i_1} & E_{i_2i_2}+1 & \cdots & E_{i_2i_k}
\\
\vdots & \vdots & \ddots & \vdots
\\
E_{i_ki_1} & E_{i_ki_2} & \cdots & E_{i_ki_k}+k-1
\end{bmatrix}
\quad \text{and} \quad
\Pi_{II}(z)=\rdet\bigl(E_{II}^{\natural}-z\cdot \Id_k\bigr).
\end{gather*}
Then,
\begin{equation*}
\begin{split}
(\Delta\Pi_{II})(z)
&=\Pi_{II}(z+1)-\Pi_{II}(z)
\\
&=\sum_{p=1}^{k}\bigl(\rdet\bigl(E_{II}^{\natural}-(z+1)\cdot\Id_k+\Lambda_{p-1}\bigr)-\rdet\bigl(E_{II}^{\natural}-(z+1)\cdot\Id_k+\Lambda_{p}\bigr)\bigr),
\end{split}
\end{equation*}
where $\Lambda_p$ is a diagonal matrix with the first $p$ diagonal entries equal to 1 and the last $k-p$ entries equal to~0. Recall that the row determinant is an additive function of a fixed row. Applying this to the $p$-th row of an $n\times n$ matrix $E^{\natural}-(z+1)\cdot\Id_k+\Lambda_{p-1}$ gives
\begin{gather*}
\rdet\bigl(E_{II}^{\natural}-(z+1)\cdot\Id_k+\Lambda_{p-1}\bigr)\\
\qquad{}=\rdet\bigl(E_{II}^{\natural}-(z+1)\cdot\Id_k+\Lambda_{p}\bigr)-\rdet\bigl(E_{I^{(p)}I^{(p)}}^{\natural}-z\cdot\Id_{k-1}\bigr),
\end{gather*}
where $I^{(p)}=I\setminus\{i_p\}$. Therefore,
\begin{equation*}
(\Delta C_k)(z)=\sum_{I\in\binom{[n]}{k}}(\Delta\Pi_{II})(z)=-\sum_{I\in\binom{[n]}{k}}\sum_{p=1}^{k}\Pi_{I^{(p)}I^{(p)}}(z)=-(n-k+1)C_{k-1}(z).
\end{equation*}
By the inductive hypothesis, $B_{k-1}(z)=C_{k-1}(z)$. Therefore, $(\Delta B_k)(z)=(\Delta C_k)(z)$ which concludes the inductive step due to the equality $B_k(0)=C_k(0)=C_k$.
\end{proof}

\section[The Yangian Y(gl\_n) and its Bethe subalgebras]{The Yangian $\boldsymbol{\operatorname{Y}(\mathfrak{gl}_{n})}$ and its Bethe subalgebras}\label{s6}

Here we recall the definitions and basic properties of the Yangian $\operatorname{Y}(\mathfrak{gl}_{n})$ and the Bethe subalgebras. Our exposition follows the monograph by Molev (see \cite[Section~1]{Molev_07} for more details).

\subsection{Notation} Let $\{e_{ij}\}_{i,j=1}^{n}$ be the standard matrix units of $\Mat(n,n)$ and let $\{E_{ij}\}_{i,j=1}^{n}$ be the corresponding generators of the universal enveloping algebra $U(\mathfrak{gl}_n)$. The matrix units $\{e_{ij}\}_{i,j=1}^{n}$ act on $\mathbb{C}^n$ spanned by $e_1,\dots,e_n$ in the usual way,
$e_{ij}e_{k}=\delta_{jk}e_{i}$.

We often work with algebras of the form $\mathcal{A}\otimes(\End\mathbb{C}^n)^{\otimes m}$, where $\mathcal{A}$ is an associative algebra. For any $C=\sum_{i,j=1}^{n}c_{ij}\otimes e_{ij}\in \mathcal{A}\otimes\End\mathbb{C}^n$ and any $D=\sum_{i,j,k,l=1}^{n}d_{ijkl}e_{ij}\otimes e_{kl}\in(\End\mathbb{C}^n)^{\otimes 2}$ (these operators might depend on some parameters), we define
\begin{equation*}
\begin{split}
&C_{a}
=\sum_{i,j=1}^{n}c_{ij}\otimes 1^{\otimes (a-1)}\otimes e_{ij}\otimes 1^{\otimes(m-a)}\in \mathcal{A}\otimes(\End\mathbb{C}^n)^{\otimes m},\qquad \text{and}
\\
&D_{ab}
=\sum_{i,j=1}^{n}d_{ijkl}\cdot 1^{\otimes (a-1)}\otimes e_{ij}\otimes 1^{\otimes(b-a-1)}\otimes e_{kl}\otimes 1^{\otimes(m-b)}\in(\End\mathbb{C}^n)^{\otimes m}.
\end{split}
\end{equation*}
Usually, we identify $(\End\mathbb{C}^n)^{\otimes m}$ with the subalgebra $1\otimes(\End\mathbb{C}^n)^{\otimes m}$ inside $\mathcal{A}\otimes(\End\mathbb{C}^n)^{\otimes m}$.

Each element $\sigma$ of the symmetric group $\mathfrak{S}_{m}$ defines an element of $(\End\mathbb{C}^n)^{\otimes m}$ whose action on $(\mathbb{C}^n)^{\otimes m}$ corresponds to permuting the tensors via $\sigma$. Namely, $\sigma\in\mathfrak{S}_{m}$ corresponds to
\begin{equation*}
\sum_{i_1,\dots,i_m=1}^{n}e_{i_1i_{\sigma(1)}}\otimes\dots\otimes e_{i_mi_{\sigma(m)}}\in(\End\mathbb{C}^n)^{\otimes m}.
\end{equation*}
Clearly, this gives rise to an embedding of the group algebra $\mathbb{C}[\mathfrak{S}_m]$ into $(\End\mathbb{C}^n)^{\otimes m}$. For any distinct $i,j\in\{1,\dots,m\}$, denote by $P_{ij}$ the image of the transposition $(i~j)\in\mathbb{C}[\mathfrak{S}_m]$ in~$(\End\mathbb{C}^n)^{\otimes m}$.
Let $A_m$ be the antisymmetrization operator, i.e.,
\begin{equation*}
A_m=\sum_{\sigma\in\mathfrak{S}_m}\sgn(\sigma)\sum_{i_1,\dots,i_m=1}^{n}e_{i_1i_{\sigma(1)}}\otimes\dots\otimes e_{i_mi_{\sigma(m)}}\in(\End\mathbb{C}^n)^{\otimes m}.
\end{equation*}
It is not difficult to see that $A_m^2=m!\cdot A_m$ and $A_mP_{ij}=P_{ij}A_m=-A_m$ for any distinct $i,j\in\{1,\dots,m\}$.

\subsection[R-matrices]{$\boldsymbol{R}$-matrices}

Define the \emph{Yang $R$-matrix} as an element \smash{$R(u)\in\operatorname{Y}(\mathfrak{gl}_n)\bigl[\hspace{-0.3mm}\bigl[u^{-1} \bigr]\hspace{-0.3mm}\bigr]\otimes\End(\mathbb{C}^n)^{\otimes 2}$} given by
\begin{equation*}
R(u)=1-u^{-1}P,\qquad \text{where} \ P=P_{12}=\sum_{i,j=1}^{n} e_{ij}\otimes e_{ji}.
\end{equation*}
The element $R(u)$ satisfies the \emph{Yang--Baxter equation}.
\begin{Proposition}[{\cite[Proposition 1.2.1]{Molev_07}}]
For any commuting indeterminates $u$, $v$ and $w$,
\begin{equation*}
R_{12}(u-v)R_{13}(u-w)R_{23}(v-w)=R_{23}(v-w)R_{13}(u-w)R_{12}(u-v).
\end{equation*}
\end{Proposition}

We also define for any $m\ge 2$ commuting indeterminates $u_1,\dots,u_m$ the following rational function:
\begin{equation*}
R(u_1,u_2,\dots,u_m)=(R_{m-1,m})(R_{m-2,m}R_{m-2,m-1})\cdots(R_{1m}\cdots R_{12}),
\end{equation*}
where we use shorthand notation $R_{ij}=R_{ij}(u_i-u_j)$. Note that for $m=1$, this is just the Yang $R$-matrix, $R(u_1,u_2)=R_{12}(u_1-u_2)$.

Using the Yang--Baxter equation, one can give an alternative definition of $R(u_1,u_2,\dots,u_m)$.
\begin{Lemma}[{\cite[Section~1.6]{Molev_07}}] In the notation above, the following identity holds:
\begin{equation*}
R(u_1,u_2,\dots,u_m)=(R_{12}\cdots R_{1m})\cdots(R_{m-2,m-1}R_{m-2,m})(R_{m-1,m}).
\end{equation*}
\end{Lemma}
In fact, the operator $A_m$ is a certain specialization of $R(u_1,\dots,u_m)$.

\begin{Proposition}[{\cite[Proposition 1.6.2]{Molev_07}}]\label{anti_sym}
If $u_i-u_{i+1}=1$ for all $i=1,\dots,m-1$, then
\begin{equation*}
R(u_1,u_2,\dots,u_m)=A_m.
\end{equation*}
\end{Proposition}

\subsection[Yangian Y(gl\_n)]{Yangian $\boldsymbol{\operatorname{Y}(\mathfrak{gl}_{n})}$}
The Yangian $\operatorname{Y}(\mathfrak{gl}_{n})$ is an associative algebra generated by elements \smash{$t_{ij}^{(r)}$}, $1\leq i,j\leq n$, $r\geq 1$, which satisfy certain relations. These relations can be succinctly written using the so-called \emph{ternary $($RTT$)$ relation}. To state it, we need to introduce more notation.

For each $i,j\in\{1,\dots,n\}$ consider a formal Laurent series
\begin{equation*}
t_{ij}(u)=\delta_{ij}+\sum_{r\geq 1}t_{ij}^{(r)}u^{-r}\in\operatorname{Y}(\mathfrak{gl}_{n})\bigl[\hspace{-0.3mm}\bigl[u^{-1}\bigr]\hspace{-0.3mm}\bigr].
\end{equation*}
Let \smash{$T(u)=\sum_{i,j=1}^{n}t_{ij}(u)\otimes e_{ij}\in\operatorname{Y}(\mathfrak{gl}_{n})\bigl[\hspace{-0.3mm}\bigl[u^{-1}\bigr]\hspace{-0.3mm}\bigr]\otimes\End\mathbb{C}^n$} be an $n\times n$ matrix whose $(i,j)$-entry equals $t_{ij}(u)$.

Now we can state the defining relations of $\operatorname{Y}(\mathfrak{gl}_{n})$ (see \cite[Proposition 1.2.2]{Molev_07}):
\begin{equation}\label{eq:rtt=ttr}
R(u-v)T_1(u)T_2(v)=T_2(v)T_1(u)R(u-v).
\end{equation}
Here both sides are viewed as elements of \smash{$\End(\mathbb{C}^{n})^{\otimes 2}\otimes\operatorname{Y}(\mathfrak{gl}_{n})\bigl[\hspace{-0.3mm}\bigl[u^{-1},v^{-1}\bigr]\hspace{-0.3mm}\bigr]$}. Expanding them in~$u^{-1}$ and~$v^{-1}$ gives infinitely many relations between generators \smash{$t_{ij}^{(r)}$}.

We need the following generalization of the ternary relation \eqref{eq:rtt=ttr}.
\begin{Proposition}\label{gen_ternary_rel}
For any $m\ge 2$ commuting indeterminates $u_1,\dots,u_m$,
\begin{equation*}
R(u_1,\dots,u_m)T_1(u_1)\cdots T_m(u_m)=T_m(u_m)\cdots T_1(u_1)R(u_1,\dots,u_m).
\end{equation*}
\end{Proposition}

\subsection{The evaluation homomorphism}
The relation between $\operatorname{Y}(\mathfrak{gl}_{n})$ and $U(\mathfrak{gl}_{n})$ is given by the so-called \emph{evaluation homomorphism} $\mathrm{ev}\colon \operatorname{Y}(\mathfrak{gl}_n)\to U(\mathfrak{gl}_n)$. The map $\mathrm{ev}$ is defined by the formula $\mathrm{ev}(t_{ij}(u))=\delta_{ij}+u^{-1}E_{ij}$, i.e., \smash{$\mathrm{ev}\bigl(t_{ij}^{(r)}\bigr)=\delta_{1r}E_{ij}$}. One can verify that $\mathrm{ev}\colon \operatorname{Y}(\mathfrak{gl}_n)\to U(\mathfrak{gl}_n)$ is a surjective algebra homomorphism (see also \cite[Proposition 1.1.3]{Molev_07}). Note that $\mathrm{ev}$ maps the $n\times n$ matrix $T(u)$ to $1+u^{-1}E$.

\subsection{Bethe subalgebras}

Define the \emph{trace} on $(\End\mathbb{C}^n)^{\otimes m}$ as the linear map $\tr_m$ defined by
\begin{equation*}
\tr_m\,(e_{i_1j_1}\otimes\dots\otimes e_{i_mj_m})=\delta_{i_1j_1}\cdots\delta_{i_mj_m}.
\end{equation*}
Observe that if we view $A\in(\End\mathbb{C}^n)^{\otimes m}$ as an element of $\End\,(\mathbb{C}^n)^{\otimes m}$, then $\tr_m(A)$ is indeed the trace of the operator $A$ acting on $(\mathbb{C}^n)^{\otimes m}$.
We extend $\tr_m$ to a map $\tr_m\colon\mathcal{A}\otimes(\End\mathbb{C}^n)^{\otimes m}\to\mathcal{A}$ where $\mathcal{A}$ is an arbitrary associative algebra,
\begin{equation*}
\tr_m\,(a\otimes e_{i_1j_1}\otimes\dots\otimes e_{i_mj_m})=\delta_{i_1j_1}\cdots\delta_{i_mj_m}a,~a\in\mathcal{A}.
\end{equation*}

Consider the elements $\tau_k(u)$ in \smash{$\operatorname{Y}(\mathfrak{gl}_n)\bigl[\hspace{-0.3mm}\bigl[u^{-1}\bigr]\hspace{-0.3mm}\bigr]$} defined as follows:
\begin{equation*}
\tau_k(u)=\frac{1}{k!}\tr_k A_k T_1(u)\cdots T_k(u-k+1).
\end{equation*}

\begin{Proposition}[{\cite[Proposition 1.14.1]{Molev_07}}]\label{bethe_commut}
The coefficients of the series $\tau_1(u),\dots,\tau_n(u)\in\operatorname{Y}(\mathfrak{gl}_n)\bigl[\hspace{-0.3mm}\bigl[u^{-1}\bigr]\hspace{-0.3mm}\bigr]$ generate a commutative subalgebra in $\operatorname{Y}(\mathfrak{gl}_{n})$.
\end{Proposition}

Observe that for any complex $n\times n$ matrix $C$ the matrix $C\cdot T(u)$ also satisfies the ternary relation \eqref{eq:rtt=ttr}. Indeed, this is a consequence of the fact that the element $C_1C_2\cdots C_m$ commutes with the image of $\mathbb{C}[\mathfrak{S}_m]$ in $\End(\mathbb{C}^n)^{\otimes m}$. Thus, one can generalize Proposition~\ref{bethe_commut} to construct commutative subalgebras of the Yangian which depend on a parameter in $\End\mathbb{C}^{n}$.
\begin{Corollary}[{\cite[Proposition~1.14.2]{Molev_07}}]\label{bethe_modif_cor1}
For any complex $n\times n$ matrix $C$, the coefficients of the series
\begin{equation*}
\tau_k(u,C)=\frac{1}{k!}\tr_{k} A_kC_1\cdots C_k T_1(u)\cdots T_k(u-k+1)
\end{equation*}
generate a commutative subalgebra of $\operatorname{Y}(\mathfrak{gl}_{n})$.
\end{Corollary}

\begin{Corollary}[{\cite[Proposition 1.14.3]{Molev_07}}]\label{bethe_modif_cor2}
For any complex $n\times n$ matrix $C$, the coefficients of the series
\begin{equation*}
\sigma_k(u,C)=\frac{1}{n!}\tr_{n} A_nT_1(u)\cdots T_k(u-k+1)C_{k+1}\cdots C_n
\end{equation*}
generate a commutative subalgebra of $\operatorname{Y}(\mathfrak{gl}_{n})$.
\end{Corollary}

\begin{Definition}
The commutative subalgebras $\operatorname{Y}(\mathfrak{gl}_{n})$ defined in Corollaries~\ref{bethe_modif_cor1} and~\ref{bethe_modif_cor2} are called the \emph{Bethe subalgebras}. (See \cite[Section~1.14]{Molev_07} and also \cite{Nazarov_Olshanski}.)
\end{Definition}

In calculations, we also use alternative formulas for the elements $\tau_k(u,C)$ and $\sigma_k(u,C)$, which are consequences of Proposition~\ref{gen_ternary_rel}.
\begin{Lemma}\label{bethe_elem_alter}
In the notation of Corollaries $\ref{bethe_modif_cor1}$ and $\ref{bethe_modif_cor2}$, we have
\begin{equation*}
\begin{split}
&\tau_k(u,C)
=\frac{1}{k!}\tr_{k} A_kC_1\cdots C_k T_k(u-k+1)\cdots T_1(u),
\\
&\sigma_k(u,C)
=\frac{1}{n!}\tr_{n} A_n T_k(u-k+1)\cdots T_1(u)C_{k+1}\cdots C_n.
\end{split}
\end{equation*}
\end{Lemma}

\subsection{Application to Capelli identities}
Let us show how the $R$-matrix formalism can be used in order to prove various identities related to the Capelli identity.

\begin{Lemma}\label{capelli_minor_char_poly_lem}
For any $k$-tuples $I=(i_1,\dots,i_k)$ and $J=(j_1,\dots,j_k)$ in $[n]^{k}$, we have the following identity in $U(\mathfrak{gl}_n)[z]$:
\begin{equation*}
\rdet\bigl[E_{i_{\alpha}j_{\beta}}+(z+\alpha-1)\delta_{i_{\alpha}j_{\beta}}\bigr]_{\alpha,\beta=1}^{k}=\cdet\bigl[E_{i_{\alpha}j_{\beta}}+(z+k-\beta)\delta_{i_{\alpha}j_{\beta}}\bigr]_{\alpha,\beta=1}^{k}.
\end{equation*}
In particular, these expressions are skew-symmetric in $i_1,\dots,i_k$ and $j_1,\dots,j_k$, respectively.
\end{Lemma}

\begin{proof}
We start with the ternary relation for $m=k$ (see Propositions \ref{anti_sym} and \ref{gen_ternary_rel}):
\begin{equation*}
A_kT_1(u)\cdots T_k(u-k+1)=T_k(u-k+1)\cdots T_1(u)A_k.
\end{equation*}
Now let us apply both sides as operators on $(\mathbb{C}^n)^{\otimes k}$ to the vector $e_{j_1}\otimes\cdots\otimes e_{j_k}$.
The left-hand side gives
\begin{gather*}
(A_kT_1(u)\cdots T_k(u-k+1))(e_{j_1}\otimes\cdots\otimes e_{j_k})
\\
\qquad{}=\sum_{l_1,\dots,l_k=1}^{n}t_{l_1j_1}(u)\cdots t_{l_kj_k}(u-k+1)\cdot A_k(e_{l_1}\otimes\cdots\otimes e_{l_k})
\\
\qquad{}=\sum_{l_1,\dots,l_k=1}^{n}\sum_{\sigma\in\mathfrak{S}_k}\sgn(\sigma)t_{l_1j_1}(u)\cdots t_{l_kj_k}(u-k+1)\cdot e_{l_{\sigma(1)}}\otimes\cdots\otimes e_{l_{\sigma(k)}}
\\
\qquad{}=\sum_{l_1,\dots,l_k=1}^{n}\sum_{\sigma\in\mathfrak{S}_k}\sgn(\sigma)t_{l_{\sigma(1)}j_1}(u)\cdots t_{l_{\sigma(k)}j_k}(u-k+1)\cdot e_{l_1}\otimes\cdots\otimes e_{l_k},
\end{gather*}
while the right-hand side gives
\begin{gather*}
(T_k(u-k+1)\cdots T_1(u)A_k)(e_{j_1}\otimes\cdots\otimes e_{j_k})
\\
\qquad{}=\sum_{\sigma\in\mathfrak{S}_k}\sgn(\sigma)(T_k(u-k+1)\cdots T_1(u))\bigl(e_{j_{\sigma(1)}}\otimes\cdots\otimes e_{j_{\sigma(k)}}\bigr)
\\
\qquad{}=\sum_{\sigma\in\mathfrak{S}_k}\sum_{l_1,\dots,l_k=1}^{n}\sgn(\sigma)t_{l_1j_{\sigma(1)}}(u-k+1)\cdots t_{l_kj_{\sigma(k)}}(u)\cdot e_{l_1}\otimes\cdots\otimes e_{l_k}.
\end{gather*}
Comparing the coefficients of $e_{i_1}\otimes\cdots\otimes e_{i_k}$, we obtain
\begin{gather}
\begin{aligned}[b]
& \sum_{\sigma\in\mathfrak{S}_k}\sgn(\sigma)t_{i_{\sigma(1)}j_1}(u)\cdots t_{i_{\sigma(k)}j_k}(u-k+1) \\
& \qquad{}=\sum_{\sigma\in\mathfrak{S}_k}\sgn(\sigma)t_{i_1j_{\sigma(1)}}(u-k+1)\cdots t_{i_kj_{\sigma(k)}}(u).
\end{aligned}\label{t_minor_skew_sym}
\end{gather}
Then, we apply the evaluation map to both sides and note that
\begin{gather*}
\mathrm{ev}(\text{LHS of \eqref{t_minor_skew_sym}})
=\bigl(u^{\underline{k}}\bigr)^{-1}\sum_{\sigma\in\mathfrak{S}_k}\sgn(\sigma)\prod_{p=1,\dots,k}\bigl(E_{i_{\sigma(p)}j_{p}}+(u-p+1)\delta_{i_{\sigma(p)}j_{p}}\bigr)
\\ \hphantom{\mathrm{ev}(\text{LHS of \eqref{t_minor_skew_sym}})}{}
=\bigl(u^{\underline{k}}\bigr)^{-1}\cdot\cdet\bigl[E_{i_{\alpha}j_{\beta}}+(u-\beta+1)\delta_{i_{\alpha}j_{\beta}}\bigr]_{\alpha,\beta=1}^{k},
\\
\mathrm{ev}(\text{RHS of \eqref{t_minor_skew_sym}})
=\bigl(u^{\underline{k}}\bigr)^{-1}\sum_{\sigma\in\mathfrak{S}_k}\sgn(\sigma)\prod_{p=1,\dots,k}\bigl(E_{i_{p}j_{\sigma(p)}}+(u-k+p)\delta_{i_pj_{\sigma(p)}}\bigr)
\\ \hphantom{\mathrm{ev}(\text{RHS of \eqref{t_minor_skew_sym}})}{}
=\bigl(u^{\underline{k}}\bigr)^{-1}\cdot\rdet\bigl[E_{i_{\alpha}j_{\beta}}+(u-k+\alpha)\delta_{i_{\alpha}j_{\beta}}\bigr]_{\alpha,\beta=1}^{k}.
\end{gather*}
Combining everything, we obtain
\begin{equation*}
\rdet\bigl[E_{i_{\alpha}j_{\beta}}+(u-k+\alpha)\delta_{i_{\alpha}j_{\beta}}\bigr]_{\alpha,\beta=1}^{k}=\cdet\bigl[E_{i_{\alpha}j_{\beta}}+(u-\beta+1)\delta_{i_{\alpha}j_{\beta}}\bigr]_{\alpha,\beta=1}^{k}.
\end{equation*}
Substituting $u=z+k-1$ finishes the proof of the lemma.
\end{proof}

\begin{Remark}
In particular, by plugging $z=0$ into the lemma, we obtain another proof of Corollary~\ref{cdet_rdet_cor}.
\end{Remark}

\section{Commutativity of big algebras}\label{s7}
The primary aim of this section is to prove the following crucial result.

\begin{Theorem}\label{thm:big_alg_commut}
The algebra $\mathscr{B}(\mathcal{P}(n,r))$ is commutative.
\end{Theorem}

Before we turn to the proof of the theorem, let us discuss its consequences for the big algebras of irreducible representations.

\subsection{Big algebras of irreducible representations}\label{subsect:big_alg_irreps}

For any dominant integral weight $\lambda\in\mathfrak{h}^*$, we set $\mathscr{B}(\lambda)\coloneqq\mathscr{B}(V(\lambda))$, i.e., the big algebra associated to the finite-dimensional irreducible $\mathrm{GL}_n$-module $V(\lambda)$ of highest weight $\lambda$.

Recall from Proposition~\ref{prop:big_alg_direct_sum_prop} that for any direct sum $V=\bigoplus_{\alpha}V_{\alpha}$ of $\mathrm{GL}_n$-modules the big algebras $\mathscr{B}(V_{\alpha})$ are homomorphic images of the big algebra $\mathscr{B}(V)$. Moreover, we have the inclusion $\mathscr{B}(V)\subset S(\mathfrak{gl}_n^{*})\otimes\bigl(\prod_{\alpha}\End V_{\alpha}\bigr)$ and the map $\mathscr{B}(V)\to\mathscr{B}(V_{\alpha})$ is induced by the projection $S(\mathfrak{gl}_n^{*})\otimes\bigl(\prod_{\alpha}\End V_{\alpha}\bigr)\to S(\mathfrak{gl}_{n}^{*})\otimes\End V_{\alpha}$. In particular, the image of an element $F\in\mathscr{B}(V)$ in $\mathscr{B}(V_{\alpha})$ can be regarded as the ``restriction of $F$ to $V_{\alpha}$''.
More precisely, for any $Y\in\mathfrak{gl}_{n}$ the linear operator $F(Y)\in\End V$ preserves each subspace $V_{\alpha}$ and its restriction to $V_{\alpha}$ is the image of $F$ in $\mathscr{B}(V_{\alpha})$ evaluated at $Y$.

Now let us apply these observations to the $\mathrm{GL}_n$-module decomposition of $\mathcal{P}(n,r)$, i.e., to $\mathcal{P}(n,r)\simeq\bigoplus_{\lambda\colon \ell(\lambda)\le r}V(\lambda)^{\oplus d_{\lambda}}$, where $V(\lambda)=V_{\mathrm{GL}_n}(\lambda)$ and $d_{\lambda}=\dim V_{\mathrm{GL}_r}(\lambda)$ (see \eqref{eq:p_gl_n_decomp} for more details). Then, Proposition~\ref{prop:big_alg_direct_sum_prop} implies that $\mathscr{B}(\lambda)$ is a quotient of $\mathscr{B}(\mathcal{P}(n,r))$ for any $\lambda$ of the form $(\lambda_1,\dots,\lambda_r,0,\dots,0)$. In particular, the big algebra of any finite-dimensional irreducible polynomial $\mathrm{GL}_n$-representation is a quotient of $\mathscr{B}(\mathcal{P}(n,n))$. The discussion above together with Theorem~\ref{thm:big_alg_commut} gives the following corollary.

\begin{Corollary}\label{cor:big_alg_commut_irrep}
The big algebra $\mathscr{B}(\lambda)$ of any polynomial finite-dimensional irreducible $\mathrm{GL}_n$-representation $V(\lambda)$ is commutative.
\end{Corollary}

\begin{Remark}
The commutativity of big algebras defined by means of the Kirillov--Wei operator was first shown by Hausel and Zveryk, see \cite[Theorem~2.3]{Hausel_24}. However, their proof involves quite non-trivial constructions such as the \emph{Feigin--Frenkel center} and \emph{Segal--Sugawara vectors}. Namely, they use explicit formulas for the Segal--Sugawara vectors for the Feigin--Frenkel center that were obtained by Yakimova in \cite{Yakimova}.
Our proof of the commutativity of big algebras in type~$A$ is different and essentially elementary. It relies on direct calculations and known facts about the Bethe subalgebras in Yangians and representation theory of $\mathfrak{gl}_n$.
\end{Remark}

\subsection{Proof of Theorem~\ref{thm:big_alg_commut}}

\begin{proof}
In view of Corollary~\ref{cor:F_operators}, it suffices to check that the operators $F_{p,q}$ commute. By \eqref{eq:F_operators_formula} and \eqref{capelli_minor_ident},
\begin{equation*}
F_{p,q}(Y)=\sum_{\substack{I_1,J_1\in\binom{[n]}{p}\\I_2,J_2\in\binom{[n]}{q}\\I_1\sqcup I_2=J_1\sqcup J_2}}\sgn\binom{I_1~I_2}{J_1~J_2}\det Y_{I_1J_1}\cdot L(\Pi_{J_2I_2}).
\end{equation*}
Note that in view of Proposition~\ref{prop:restr_cartan}, we only need to check the commutativity for $Y\in\mathfrak{h}$, i.e., for diagonal matrices. According to Proposition~\ref{f_cartan_restr_prop}, for $Y=\diag(z_1,\dots,z_n)$, we have
\begin{equation*}
F_{p,q}(Y)=\sum_{\substack{I\in\binom{[n]}{p}\\J\in\binom{[n]}{q}\\I\cap J=\varnothing}}\det Y_{II}\cdot L(\Pi_{JJ})=\sum_{\substack{I\in\binom{[n]}{p}\\J\in\binom{[n]}{q}\\I\cap J=\varnothing}}\prod_{i\in I}z_i\cdot L(\Pi_{JJ}).
\end{equation*}
To show that the operators $F_{p,q}$ commute, we will show that (see also \eqref{eq:F_pl_via_bethe2} below)
\begin{equation}\label{eq:F_pl_via_bethe1}
L(\mathrm{ev}(\sigma_{n-p}(v,Y)))=\sum_{\ell=0}^{n-p}(v(v-1)\cdots(v-\ell+1))^{-1}\cdot F_{p,\ell}(Y),\qquad Y\in\mathfrak{h}
\end{equation}
and then use the commutativity of the Bethe subalgebra. (Here we use the notation from Corollary~\ref{bethe_modif_cor2}.)

Indeed, applying Lemma~\ref{bethe_elem_alter} to the diagonal matrix $C=Y=\diag(z_1,\dots,z_n)$ yields
\begin{equation*}
\begin{aligned}
n!\cdot\sigma_k(u,Y)
&=\tr_n A_n T_k(u-k+1)\cdots T_1(u)Y_{k+1}\cdots Y_n
\\
&=\sum_{i_1,\dots,i_n}\sum_{\sigma\in\mathfrak{S}_{n}}\sgn(\sigma)t_{i_1i_{\sigma(1)}}(u-k+1)\cdots t_{i_ki_{\sigma(k)}}(u)\prod_{j=k+1}^{n}\delta_{i_{j}i_{\sigma(j)}}z_{i_j}
\\
&=\sum_{(i_1,\dots,i_k)\in[n]^{\underline{k}}}\sum_{\sigma\in\mathfrak{S}_k}\sgn(\sigma)t_{i_1i_{\sigma(1)}}(u-k+1)\cdots t_{i_ki_{\sigma(k)}}(u)\cdot\prod_{i\in I^c}z_{i}
\\
&=k!\sum_{i_1<\dots<i_k}\sum_{\sigma\in\mathfrak{S}_k}\sgn(\sigma)t_{i_1i_{\sigma(1)}}(u-k+1)\cdots t_{i_ki_{\sigma(k)}}(u)\cdot\prod_{i\in I^c}z_{i},
\end{aligned}
\end{equation*}
where by $I^c$ we denote the complement of $I=\{i_1,\dots,i_k\}$ in $[n]$. Therefore,
\begin{equation*}
n!\cdot\mathrm{ev}(\sigma_k(u,Y))=k!\cdot(u(u-1)\cdots(u-k+1))^{-1}\sum_{I\in\binom{[n]}{k}}\Pi_{II}(k-1-u)\cdot\prod_{i\in I^c}z_i.
\end{equation*}
Here we used that the expression{\samepage
\begin{gather*}
\mathrm{ev}\biggl(\sum_{\sigma\in\mathfrak{S}_k}\sgn(\sigma)t_{i_1i_{\sigma(1)}}(u-k+1)\cdots t_{i_ki_{\sigma(k)}}(u)\biggr)\\
\qquad{}=\bigl(u^{\underline{k}}\bigr)^{-1}\cdot\rdet\bigl[E_{i_{\alpha}i_{\beta}}+(u-k+\alpha)\delta_{i_{\alpha}i_{\beta}}\bigr]_{\alpha,\beta=1}^{k}
\end{gather*}
is skew-symmetric in $i_1,\dots,i_k$ (see Lemma~\ref{capelli_minor_char_poly_lem} and the formula~\eqref{t_minor_skew_sym}).}

Thus,
\begin{gather*}
n!\cdot(-1)^{n-p}\binom{u}{n-p}\cdot\mathrm{ev}(\sigma_{n-p}(n-p-1-u,Y))=\sum_{I\in\binom{[n]}{n-p}}\Pi_{II}(u)\cdot\prod_{i\in I^c}z_i
\\
\qquad{}=\sum_{I\in\binom{[n]}{n-p}}\sum_{\ell=0}^{n-p}(-1)^{n-p-\ell}u^{\underline{n-p-\ell}}\sum_{K\in\binom{I}{\ell}}\Pi_{KK}\cdot\prod_{i\in I^c}z_i\\
\qquad{}=\sum_{\ell=0}^{n-p}(-1)^{n-p-\ell}u^{\underline{n-p-\ell}}\sum_{\substack{J\in\binom{n}{p}\\K\in\binom{[n]}{\ell}\\J\cap K=\varnothing}}\Pi_{KK}\prod_{j\in J}z_j.
\end{gather*}
Applying the map $L$, we obtain
\begin{gather}
n!\cdot(-1)^{n-p}\binom{u}{n-p}\cdot L(\mathrm{ev}(\sigma_{n-p}(n-p-1-u,Y)))\nonumber\\
\qquad{}=\sum_{\ell=0}^{n-p}(-1)^{n-p-\ell}u^{\underline{n-p-\ell}}\cdot F_{p,\ell}(Y),\label{eq:F_pl_via_bethe2}
\end{gather}
which is equivalent to \eqref{eq:F_pl_via_bethe1}.

Now Corollary~\ref{bethe_modif_cor2} implies that for any $0\le p_1,p_2\le n$, we have the identity
\begin{equation*}
\Biggl[\sum_{\ell_1=0}^{n-p_1}(-1)^{n-p_1-\ell_1}u^{\underline{n-p_1-\ell_1}}\cdot F_{p_1,\ell_1}(Y),\sum_{\ell_2=0}^{n-p_2}(-1)^{n-p_2-\ell_2}v^{\underline{n-p_2-\ell_2}}\cdot F_{p_2,\ell_2}(Y)\Biggr]=0.
\end{equation*}
Finally, the polynomials $\bigl\{u^{\underline{k}}v^{\underline{l}}\bigr\}_{k,l\ge 0}$ are linearly independent and hence ${[F_{p_1,q_1}(Y),F_{p_2,q_2}(Y)]\!=\!0}$ for all $p_1$, $p_2$, $q_1$, $q_2$.
\end{proof}

\subsection{Big algebras and Bethe subalgebras}

In fact, the images of elements of the Bethe subalgebra under the evaluation homomorphism can be regarded as elements of Kirillov algebra in the sense of the lemma below. This observation together with the formula \eqref{eq:F_pl_via_bethe1} implies that, roughly speaking, big algebras are homomorphic images of Bethe subalgebras of the Yangian $\operatorname{Y}(\mathfrak{gl}_{n})$.
The connection with Bethe subalgebras becomes more transparent if we recall that big algebras $\mathscr{B}(V)$ are homomorhic images of the Feigin--Frenkel center $\mathfrak{z}(\mathfrak{g}_{-})$. Then, Talalaev in \cite{Talalaev} has shown how one can obtain explicit formulas for generators of $\mathfrak{z}(\mathfrak{g}_{-})$ in type $A$ from the description of Bethe subalgebra of the Yangian. To~avoid the reference to the Feigin--Frenkel center, one may regard the big algebra~$\mathscr{B}(V)$ as a~family of commutative subalgebras in $\End V$ that depend on a parameter $\chi\in\Spec I(\mathfrak{g})=\mathfrak{g}\sslash G$. For a fixed $\chi$, the corresponding algebra is the homomorphic image of the \emph{quantum shift of argument subalgebra} $\mathcal{A}_{\chi}\subset U(\mathfrak{g})$ (see \cite[Section~4]{Rybnikov}). The work of Ilin and Rybnikov \cite[Theorem~B]{Ilin_Rybnikov} implies\footnote{We thank Leonid Rybnikov for explaining this fact to us.} that the latter can be obtained as a certain limit of Bethe subalgebras, namely as \smash{$\mathcal{A}_{\chi}=U(\mathfrak{g})\cap\lim_{\varepsilon\to 0}B(\exp(\varepsilon\chi))$}.
Our formulas establish a more direct link between big algebras and Bethe subalgebras of $\operatorname{Y}(\mathfrak{gl}_n)$.

\begin{Lemma}
Consider the Laurent expansion \smash{$\sigma_k(u,C)=\sum_{r=0}^{\infty}\sigma_k^{(r)}(C)u^{-r}$}. Then, the maps \smash{$C\mapsto L\bigl(\mathrm{ev}\bigl(\sigma_k^{(r)}\bigl(C^{t}\bigr)\bigr)\bigr)$} belong to the Kirillov algebra $\mathscr{C}(\mathcal{P}(n,r))$.
\end{Lemma}

\begin{proof}
Denote \smash{$\Phi(u,C)=(L\circ\mathrm{ev})\bigl(\sigma_k\bigl(u,C^{t}\bigr)\bigr)$}. Then we have to verify the condition \eqref{eq:equivar_condition} for the $\mathrm{GL}_n$-representation $\widetilde{\pi}=\widetilde{L}$, i.e., that for any $g\in\mathrm{GL}_n$ the identity
\begin{equation*}
\Phi(u,\Ad(g)(C))=L(\Ad(g)(\mathrm{ev}(\sigma_k(u,C))))=\widetilde{L}(g)\Phi(u,C)\widetilde{L}(g)^{-1},
\end{equation*}
holds. This can be checked via direct calculation using the definition of $\sigma_k$ and the cyclicity property of trace.
\end{proof}

\section[Big algebras of the symmetric powers of the vector representation]{Big algebras of the symmetric powers\\of the vector representation}\label{sect:sym_pow_vect_rep}

The goal of this section is to give a more explicit description of the big algebra for the case of symmetric powers $S^m(\mathbb{C}^{n})$. These representations are also particular examples of weight multiplicity free representations.

\subsection{Kirillov algebras in the weight multiplicity free case}\label{subsect:kirillov_wmf}
Note that Kirillov algebra $\mathscr{C}(V)$ is not in general commutative. However, if $V$ is weight multiplicity free, then it is the case (see \cite[Corollary 1]{Kirillov_01}, \cite[Introduction]{Panyushev} and \cite[Theorem~4.1]{Rozhkovskaya}). Moreover, the following holds.

\begin{Proposition}
Let $V$ be a weight multiplicity free representation of $\mathfrak{g}=\mathfrak{gl}_n$. Then, for any $A\in\mathscr{C}(V)$ and any $Y\in\mathfrak{h}$, the operator $A(Y)\in\End V$ is diagonal in the weight basis of $V$. In~particular, the Kirillov algebra $\mathscr{C}(V)$ is commutative.
\end{Proposition}

\begin{proof}
This is an immediate consequence of Proposition~\ref{prop:restr_cartan}.
\end{proof}

\begin{Remark}
Note that since $V$ is weight multiplicity free, the elements of the weight basis of~$V$ are determined uniquely up to multiplication by a non-zero scalar.
\end{Remark}

As was discussed in Section~\ref{subsect:big_alg_irreps}, we can treat the elements of $\mathscr{B}(\mathcal{P}(n,r))$ as elements of the big algebra $\mathscr{B}(\lambda)$ by restricting the corresponding linear operators to a subrepresentation~$V$ of~$\mathcal{P}(n,r)$ which is isomorphic to $V(\lambda)$. For such a subrepresentation $V$ and an element ${F\in\mathscr{B}(n,r)}$, we denote the image of $F$ in $\mathscr{B}(\lambda)$ by $F|_{V(\lambda)}$ (or simply by~$F$ if the context is clear).

\begin{Corollary}\label{cor:sym_pow_rep}
The big algebra $\mathscr{B}(\mathcal{P}(n,1))$ coincides with the medium algebra $\mathscr{M}(\mathcal{P}(n,1))$.
\end{Corollary}

\begin{proof}
By Corollary~\ref{cor:F_generators_r} applied to $r=1$, both algebras $\mathscr{M}(\mathcal{P}(n,1))$ and $\mathscr{B}(\mathcal{P}(n,1))$ are generated by the operators $\{F_{p,0}\}_{1\le p\le n}\cup\{F_{p,1}\}_{0\le p\le n-1}$. Hence, $\mathscr{M}(\mathcal{P}(n,1))=\mathscr{B}(\mathcal{P}(n,1))$.
\end{proof}

Recall that the medium algebra $\mathscr{M}(V)$ is the subalgebra of $\mathscr{C}(V)$ generated by $S(\mathfrak{gl}_n^{*})^{\mathfrak{gl}_n}$ and~$\DKir\bigl(S(\mathfrak{gl}_n^{*})^{\mathfrak{gl}_n}\bigr)$ (see \eqref{eq:medium_algebra_def}). For a dominant integral weight $\lambda$, we denote $\mathscr{M}(\lambda)\coloneqq\mathscr{M}(V(\lambda))$. The following description of medium algebra was suggested by Hausel in \cite{Hausel_24}.

\begin{Conjecture}[{\cite[Theorem~1.1.3]{Hausel_24}}]\label{conj:kirillov_alg_center}
For any dominant integral weight $\lambda$, the medium algebra~$\mathscr{M}(\lambda)$ is the center of the Kirillov algebra $\mathscr{C}(\lambda)$.
\end{Conjecture}

A consequence of this conjecture is the following fact about big algebras for weight multiplicity free representations.

\begin{Corollary}[{conditional on Conjecture \ref{conj:kirillov_alg_center}}]\label{cor:wmf_medium_equal_kirillov}
Let $V$ be an irreducible weight multiplicity free representation. Then, $\mathscr{M}(V)=\mathscr{B}(V)=\mathscr{C}(V)$. In other words, the medium, big and Kirillov algebras of $V$ all coincide.
\end{Corollary}

\begin{Remark}\label{rem:conj_status}
We were informed by Hausel that the intended argument for \cite[Theorem~1.1.3]{Hausel_24} had a gap. However, the conjecture is supported by large computational evidence. Besides that, in \cite[Theorem~1.2]{Hausel_25} an analogous statement is shown for the center of the algebra $(U(\mathfrak{g})\otimes\End V(\lambda))^{G}$. The only consequence of Conjecture \ref{conj:kirillov_alg_center} that is actually relevant for us in this section is that the big algebra $\mathscr{B}(m\varpi_1)$ coincides with the medium algebra $\mathscr{M}(m\varpi_1)=Z(\mathscr{C}(m\varpi_1))$. However, the latter fact holds as a consequence of Corollary~\ref{cor:sym_pow_rep}.
\end{Remark}

\subsection[Description of the algebra B(S\^{}m(C\^{}n))]{Description of the algebra $\boldsymbol{\mathscr{B}(S^m(\mathbb{C}^n))}$}
Now we apply our results to symmetric powers of the vector representation of $\mathfrak{gl}_{n}$.
From now~on, we only consider $\mathcal{P}(n,1)$, i.e., we set $r=1$.
Denote for brevity $x_{j}=x_{j1}$ and $\partial_{j}=\partial_{j1}$, $1\le j\le n$.
Then, $\mathcal{P}(n,r)=\mathcal{P}(n,1)=\mathbb{C}[x_1,\dots,x_n]$ is isomorphic to the symmetric algebra ${S(\mathbb{C}^{n})=\bigoplus_{m\ge 0}S^m(\mathbb{C}^{n})}$. Note that the latter decomposition is also compatible with the $\mathrm{GL}_n$-action. Namely, the $m$-th summand ${S^m(\mathbb{C}^n)\simeq\mathbb{C}[x_1,\dots,x_n]_{m}}$ is isomorphic to $V(m\varpi_1)$. In~particular, each irreducible representation of $\mathrm{GL}_n$ occurs in $\mathcal{P}(n,1)$ with the multiplicity at most $1$.

Denote by $\{\varepsilon_1,\dots,\varepsilon_n\}$ the basis of the weight lattice in~$\mathfrak{h}^*$ that is dual to $\{E_{11},\dots,E_{nn}\}$ in~$\mathfrak{h}$. The fundamental weights $\varpi_1,\dots,\varpi_n$ are given by $\varpi_k=\varepsilon_1+\dots+\varepsilon_k$, $1\le k\le n$.
Recall that the Cartan subalgebra $\mathfrak{h}$ consists of diagonal matrices with $n$ entries. We identify $S(\mathfrak{h}^{*})$ with a polynomial ring $\mathbb{C}[t_1,\dots,t_n]$ using the natural coordinates on $\mathfrak{h}$.
Then, by the Chevalley restriction theorem \cite[Chapter~VI, \S23.1]{Humphreys}, $S(\mathfrak{gl}_n^*)^{\mathfrak{gl_n}}\simeq S(\mathfrak{h}^{*})^{W}=\mathbb{C}[t_1,\dots,t_n]^{\mathfrak{S}_n}$.

\begin{Proposition}
For any positive integer $m$, the big algebra $\mathscr{B}(m\varpi_1)$ is generated by the operators $\{F_{p,1}\}_{0\le p\le n-1}$ restricted to $V(m\varpi_1)$ over $S(\mathfrak{gl}_n^{*})^{\mathfrak{gl}_n}$.
\end{Proposition}

\begin{proof}
Observe that $\mathcal{P}(n,1)\simeq S(\mathbb{C}^n)$ decomposes into the direct sum $\bigoplus_{m\geq 0}V(m\varpi_1)$ as a $\mathfrak{gl}_n$-module. Thus, the second part of the proposition now follows from the proof of Corollary~\ref{cor:sym_pow_rep}.
\end{proof}

The following proposition gives a more explicit description of the big algebra for symmetric powers of the vector representation.

\begin{Proposition}\label{prop:big_alg_sym_power_descr}
Assume $\Func(\wt(m\varpi_1),S(\mathfrak{h}^{*}))$ is the algebra of $S(\mathfrak{h}^{*})$-valued functions on $\wt(m\varpi_1)$. There exists an injective $S(\mathfrak{h}^{*})^{W}$-algebra homomorphism
\[
\mathfrak{i}_{m}\colon \ \mathscr{B}(m\varpi_1)\to\Func(\wt(m\varpi_1),S(\mathfrak{h}^{*})).
\] Namely, the map $\mathfrak{i}_{m}$ sends an element $F\in\mathscr{B}(m\varpi_{1})$ to a function $\hat{F}\in\Func(\wt(m\varpi_1),S(\mathfrak{h}^{*}))$ such that the operator $F(Y)$ acts as multiplication by \smash{$\bigl(\hat{F}(\mu)\bigr)(Y)$} on $V_{\mu}(m\varpi_1)$ for all $\mu\in\wt(m\varpi_1)$ and $Y\in\mathfrak{h}$.

The image $\mathfrak{i}_{m}(\mathscr{B}(m\varpi_1))\simeq\mathscr{B}(m\varpi_1)$ is the subalgebra of $\Func(\wt(m\varpi_1),S(\mathfrak{h}^{*}))$ generated by the following two subalgebras:
\begin{itemize}\itemsep=0pt
\item
the subalgebra $\mathcal{F}^{(m)}_{0}$ of constant maps $\mu\mapsto f$ for $f\in S(\mathfrak{h}^{*})^{W}\simeq\mathbb{C}[t_1,\dots,t_n]^{\mathfrak{S}_{n}}$, and
\item
the subalgebra $\mathcal{F}^{(m)}_{1}$ of maps $\mu\mapsto f(\underbrace{t_1,\dots,t_1}_{\mu_1},\dots,\underbrace{t_n,\dots,t_n}_{\mu_n})$ for $f\in\Lambda_m$.\vspace{-3mm}
\end{itemize}
Here $\Lambda_{m}$ stands for the ring of symmetric polynomials in $m$ variables.
\end{Proposition}
\begin{Remark}\label{remark:plethystic_subst}
{\samepage Using the language of symmetric functions, one can interpret the expression
\[
f(\underbrace{t_1,\dots,t_1}_{\mu_1},\dots,\underbrace{t_n,\dots,t_n}_{\mu_n})
\]
as the \emph{plethystic substitution} $f[\mu_1t_1+\dots+\mu_nt_n]$, where $f$ is a symmetric function (see \cite[Section~I.8]{Macdonald} for definitions).}
\end{Remark}

\begin{proof}
Using Proposition~\ref{prop:restr_cartan}, we can view $\mathscr{B}(m\varpi_1)$ as a subalgebra of $S(\mathfrak{h}^{*})\otimes\End_{\mathfrak{h}}S^m(\mathbb{C}^n)$.
According to Proposition~\ref{f_cartan_restr_prop}, for $Y=\diag(t_1,\dots,t_n)\in\mathfrak{h}$, we have
\begin{equation*}
\begin{split}
&F_{p,0}(Y)
=\sum_{I\in\binom{[n]}{p}}\prod_{i\in I}t_{i}=e_{p}(t_1,\dots,t_n),\qquad\text{and}
\\
&F_{p,1}(Y)=\sum_{\substack{I\in\binom{[n]}{p}\\j\in [n]\setminus I}}\biggl(\prod_{i\in I}t_i\biggr)\cdot x_{j}\partial_{j}=\sum_{j=1}^{n}e_{p}\bigl(t_1,\dots,\widehat{t}_{j},\dots,t_n\bigr)\cdot x_{j}\partial_{j}.
\end{split}
\end{equation*}
Now we restrict these operators to the subrepresentation $S^{m}(\mathbb{C}^n)\subset\mathcal{P}(n,1)\simeq\mathbb{C}[x_1,\dots,x_n]$.
Consider the monomial basis $\{x^{\mu}\}_{\mu\in\wt(m\varpi_1)}$ of $S^{m}(\mathbb{C}^n)\simeq\mathbb{C}[x_1,\dots,x_n]_{m}$. Here we denote $x^{\mu}\coloneqq x_{1}^{\mu_1}\dots x_{n}^{\mu_{n}}$ and $\mu=(\mu_1,\dots,\mu_n)$ runs over $\wt(m\varpi_{1})$, that is, over all $n$-tuples of non-negative integers that add up to $m$. Clearly, the $\mu$-eigenspace of $\mathbb{C}[x_1,\dots,x_n]_{m}$ is spanned by~$x^{\mu}$ for all $\mu\in\wt(m\varpi_1)$.

The operator $x_{j}\partial_{j}$ acts on $\{x^{\mu}\}_{\mu\in\wt(m\varpi_{1})}$ as $x_{j}\partial_{j}(x^{\mu})=\mu_{j}x^{\mu}$. Thus, the operators $F_{p,0}$ and~$F_{p,1}$ are diagonal in the basis $\{x^{\mu}\}_{\mu\in\wt(m\varpi_{1})}$ and the corresponding matrices are as follows:
\begin{equation*}
\begin{split}
& F_{p,0}(Y)\leftrightarrow(e_p(t_1,\dots,t_n)\mid \mu\in\wt(m\varpi_1)),
\\
& F_{p,1}(Y)\leftrightarrow\Biggl(\sum_{j=1}^{n}\mu_{j}e_p\bigl(t_1,\dots,\widehat{t}_{j},\dots,t_n\bigr)\mid \mu\in\wt(m\varpi_1)\Biggr).
\end{split}
\end{equation*}
Let $\hat{F}_{p,q}$ (here $q=0,1$) be the element of $\Func(\wt(m\varpi_1),S(\mathfrak{h}^{*}))$ that corresponds to $F_{p,q}$. In~other words, we set $\hat{F}_{p,0}(\mu)=e_{p}(t_1,\dots,t_n)$ and \smash{$\hat{F}_{p,1}(\mu)=\sum_{j=1}^{n}\mu_{j}e_p\bigl(t_1,\dots,\widehat{t}_{j},\dots,t_n\bigr)$}.
Using the fact that the big algebra $\mathscr{B}(m\varpi_1)$ is generated by $\{F_{p,0}\}_{1\le p\le n}$ and $\{F_{p,1}\}_{0\le p\le n-1}$, we can identify $\mathscr{B}(m\varpi_{1})$ with the subalgebra $\mathcal{A}$ of the algebra $\Func(\wt(m\varpi_1),S(\mathfrak{h}^{*}))$ generated by the elements \smash{$\bigl\{\hat{F}_{p,0}\bigr\}_{1\le p\le n}$} and \smash{$\bigl\{\hat{F}_{p,1}\bigr\}_{0\le p\le n-1}$}. Our goal is to prove that $\mathcal{A}$ is generated by \smash{$\mathcal{F}^{(m)}_0$} and~\smash{$\mathcal{F}^{(m)}_1$}.

The elements \smash{$\bigl\{\hat{F}_{p,0}\bigr\}_{1\le p\le n}$} are constant functions on $\wt(m\varpi_{1})$ and they generate the subalgebra \smash{$\mathcal{F}^{(m)}_{0}\subset \Func(\wt(m\varpi_1),S(\mathfrak{h}^{*}))$}. Thus, \smash{$\mathcal{F}^{(m)}_{0}\subset\mathcal{A}$}.
Next, observe that
\begin{equation*}
\hat{F}_{p,1}(\mu)=\sum_{j=1}^{n}\mu_{j}e_p\bigl(t_1,\dots,\widehat{t}_{j},\dots,t_n\bigr)=\biggl(\mu_1\frac{\partial}{\partial t_{1}}+\dots+\mu_n\frac{\partial}{\partial t_{n}}\biggr)(e_{p+1}(t_1,\dots,t_n)).
\end{equation*}
Recall that the elementary symmetric polynomials $e_q(t_1,\dots,t_n)$ generate the ring $\mathbb{C}[t_1,\dots,t_n]^{\mathfrak{S}_{n}}$. Then, the Leibniz rule implies that $\mathcal{A}$ contains all functions of the form
\begin{equation*}
\mu\mapsto\biggl(\mu_1\frac{\partial}{\partial t_{1}}+\dots+\mu_n\frac{\partial}{\partial t_{n}}\biggr)(f) \qquad \text{for}~f\in\mathbb{C}[t_1,\dots,t_n]^{\mathfrak{S}_{n}}.
\end{equation*}
In particular, by setting $f=p_r$, where $p_r(t_1,\dots,t_n)=t_1^{r}+\dots+t_n^{r}$ is the $r$-th power sum, we obtain the function
\begin{align}
\mu\mapsto \biggl(\mu_1\frac{\partial}{\partial t_{1}}+\dots+\mu_n\frac{\partial}{\partial t_{n}}\biggr)(p_r(t_1,\dots,t_n))&{}=\sum_{j=1}^{n}\mu_{j}rt_{j}^{r-1}\nonumber\\
&{}=rp_{r-1}(\underbrace{t_1,\dots,t_1}_{\mu_1},\dots,\underbrace{t_n,\dots,t_n}_{\mu_n}).\label{F_1_power_sum}
\end{align}
Since the power sums also generate the ring of symmetric functions, we conclude that \smash{$\mathcal{F}^{(m)}_{1}\subset\mathcal{A}$}.

It remains to show that the elements \smash{$\bigl\{\hat{F}_{p,0}\bigr\}{}_{1\le p\le n}$} and \smash{$\bigl\{\hat{F}_{p,1}\bigr\}{}_{0\le p\le n-1}$} also belong to the subalgebra of $\Func(\wt(m\varpi_1),S(\mathfrak{h}^{*}))$ generated by \smash{$\mathcal{F}^{\smash{(m)}}_{0}$} and \smash{$\mathcal{F}^{\smash{(m)}}_{1}$}. Indeed, this is clear for~\smash{$\hat{F}_{p,0}$} as it belongs to \smash{$\mathcal{F}^{\smash{(m)}}_{0}$}. As for the \smash{$\hat{F}_{p,1}$}, note first that the function \smash{$\mu\mapsto\bigl(\sum_{i=1}^{n}\mu_i\frac{\partial}{\partial t_{i}}\bigr)(p_{r}(t_1,\dots,t_n))$} belongs to \smash{$\mathcal{F}^{\smash{(m)}}_{1}$} due to the formula~\eqref{F_1_power_sum}. Since the power sums $p_r(t_1,\dots,t_n)$ generate the ring $\mathbb{C}[t_1,\dots,t_n]^{\mathfrak{S}_{n}}$, the Leibniz rule implies that the function \smash{$\mu\mapsto\bigl(\sum_{i=1}^{n}\mu_i\frac{\partial}{\partial t_{i}}\bigr)(e_{p+1}(t_1,\dots,t_n))$} belongs to the subalgebra generated by \smash{$\mathcal{F}^{(m)}_{0}$} and \smash{$\mathcal{F}^{(m)}_{1}$}.
\end{proof}

\begin{Remark}
In fact, the proof of the proposition above can be easily modified to obtain a similar description of the medium algebra $\mathscr{M}(\lambda)=\mathscr{M}(V(\lambda))$ for any $\lambda=(\lambda_1,\dots,\lambda_n)$ with $\lambda_1\ge\dots\ge\lambda_n\ge 0$.
\end{Remark}

As an application of Proposition~\ref{prop:big_alg_sym_power_descr}, one can explicitly describe the big algebra of the vector representation of $\mathfrak{gl}_n$. Since $\mathscr{B}(\varpi_1)=\mathscr{C}(\varpi_1)$, this description recovers a special case of a~result due to Panyushev \cite{Panyushev}. For later use in the proof of Theorem~\ref{thm:big_alg_sym_pow_iso}, we state it in the following~form.

\begin{Corollary}[{\cite[Theorem~2.6 for $\lambda=\varpi_1$]{Panyushev}}]\label{cor:big_alg_vect_rep}
The $S(\mathfrak{gl}_{n}^{*})^{\mathfrak{gl}_{n}}$-algebra $\mathscr{B}(\varpi_{1})$ is isomorphic to the $S(\mathfrak{h}^{*})^{W}$-algebra $S(\mathfrak{h})^{W(\varpi_{1})}$, where $W(\varpi_{1})$ denotes the stabilizer of $\varpi_{1}\in\mathfrak{h}^{*}$ in the Weyl group $W\simeq\mathfrak{S}_{n}$. More explicitly, we have isomorphisms
\begin{equation*}
\mathscr{B}(\varpi_{1})\simeq \mathbb{C}[t_{1},\dots,t_{n}]^{\mathfrak{S}_{1}\times\mathfrak{S}_{n-1}}\simeq\mathbb{C}[t_{1},\dots,t_{n}]^{\mathfrak{S}_{n}}[T]/\bigl(T^n-c_1T^{n-1}+\dots+(-1)^nc_n\bigr).
\end{equation*}
Here, $c_1,\dots,c_n\in S(\mathfrak{h}^{*})^{W}$ are defined as in Section~$\ref{subsect:inv_polys}$ and $T$ can be identified with $t_1$.
\end{Corollary}

\begin{proof}
For $m=1$, the set $\wt(m\varpi_1)=\wt(\varpi_1)$ consists of $n$ weights which lie in the $W$-orbit of~$\varpi_1$. Then, it is not difficult to see that the subalgebra \smash{$\mathcal{F}^{(1)}_{1}$} is isomorphic to $\mathbb{C}[t_1]$ in this case (namely, $F\mapsto F(\varpi_1)$ is an isomorphism). Therefore, $\mathscr{B}(\varpi_1)$ is isomorphic to the subalgebra of~${\mathbb{C}[t_1,\dots,t_{n}]}$ generated by \smash{$\mathbb{C}[t_1,\dots,t_n]^{\mathfrak{S}_{n}}$} and $\mathbb{C}[t_1]$, i.e., to \smash{$\mathbb{C}[t_{1},\dots,t_{n}]^{\mathfrak{S}_{1}\times\mathfrak{S}_{n-1}}$}.

The last assertion of the corollary follows from the fact that \smash{$\mathbb{C}[t_1,\dots,t_n]^{\mathfrak{S}_{1}\times\mathfrak{S}_{n-1}}$} is a free $\mathbb{C}[t_1,\dots,t_n]^{\mathfrak{S}_{n}}$-module generated by $1,t_1,\dots,t_1^{n-1}$ and the identity $t_1^n-c_1t_1^{n-1}+\dots+(-1)^nc_n=0$ in $\mathbb{C}[t_1,\dots,t_n]$.
\end{proof}

\begin{Remark}
The element $T$ can be described in terms of the isomorphism from Proposition~\ref{prop:big_alg_sym_power_descr}. Namely, it corresponds to the function on $\wt(\varpi_1)=\{\varepsilon_1,\dots,\varepsilon_n\}$ that sends $\varepsilon_i$ to $t_i$ for~${i=1,\dots,n}$.
\end{Remark}

\subsection[An isomorphism between B(m varpi\_1) and S\^{}m(B(varpi\_1)]{An isomorphism between $\boldsymbol{\mathscr{B}(m\varpi_1)}$ and $\boldsymbol{S^m(\mathscr{B}(\varpi_1))}$}

For each positive integer $\alpha$, define an element $P_{\alpha}\in\mathscr{C}(\mathcal{P}(n,1))$ by the formula
\begin{equation}\label{eq:P_alpha_formula}
P_{\alpha}(Y)=L(Y^{\alpha}),
\end{equation}
where $Y\in\mathfrak{gl}_{n}$ is viewed as an $n\times n$ matrix. One checks using \eqref{eq:equivar_condition} that $P_{\alpha}$ is indeed an element of the Kirillov algebra $\mathscr{C}(\mathcal{P}(n,1))$.

\begin{Proposition}\label{prop:power_sums}
The algebra $\mathscr{B}(m\varpi_1)$ is generated over the ring $S(\mathfrak{gl}_n^{*})^{\mathfrak{gl}_{n}}$ by the elements $P_1,\dots,P_m$ restricted to $V(m\varpi_1)$.
\end{Proposition}

\begin{proof}
Note that for $Y=\diag(t_1,\dots,t_n)\in\mathfrak{h}$, we have $P_{\alpha}(Y)=\sum_{i=1}^{n}t_i^{\alpha}\cdot x_i\partial_i$.
It now follows from Proposition~\ref{prop:restr_cartan} and the proof of Proposition~\ref{prop:big_alg_sym_power_descr} that the restriction of $P_{\alpha}$ to $V(m\varpi_1)$ coincides with the preimage of the function $\mu\mapsto\sum_{i=1}^{n}\mu_i t_i^{\alpha}$ under the map $\mathfrak{i}_{m}$. In that proof we also showed that these functions generate the algebra $\mathfrak{i}_{m}(\mathscr{B}(m\varpi_1))$ over $S(\mathfrak{h}^{*})^{W}$.
Therefore, the elements $P_{\alpha}|_{V(m\varpi_1)}$ generate the algebra $\mathscr{B}(m\varpi_1)$ over $S(\mathfrak{gl}_n^{*})^{\mathfrak{gl}_{n}}$.
\end{proof}

Using the obtained facts, we can relate the big algebras of $\mathscr{B}(\varpi_{1})$ and $\mathscr{B}(m\varpi_1)$.

\begin{Theorem}\label{thm:big_alg_sym_pow_iso}
The $S(\mathfrak{gl}_{n}^{*})^{\mathfrak{gl}_{n}}$-algebra $\mathscr{B}(m\varpi_{1})$ is isomorphic to the subalgebra $S^{m}(\mathscr{B}(\varpi_{1}))$ of $\mathfrak{S}_{m}$-invariant elements of the tensor product $S(\mathfrak{gl}_{n}^{*})^{\mathfrak{gl}_{n}}$-algebra $\mathscr{B}(\varpi_1)^{\otimes m}$.
\end{Theorem}

\begin{proof}
Throughout this proof, all tensor products are considered over the ring $S(\mathfrak{gl}_{n}^{*})^{\mathfrak{gl}_{n}}\simeq S(\mathfrak{h}^{*})^{W}$.
By Proposition~\ref{prop:big_alg_sym_power_descr}, for any positive integer $k$, we can view $\mathscr{B}(k\varpi_1)$ as a subalgebra $\Func(\wt(k\varpi_1),S(\mathfrak{h}^*))$. We use this identification for $k=1$ and $k=m$.

Consider the map $\Upsilon\colon\Func(\wt(\varpi_1),S(\mathfrak{h}^*))^{\otimes m}\to\Func(\wt(m\varpi_1),S(\mathfrak{h}^*))$ defined by the following formula:
\begin{gather*}
\Upsilon(f_1\otimes\dots\otimes f_m)(\mu)=\prod_{j=1}^{n}\prod_{l=1}^{\mu_j}f_{\mu_1+\dots+\mu_{j-1}+l}(\varepsilon_{j})
\end{gather*}
for $f_1,\dots,f_m\in\Func(\wt(\varpi_1),S(\mathfrak{h}^*))$, $\mu\in\wt(m\varpi_1)$.
It is clear that $\Upsilon$ defined in this way is an $S(\mathfrak{gl}_{n}^{*})^{\mathfrak{gl}_{n}}$-algebra homomorphism.

Now let us show that $\Upsilon$ maps $S^m(\mathscr{B}(\varpi_1))$ to $\mathscr{B}(m\varpi_1)$ under the identification above. Using the notation of Corollary~\ref{cor:big_alg_vect_rep}, any $f\in\mathscr{B}(\varpi_1)$ can be uniquely expressed as $a_{0}+a_{1}T+\dots+a_{n-1}T^{n-1}$ with $a_i\in S(\mathfrak{gl}_n^{*})^{\mathfrak{gl}_n}$.
Define the following elements:
\smash{$X_i=1\otimes\dots\otimes \underset{i}{T}\otimes\dots\otimes 1$}, ${i\in\{1,\dots,m\}}$.

\begin{Claim}
The subalgebra $S^{m}(\mathscr{B}(\varpi_{1}))$ is spanned by symmetric polynomials in $X_1,\dots,X_m$ as a~$S(\mathfrak{gl}_{n}^{*})^{\mathfrak{gl}_{n}}$-module.
\end{Claim}

\begin{proof}
As an $S(\mathfrak{gl}_{n}^{*})^{\mathfrak{gl}_{n}}$-module, the algebra $\mathscr{B}(\varpi_{1})^{\otimes m}$ is a free module spanned by \smash{$X_1^{i_1}\!\cdots\! X_m^{i_m}$} with $i_1,\dots,i_m\in\{0,1,\dots,n-1\}$. Since the action of $\mathfrak{S}_{m}$ permutes $X_1,\dots,X_m$, the $S(\mathfrak{gl}_{n}^{*})^{\mathfrak{gl}_{n}}$-module $S^m(\mathscr{B}(\varpi_{1})^{\otimes m})$ is free with generators \smash{$\sum_{\sigma\in\mathfrak{S}_{m}}X_{1}^{\smash{i_{\sigma(1)}}\vphantom{{}_0}}\cdots X_{m}^{\smash{i_{\sigma(m)}}\vphantom{{}_0}}$}, where $(i_1,\dots,i_m)$ runs over all $m$-tuples of integers such that $0\le i_1\le\dots\le i_m\le n-1$. \big(In particular, the rank of this module equals \smash{$\binom{n+m-1}{m}$}.\big) It is clear that these expressions are symmetric polynomials in~${X_1, \dots, X_m}$ and generate all of them since $e_k(X_1,\dots,X_m)$ are among the generators.
\end{proof}

Now let us prove that the image of the restriction of $\Upsilon$ to $S^m(\mathscr{B}(\varpi_1))$ coincides with $\mathscr{B}(m\varpi_1)$. To start with, observe that for any $\mu\in\wt(m\varpi_1)$, we have
\begin{equation*}
(\Upsilon(X_1)(\mu),\Upsilon(X_2)(\mu),\dots,\Upsilon(X_m)(\mu))=(\underbrace{t_1,\dots,t_1}_{\mu_1},
\dots,\underbrace{t_n,\dots,t_n}_{\mu_n}).
\end{equation*}
Therefore, for any polynomial $f$ in $m$ variables with coefficients in $S(\mathfrak{gl}_n^{*})^{\mathfrak{gl}_{n}}$, we have
\begin{equation}\label{eq:upsilon_image}
\Upsilon(f(X_1,\dots,X_m))(\mu)=
f(\underbrace{t_1,\dots,t_1}_{\mu_1},\dots,\underbrace{t_n,\dots,t_n}_{\mu_n}).
\end{equation}
In particular, if $f$ is symmetric, then \smash{$\Upsilon(f(X_1,\dots,X_m))\in\mathcal{F}_{1}^{(m)}\subset\mathscr{B}(m\varpi_1)$}. Hence, $\Upsilon$ maps elements of $S^m(\mathscr{B}(\varpi_1))$ onto $\mathscr{B}(m\varpi_1)$ (see also Proposition~\ref{prop:big_alg_sym_power_descr}).

Finally, let us show that $\Upsilon$ is injective on $S^m(\mathscr{B}(\varpi_1))$. Suppose that $f$ is a symmetric polynomial in $m$ variables with coefficients in $S(\mathfrak{gl}_n^{*})^{\mathfrak{gl}_{n}}$ such that $\Upsilon(f(X_1,\dots,X_m))\equiv 0$ on $\wt(m\varpi_1)$. Since $S^m(\mathscr{B}(\varpi_1))$ is a free $S(\mathfrak{gl}_n^{*})^{\mathfrak{gl}_{n}}$-module generated by \smash{$\sum_{\sigma\in\mathfrak{S}_{m}}X_{1}^{i_{\sigma(1)}}\cdots X_{m}^{i_{\sigma(m)}}$} with $0\le i_1\le\dots \le i_m\le n-1$,
we may assume that $f$ has degree at most $n-1$ in each of the variables. The symmetry and \eqref{eq:upsilon_image} now imply that $f(a_1,\dots,a_m)$ vanishes for any $(a_1,\dots,a_m)\in\{t_1,\dots,t_n\}^{m}$. Now the following lemma, proved in Appendix~\ref{appendix}, implies that $f\equiv 0$.

\begin{Lemma}\label{lemma:vanishing_lemma}
Let $f$ be a polynomial in $m$ variables $u_1,\dots,u_m$ with coefficients in $S(\mathfrak{gl}_n^{*})^{\mathfrak{gl}_{n}}$ which has degree at most $n-1$ in each of the variables. If $f(a_1,\dots,a_m)$ vanishes for all $(a_1,\dots,a_m)\in\{t_1,\dots,t_n\}^{m}$, then $f\equiv 0$.
\end{Lemma}

Thus, the map $\Upsilon\colon S^m(\mathscr{B}(\varpi_1))\!\to\!\mathscr{B}(m\varpi_1)$ is indeed an isomorphism of $S(\mathfrak{gl}_n^{*})^{\mathfrak{gl}_n}$-algebras.\looseness=-1
\end{proof}

\begin{Remark}
Let us explain why the isomorphism described in the proof of Theorem~\ref{thm:big_alg_sym_pow_iso} is in fact ``natural''.
Indeed, $S^m(\mathscr{B}(\varpi_1))$ and $\mathscr{B}(m\varpi_1)$ can be regarded as certain subalgebras of~${S^m(S(\mathfrak{gl}_{n}^{*})\otimes\End V(\varpi_1))}$ and $S(\mathfrak{gl}_{n}^{*})\otimes\End V(m\varpi_1)$, respectively. If we evaluate these algebras at an element $Y=\diag(t_1,\dots,t_n)\in\mathfrak{h}$, we will obtain subalgebras of $S^{m}(\End_{\mathfrak{h}} V(\varpi_1))$ and $\End_{\mathfrak{h}} V(m\varpi_1)$, respectively, that depend on parameters $t_1,\dots,t_n$. The formula for the isomorphism $\Upsilon\colon S^m(\mathscr{B}(\varpi_1))\to\mathscr{B}(m\varpi_1)$ essentially comes from identifications $S^m(\End_{\mathfrak{h}} V(\varpi_1))\simeq\End_{\mathfrak{h}} S^m(V(\varpi_1))\simeq \End_{\mathfrak{h}} V(m\varpi_1)$.
\end{Remark}

\subsection{Generators and relations} In this subsection, we give a presentation of the big algebra $\mathscr{B}(m\varpi_1)$ in terms of generators and relations.
Recall that we identify $S(\mathfrak{gl}_{n}^{*})^{\mathfrak{gl}_n}$ with $S(\mathfrak{h}^*)^{W}=\mathbb{C}[t_1,\dots,t_n]^{\mathfrak{S}_n}$.

\begin{Theorem}\label{thm:big_sym_power_gen_rel}
The $S(\mathfrak{gl}_{n}^{*})^{\mathfrak{gl}_n}$-algebra $\mathscr{B}(m\varpi_1)$ is isomorphic to
\begin{equation*}
\mathbb{C}[t_1,\dots,t_n]^{\mathfrak{S}_n}[u_1,\dots,u_m]^{\mathfrak{S}_m}/J,
\end{equation*}
where $J$ is the ideal of $R[u_1,\dots,u_m]^{\mathfrak{S}_m}$ given by
\begin{equation}\label{eq:J_ideal_def}
J=\bigl\langle f_1+\dots+f_m,u_1f_1+\dots+u_mf_m,\dots,u_1^{m-1}f_1+\dots+u_m^{m-1}f_m\bigr\rangle
\end{equation}
and $f_i=(u_i-t_1)\cdots(u_i-t_n)$ for each $i=1,\dots,m$.
\end{Theorem}
\begin{Remark}
The group $\mathfrak{S}_m$ acts on $\mathbb{C}[t_1,\dots,t_n]^{\mathfrak{S}_n}[u_1,\dots,u_m]$ by permuting the generators $u_1,\dots,u_m$.
\end{Remark}
\begin{proof}[Proof of Theorem~\ref{thm:big_sym_power_gen_rel}]
Denote $R=S(\mathfrak{h}^{*})^{W}=\mathbb{C}[t_1,\dots,t_n]^{\mathfrak{S}_n}$.
By Proposition~\ref{prop:big_alg_sym_power_descr}, the big algebra $\mathscr{B}(m\varpi_1)$~is isomorphic to the subalgebra of $\Func(\wt(m\varpi_1),S(\mathfrak{h}^*))$ generated by the subalgebras~\smash{$\mathcal{F}^{\smash{(m)}}_{0}$} and~\smash{$\mathcal{F}^{\smash{(m)}}_{1}$}. Define a homomorphism $\Psi\colon R[u_1,\dots,u_m]^{\mathfrak{S}_m}\to\Func(\wt(m\varpi_1),S(\mathfrak{h}^*))$ which sends any symmetric polynomial $f(u_1,\dots,u_m)$ with coefficients in $R$ to the function\footnote{Compare this to \eqref{eq:upsilon_image}.}
\begin{equation*}
\mu\mapsto f(\underbrace{t_1,\dots,t_1}_{\mu_1},\dots,\underbrace{t_n,\dots,t_n}_{\mu_n}),\qquad\mu=(\mu_1,\dots,\mu_n)\in\wt(m\varpi_1).
\end{equation*}
It follows from the definitions of \smash{$\mathcal{F}_{0}^{(m)}$} and \smash{$\mathcal{F}_{1}^{(m)}$} that $\Psi$ defines a surjective homomorphism from \smash{$ R[u_1,\dots,u_m]^{\mathfrak{S}_m}$ to $\mathscr{B}(m\varpi_1)\subset\Func(\wt(m\varpi_1),S(\mathfrak{h}^*))$}. Therefore, to finish the proof of the theorem, it remains to show that $\ker\Psi=J$.

Indeed, $\Psi$ sends a symmetric polynomial $f(u_1,\dots,u_m)$ with coefficients in $R$ to the zero function on $\wt(m\varpi_1)$ if and only if
\[
f(\underbrace{t_1,\dots,t_1}_{\mu_1},\dots,\underbrace{t_n,\dots,t_n}_{\mu_n})\equiv 0 \qquad
\text{for any}\ \mu=(\mu_1,\dots,\mu_n)\in\wt(m\varpi_1).
\]
Since $f$ is symmetric, the last condition is equivalent to the fact that $f$ vanishes on the Cartesian product $\{t_1,\dots,t_n\}^{m}$. Thus, $\ker\Psi$ is the ideal of those symmetric polynomials in~$R[u_1,\dots,u_m]^{\mathfrak{S}_m}$ which vanish whenever we substitute $(u_1,\dots,u_m)\in\{t_1,\dots,t_n\}^{m}$. This ideal coincides with $J$ by Lemma~\ref{lemma:generators_J_ideal} below.
\end{proof}

\begin{Lemma}\label{lemma:generators_J_ideal}
In the notation of the proof of Theorem~$\ref{thm:big_sym_power_gen_rel}$, we have
\begin{enumerate}\itemsep=0pt
\item[$(i)$] $J=\langle f_1,\dots,f_m\rangle\cap R[u_1,\dots,u_m]^{\mathfrak{S}_m}$, where $\langle f_1,\dots,f_m\rangle$ is the ideal of $R[u_1,\dots,u_m]$ generated by $\{f_i\}_{i=1}^{m}$;

\item[$(ii)$] $J=\bigl\{f(u_1,\dots,u_m)\in R[u_1,\dots,u_m]^{\mathfrak{S}_m}\mid f|_{\{t_1,\dots,t_n\}^{m}}\equiv 0\bigr\}$.
\end{enumerate}
\end{Lemma}
\begin{proof}[Proof]
(i) It is clear from \eqref{eq:J_ideal_def} that $J\subset\langle f_1,\dots,f_m\rangle\cap R[u_1,\dots,u_m]^{\mathfrak{S}_m}$. Take any element $F\in\langle f_1,\dots,f_m\rangle\cap R[u_1,\dots,u_m]^{\mathfrak{S}_m}$. Then, we can write it as
\begin{equation*}
F=a_1f_1+\dots+a_mf_m,
\end{equation*}
where $a_1,\dots,a_m\in R[u_1,\dots,u_m]$. Averaging this expression over $\mathfrak{S}_m$ gives
\begin{equation*}
F=\frac{1}{m!}\sum_{\sigma\in\mathfrak{S}_m}\sum_{i=1}^{m}(\sigma\cdot a_i)(\sigma\cdot f_i).
\end{equation*}
However, it is clear from the equality $f_i=(u_i-t_1)\cdots(u_i-t_n)$ that $\mathfrak{S}_m$ permutes $f_1,\dots,f_m$, i.e., $\sigma\cdot f_i=f_{\sigma(i)}$ for all $i\in\{1,\dots,m\}$. Therefore, we can rewrite the equality above as follows:
\begin{equation*}
F=b_1f_1+\dots+b_mf_m,
\qquad \text{where}~b_i=\frac{1}{m!}\sum_{j=1}^{m}\sum_{\substack{\sigma\in\mathfrak{S}_m\\\sigma(j)=i}}(\sigma\cdot a_j)=\frac{1}{m!}\sum_{\tau\in\mathfrak{S}_m}\bigl(\tau^{-1}\cdot a_{\tau(i)}\bigr).
\end{equation*}
It is not difficult to see that for each $i\in\{1,\dots,m\}$ the element $b_i$ belongs to \smash{$R[u_1,\dots,u_m]^{\mathfrak{S}_m(i)}$}, where $\mathfrak{S}_m(i)=\Stab_{\mathfrak{S}_m}(i)$ is the stabilizer of $i\in\{1,\dots,m\}$ in $\mathfrak{S}_m$. Indeed, it follows from a~more general observation that $\sigma\cdot b_i=b_{\sigma(i)}$ for any $i\in\{1,\dots,m\}$ and any $\sigma\in \mathfrak{S}_m$.

Note that \smash{$R[u_1,\dots,u_m]^{\mathfrak{S}_m(1)}$} is in fact a free rank $m$ module over \smash{$R[u_1,\dots,u_m]^{\mathfrak{S}_m}$}. Moreover, the elements \smash{$1,u_1,\dots,u_1^{m-1}$} are free generators of this module. Therefore, there exist \smash{$g_0,\dots,g_{m-1}\in R[u_1,\dots,u_m]^{\mathfrak{S}_m}$} such that
\begin{equation*}
b_1=g_0+g_1u_1+\dots+g_{m-1}u_1^{m-1}.
\end{equation*}
The discussion above implies that $g_i=g_0+g_1u_i+\dots+g_{m-1}u_i^{m-1}$ for all $i$ and thus
\begin{equation*}
F=\sum_{i=1}^{m}b_if_i=\sum_{i=1}^{m}\bigl(g_0+g_1u_i+\dots+g_{m-1}u_i^{m-1}\bigr)f_i=\sum_{j=0}^{m-1}g_{j}\bigl(u_1^{j}f_1+\dots+u_{m}^{j}f_m\bigr)\in J.
\end{equation*}
Thus, $F\in J$ which concludes the proof of the first part of the lemma.

(ii) Observe that for each $i=1,\dots,m$ we have $f_i|_{\{t_1,\dots,t_n\}^{m}}\equiv 0$, where we regard $f_1,\dots,f_m$ as polynomials in variables $u_1,\dots,u_m$ with coefficients in $R$. Then, it follows from (i) that
\begin{equation*}
J\subset\bigl\{f(u_1,\dots,u_m)\in R[u_1,\dots,u_m]^{\mathfrak{S}_m}\mid f|_{\{t_1,\dots,t_n\}^{m}}\equiv 0\bigr\}.
\end{equation*}
To show the other inclusion, take any element $f(u_1,\dots,u_m)\in R[u_1,\dots,u_m]^{\mathfrak{S}_m}$ that vanishes on $\{t_1,\dots,t_n\}^{m}$. It is clear that there exists a polynomial $g(u_1,\dots,u_m)\in R[u_1,\dots,u_m]$ which has degree at most $n-1$ in each the variables $u_1,\dots,u_m$ and such that $f\equiv g\pmod{\langle f_1,\dots,f_m\rangle}$ in $R[u_1,\dots,u_m]$. Then, $g(u_1,\dots,u_m)$ vanishes on $\{t_1,\dots,t_n\}^{m}$ as well and Lemma~\ref{lemma:vanishing_lemma} implies that $g\equiv 0$. Thus, $f\in\langle f_1,\dots,f_m\rangle$. As $f\in R[u_1,\dots,u_m]^{\mathfrak{S}_m}$ by the assumption, we obtain $f\in\langle f_1,\dots,f_m\rangle\cap R[u_1,\dots,u_m]^{\mathfrak{S}_m}=J$, as needed.
\end{proof}

One can give a more explicit description of $\mathscr{B}(m\varpi_1)$.

\begin{Corollary}\label{cor:big_sym_power_gen_rel}
We have the following isomorphism of $S(\mathfrak{gl}_n^{*})^{\mathfrak{gl}_n}$-algebras:
\begin{gather}
\mathscr{B}(m\varpi_1)\label{eq:medium_ring_presentation_2}\\
\quad{}\simeq\mathbb{C}[c_1,\dots,c_n][P_1,\dots,P_m]/\big\langle P_{i+n}-c_1P_{i+n-1}+\dots+(-1)^nc_nP_{i}\mid i=0,\dots,m-1\big\rangle,\nonumber
\end{gather}
where $S(\mathfrak{gl}_n^{*})^{\mathfrak{gl}_n}$ is identified with $\mathbb{C}[c_1,\dots,c_n]$ $($see Section~$\ref{subsect:inv_polys})$ and $P_k$ is the element of $\mathscr{B}(m\varpi_1)$ given by the formula \eqref{eq:P_alpha_formula}.
\end{Corollary}
\begin{proof}
Denote as before $R=\mathbb{C}[c_1,\dots,c_n]$. We have \smash{$f_i=u_i^{n}-c_1u_i^{n-1}+\dots+(-1)^nc_n$} and hence
\begin{align*}
u_1^{j}f_1+\dots+u_m^{j}f_m={}&p_{j+n}(u_1,\dots,u_m)-c_1p_{j+n-1}(u_1,\dots,u_m)+\cdots\\
&{}{+}\,(-1)^nc_np_{j}(u_1,\dots,u_m).
\end{align*}
Since $R[u_1,\dots,u_m]^{\mathfrak{S}_m}$ is freely generated over $R$ by the power sums
\[
p_{k}(u_1,\dots,u_m)=u_1^{k}+\dots+u_m^{k},
\]
where $k=1,\dots,m$, the statement of the corollary follows from Theorem~\ref{thm:big_sym_power_gen_rel}.
\end{proof}

\begin{Remark}
It follows from \eqref{eq:P_alpha_formula} and the proofs of Proposition~\ref{prop:power_sums} and Theorem~\ref{thm:big_sym_power_gen_rel} the image of the $k$-th power sum $p_{k}(u_1,\dots,u_m)=u_1^k+\dots+u_m^k$ in $\mathscr{B}(m\varpi_1)$ coincides with $P_{k}$. In particular, any $P_{k}$ with $k>m$ is a polynomial expression in $P_1,\dots,P_m$ (also note that $P_0\equiv m$). Therefore, to write down the explicit relations on $c_1,\dots,c_n,P_1,\dots,P_m$, one has to express $P_{m+1},\dots,P_{n+m-1}$ in terms of $P_1,\dots,P_m$. This can be done using the identities for power sum symmetric polynomials in~$m$ variables.
\end{Remark}

\begin{Remark}
The algebra $\mathscr{B}(m\varpi_1)$ coincides with the medium algebra $\mathscr{M}(m\varpi_1)$ (see \eqref{eq:medium_algebra_def} and Corollary~\ref{cor:sym_pow_rep}). In a private communication, Jakub L\"{o}wit 
informed us that one can prove the isomorphism from Theorem~\ref{thm:big_sym_power_gen_rel} using the geometric interpretation of the medium algebra~$\mathscr{M}(m\varpi_1)$. Namely, one can use the fact that for any dominant integral weight $\lambda$ the algebra~$\mathscr{M}(\lambda)$ is isomorphic to the $\mathrm{GL}_n$-equivariant cohomology of the affine Schubert variety~$\operatorname{Gr}^{\leqslant\lambda}$ associated to $\lambda$.
\end{Remark}

\subsection{Examples}\label{subsect:examples} In this subsection, we use Corollary~\ref{cor:big_sym_power_gen_rel} to compute in several examples the explicit presentations of the big algebra $\mathscr{B}(m\varpi_1)$ for $\mathfrak{gl}_n$. Then, we deduce the presentations of the corresponding big algebras for $\mathfrak{sl}_n$ and compare them with the computations of Hausel--Rychlewicz \cite[Example~4.6]{Hausel_Rychlewicz}, Rozhkovskaya \cite[Proposition 4.2]{Rozhkovskaya} and Hausel \cite[Section~4.2.3]{Hausel_24}.

Recall that $V(m\varpi_1)$ is weight multiplicity free as a $\mathfrak{gl}_n$-module, hence $\mathscr{B}(m\varpi_1)$ coincides with the medium algebra $\mathscr{M}(m\varpi_1)$. The latter algebra is generated over $S(\mathfrak{gl}_n^{*})^{\mathfrak{gl}_n}$ by the operators $M_{0,1}|_{V(m\varpi_1)},\dots,M_{n-1,1}|_{V(m\varpi_1)}$ (see \eqref{eq:M_pq_operators}).
We denote $M_i\coloneqq M_{i,1}|_{V(m\varpi_1)}$ and call these elements \emph{medium operators}. Note that $M_{0,1}|_{V(m\varpi_1)}\equiv m\cdot\Id_{V(m\varpi_1)}$ is a scalar operator and it can be discarded from the list of generators.

In all three examples below, we present $\mathscr{B}(m\varpi_1)$ using as generators over the ring $S(\mathfrak{gl}_n^{*})^{\mathfrak{gl}_n}=\mathbb{C}[c_1,\dots,c_n]$ either $\{P_1,\dots,P_m\}$, or $\{M_{1},\dots,M_{n-1}\}$.
The transition between these presentations is done by means of the following identities for symmetric polynomials (the \emph{Girard--Waring formula} and its inverse):
\begin{gather}\label{eq:girard_waring_1}
(-1)^ke_k(w_1,\dots,w_N)=\sum_{\substack{i_1,\dots,i_k\ge 0\\i_1+2i_2+\dots+ki_k=k}}\frac{(-1)^{i_1+i_2+i_3+\cdots}}{i_1!\cdots i_k!}\prod_{j=1}^{k}\biggl(\frac{1}{j}p_{j}(w_1,\dots,w_N)\biggr)^{i_j},
\\
\frac{1}{k}p_k(w_1,\dots,w_N)\nonumber\\
\qquad{}=\sum_{\substack{i_1,\dots,i_N\ge 0\\i_1+2i_2+\dots+Ni_N=k}}\frac{(-1)^{i_2+i_4+\cdots}(i_1+\dots+i_N-1)!}{i_1!\cdots i_N!}\prod_{j=1}^{N}e_{j}(w_1,\dots,w_N)^{i_j}, \label{eq:girard_waring_2}
\end{gather}
where $k$ and $N$ are arbitrary positive integers (recall that $e_k$ and $p_k$ are the elementary and the power sum symmetric polynomials, respectively).
These equalities are consequences of the following power series identities:
\begin{gather*}
1-e_1(w_1,\dots,w_N)z+\dots+(-1)^{N} e_N(w_1,\dots,w_N)\\
\qquad=(1-w_1z)\cdots(1-w_Nz)=\exp\Biggl\{\sum_{i=1}^{N}\log(1-w_iz)\Biggr\}\\
\qquad=\exp\Biggl\{-\sum_{i=1}^{N}\sum_{k=1}^{\infty
}\frac{1}{k}w_i^kz^k\Biggr\}=\exp\Biggl\{-\sum_{k=1}^{\infty}\frac{1}{k}p_k(w_1,\dots,w_N)z^k\Biggr\}.
\end{gather*}

For any positive integer $N$ and non-negative integer $k$, define the polynomial \smash{$\Theta^{(N)}_{k}(v_1,\dots,v_N)$} in $N$ variables as the polynomial satisfying the equality
\begin{equation*}
\Theta^{(N)}_{k}(p_1(w_1,\dots,w_N),p_2(w_1,\dots,w_N),\dots,p_{N}(w_1,\dots,w_N))=p_{k}(w_1,\dots,w_N).
\end{equation*}
The existence of such polynomials follows from formulas \eqref{eq:girard_waring_1} and \eqref{eq:girard_waring_2} while the uniqueness follows from the fact that the ring of symmetric polynomial in $w_1,\dots,w_N$ is freely generated~by \smash{$\{p_i(w_1,\dots,w_N)\}_{i=1}^{N}$}.
For instance, for any $k\in\{1,\dots,N\}$, then \smash{$\Theta_{k}^{\smash{(N)}\vphantom{_0}}(v_1,\dots,v_N)\equiv v_k$} and for~${k=0}$, we have \smash{$\Theta^{\smash{(N)}\vphantom{_0}}_{0}(v_1,\dots,v_N)\equiv N$}.

The formula \eqref{eq:girard_waring_1} allows us to express medium operators $M_1,\dots,M_{n-1}$ as a $S(\mathfrak{gl}_n^{*})^{\mathfrak{gl}_n}$-linear combination of elements $P_1,\dots,P_{n-1}$.
\begin{Lemma}
The $k$-th medium operator $M_k$ in $\mathscr{B}(m\varpi_1)$ can be expressed as
\begin{gather}
M_k=\sum_{\substack{i_1,\dots,i_k,i_{k+1}\ge 0\\i_1+2i_2+\dots+(k+1)i_{k+1}=k+1}}\frac{(-1)^{(k+1)+(i_1+i_2+\dots+i_{k+1})}}{i_1!\cdots i_{k+1}!}\nonumber\\
\hphantom{M_k=\sum_{\substack{i_1,\dots,i_k,i_{k+1}\ge 0\\i_1+2i_2+\dots+(k+1)i_{k+1}=k+1}}}{}
\times\Biggl(\sum_{j=1}^{k+1}i_jP_{j-1}\biggl(\frac{1}{j}\theta_j\biggr)^{i_j-1}\cdot\prod_{\substack{1\le l\le k+1\\l\neq j}}\biggl(\frac{1}{l}\theta_{l}\biggr)^{i_l}\Biggr),\label{eq:M_P_transfer}
\end{gather}
where $k=1,\dots,n-1$ and $\theta_l$ is the element of $S(\mathfrak{gl}_n^*)^{\mathfrak{gl}_n}\simeq\mathbb{C}[t_1,\dots,t_n]^{\mathfrak{S}_n}$ given by the formula $\theta_l=t_1^l+\dots+t_n^l$.
\end{Lemma}
\begin{proof}
The formula \eqref{eq:girard_waring_1} implies the following identity in $\mathbb{C}[t_1,\dots,t_n]^{\mathfrak{S}_n}$:
\begin{gather}
e_{k+1}(t_1,\dots,t_n)=\sum_{\substack{i_1,\dots,i_k,i_{k+1}\ge 0\\i_1+2i_2+\dots+(k+1)i_{k+1}=k+1}}\frac{(-1)^{(k+1)+(i_1+i_2+\dots+i_{k+1})}}{i_1!\cdots i_{k+1}!}\nonumber\\
\hphantom{e_{k+1}(t_1,\dots,t_n)=\sum_{\substack{i_1,\dots,i_k,i_{k+1}\ge 0\\i_1+2i_2+\dots+(k+1)i_{k+1}=k+1}}}{}
\times\prod_{j=1}^{k+1}\biggl(\frac{1}{l}\bigl(t_1^l+\dots+t_n^l\bigr)\biggr)^{i_l},\label{eq:elem_power_expansion}
\end{gather}
which is equivalent to the equality
\begin{equation*}
c_{k+1}=\sum_{\substack{i_1,\dots,i_k,i_{k+1}\ge 0\\i_1+2i_2+\dots+(k+1)i_{k+1}=k+1}}\frac{(-1)^{(k+1)+(i_1+i_2+\dots+i_{k+1})}}{i_1!\cdots i_{k+1}!}\prod_{j=1}^{k+1}\biggl(\frac{1}{l}\theta_{l}\biggr)^{i_l}
\end{equation*}
in $S(\mathfrak{gl}_n^{*})^{\mathfrak{gl}_n}$. By Corollary~\ref{cor:F_operators}, we have $M_k=M_{k,1}=F_{k,1}$. The proof of Proposition~\ref{prop:big_alg_sym_power_descr} implies that the embedding of $\mathscr{B}(m\varpi_1)$ into $\Func(\wt(m\varpi_1),S(\mathfrak{h}^*))$ identifies $c_{k+1}$ and $M_k=\DKir(c_{k+1})$ with functions $\mu\mapsto e_{k+1}(t_1,\dots,t_n)$ and \smash{$\mu\mapsto \bigl(\sum_{i=1}^{n}\mu_i\frac{\partial}{\partial t_i}\bigr)(e_{k+1}(t_1,\dots,t_n))$}, respectively. Note also that $P_l$ corresponds to the function \smash{$\mu\mapsto\sum_{i=1}^{n}\mu_i t_i^l$} (see \eqref{eq:P_alpha_formula} and Proposition~\ref{prop:big_alg_sym_power_descr}). Applying the differential operator \smash{$\sum_{i=1}^{n}\mu_{i}\frac{\partial}{\partial t_{i}}$} to both sides of \eqref{eq:elem_power_expansion} gives \eqref{eq:M_P_transfer}.
\end{proof}

\begin{Remark}
It follows from \eqref{eq:M_P_transfer} that $\{M_1,\dots,M_{n-1}\}$ and $\{P_1,\dots,P_{n-1}\}$ are related to each other by an invertible triangular matrix with entries in $S(\mathfrak{gl}_n^{*})^{\mathfrak{gl}_n}$. For instance, for any $n$, we have $M_1=-P_1+mc_1$ and $M_2=P_2-c_1P_1+mc_2$.
\end{Remark}

\subsubsection{A general algorithm}
Now let us describe a procedure for computing a presentation of $\mathscr{B}(m\varpi_1)$.
\begin{itemize}\itemsep=0pt
\item[(1)]
Using relations \eqref{eq:girard_waring_1} for $N=m$ and $k>m$, express $P_{m+1},\dots,P_{m+n-1}$ via $P_1,\dots,P_m$. Plugging these expressions into \eqref{eq:medium_ring_presentation_2} gives relations for $c_1,\dots,c_n$, $P_1,\dots,P_m$ in $\mathscr{B}(m\varpi_1)$.

\item[(2)]
Write $P_{1},\dots,P_{m}$ in terms of $\theta_1,\dots,\theta_{m+1}$ and $M_1,\dots,M_{m}$ using \eqref{eq:M_P_transfer}. Then, apply \eqref{eq:girard_waring_2} for $N=n$ and $k=1,\dots,m+1$ to express $\theta_1,\dots,\theta_{m+1}$ as polynomials in $c_1,\dots,c_n$. This allows us to rewrite $P_1,\dots,P_m$ as polynomial expressions in $c_1,\dots,c_n$ and $M_1,\dots,M_m$.

\item[(3)]
Substituting expressions from (2) into relations which were obtained in (1) yields relations for $c_1,\dots,c_n$, $M_1,\dots,M_m$ in $\mathscr{B}(m\varpi_1)$.
\end{itemize}

In order to compare our computations with those of Hausel \cite[Section~4.2]{Hausel_24} and Rozhkovskaya \cite[Proposition~4.2]{Rozhkovskaya}, one has to calculate the $\mathfrak{sl}_n$-version\footnote{To avoid a confusion, we denote from now on in this subsection the two variants of big algebras by $\mathscr{B}_{\mathfrak{sl}_n}(m\varpi_1)$ and $\mathscr{B}_{\mathfrak{gl}_n}(m\varpi_1)$, respectively.} of $\mathscr{B}(m\varpi_1)$. It follows from a general construction that the big algebra $\mathscr{B}_{\mathfrak{sl}_n}(m\varpi_1)$ can be obtained from its $\mathfrak{gl}_n$-version $\mathscr{B}_{\mathfrak{gl}_n}(m\varpi_1)$ simply by quotienting out the \emph{trace}, i.e., the element $c_1\in S(\mathfrak{gl}_n^{*})^{\mathfrak{gl}_n}$ (see \cite[Section~2.1]{Hausel_24}). However, to match the presentations one has to take into account that the Kirillov--Wei operators for $\mathfrak{sl}_n$ and $\mathfrak{gl}_n$ are different. In particular, the medium operators $M_{1},\dots,M_{n-1}$ in $\mathscr{B}_{\mathfrak{gl}_n}(m\varpi_1)$ are replaced by the following operators in $\mathscr{B}_{\mathfrak{sl}_n}(m\varpi_1)\simeq\mathscr{B}_{\mathfrak{gl}_n}(m\varpi_1)/\langle c_1\rangle$:
\begin{equation*}
\widetilde{M}_i=M_i-\frac{m(n-i)}{n}c_i\pmod{c_1},\qquad i=1,\dots,n-1.
\end{equation*}
In the examples below, we give presentations of the algebra $\mathscr{B}_{\mathfrak{sl}_2}(m\varpi_1)$ in terms of $c_2,\dots,c_n$ and \smash{$\widetilde{M}_1,\dots,\widetilde{M}_{n-1}$}.

\subsubsection{Explicit form of relations} We apply the algorithm above to the following three cases: $(m,n)=(m,2)$, where $m$ is arbitrary, $(m,n)=(2,3)$ and $(m,n)=(3,3)$.

\begin{Example}[$n=2$ and $m$ is arbitrary] This example was essentially computed by Hausel and Rychlewicz, see \cite[Example 4.6]{Hausel_Rychlewicz} and \cite[Section~4.2.2]{Hausel_24}. More precisely, they computed the $\mathrm{SL}_2$-equivariant cohomology of the projective space $\mathbb{P}^m$ and the latter is known to be isomorphic to the medium algebra $\mathscr{M}_{\mathfrak{sl}_2}(m\varpi_1)$.

In this case, $\mathscr{B}_{\mathfrak{gl}_2}(m\varpi_1)$ is isomorphic to the quotient of the ring $\mathbb{C}[c_1,c_2][P_1,\dots,P_m]$ by the ideal with $m$ generators
\begin{equation*}
P_{2}-c_{1}P_{1}+c_{2}P_{0}, \ \dots,\ P_{m}-c_{1}P_{m-1}+c_{2}P_{m-2},\ P_{m+1}-c_{1}P_{m}+c_{2}P_{m-1},
\end{equation*}
where we set $P_0=m$ and \smash{$P_{m+1}=\Theta^{(m+1)}_{m}(P_1,\dots,P_m)$}.
Since $\mathscr{B}_{\mathfrak{sl}_2}(m\varpi_1)\simeq\mathscr{B}_{\mathfrak{gl}_2}(m\varpi_1)/\langle c_1\rangle$, we~obtain
\begin{gather*}
\mathscr{M}_{\mathfrak{sl}_2}(m\varpi_1)\simeq\mathbb{C}[c_2][P_1,\dots,P_m]/\langle P_2+mc_2,P_3+c_2P_1\dots,P_m+c_2P_{m-2},\\
\hphantom{\mathscr{M}_{\mathfrak{sl}_2}(m\varpi_1)\simeq\mathbb{C}[c_2][P_1,\dots,P_m]/\langle }{}\,
\Theta^{(m)}_{m+1}(P_1,\dots,P_m)+c_2P_{m-1}\rangle.
\end{gather*}
It follows from the proof of Theorem~\ref{thm:big_sym_power_gen_rel} and Corollary~\ref{cor:big_sym_power_gen_rel} that the equality ${P_{l+2}+c_2P_{l}=0}$ holds in the ring $\mathscr{M}_{\mathfrak{sl}_2}(m\varpi_1)$ for all $l$. Hence, $P_{2l}=(-1)^{l}mc_2^l$ and $P_{2l+1}=(-1)^lc_2^lP_1$, or equivalently,
\begin{equation*}
P_k=\Theta^{(m)}_{k}\bigl(P_1,-c_2,-c_2P_1,c_2^2,\dots\bigr)=
\begin{cases}
(-1)^lmc_2^l, &k=2l,~\text{where}~l\in\mathbb{Z}_{\ge 0},
\\
(-1)^lc_2^lP_1, &k=2l+1,~\text{where}~l\in\mathbb{Z}_{\ge 0}.
\end{cases}
\end{equation*}
In particular, we can eliminate $P_2,\dots,P_m$ from the list of generators for the $\mathbb{C}[c_2]$-algebra $\mathscr{M}_{\mathfrak{sl}_2}(m\varpi_1)$. We obtain that the medium algebra $\mathscr{M}_{\mathfrak{sl}_2}(m\varpi_1)$ is actually a quotient of $\mathbb{C}[c_2][P_1]$:
\begin{gather}\label{eq:medium_sl2_relation}
\mathscr{B}_{\mathfrak{sl}_2}(m\varpi_1)\simeq
\begin{cases}
\mathbb{C}[c_2][P_1]\big/\bigl\langle \Theta^{(m)}_{m+1}(P_1,-c_2,-c_2P_1,\dots)-(-1)^{m/2}c_2^{m/2}P_1\bigr\rangle, \\
\quad{}\text{if}~m~\text{is even},
\\
\mathbb{C}[c_2][P_1]\big/\bigl\langle \Theta^{(m)}_{m+1}(P_1,-c_2,-c_2P_1,\dots)-(-1)^{(m+1)/2}mc_2^{(m+1)/2}\bigr\rangle, \\
\quad{}\text{if}~m~\text{is odd}.
\end{cases}
\end{gather}
Alternatively, consider the sequence $\{\psi_n(a,b)\}_{n=1}^{\infty}$ of polynomials in two variables defined by the recurrence relation $\psi_{n+2}(a,b)=-a\psi_{n}(a,b)$ with initial conditions $\psi_{0}(a,b)=m$ and $\psi_1(a,b)=b$. Then, $P_l=\psi_l(c_2,P_1)$ for all $l$ and therefore,
\begin{equation*}
\mathscr{B}_{\mathfrak{sl}_2}(m\varpi_1)\simeq\mathbb{C}[c_2,P_1]\big/\bigl\langle\Theta^{(m)}_{m+1}(\psi_1(c_2,P_1),\dots,\psi_m(c_2,P_1))-\psi_{m+1}(c_2,P_1)\bigr\rangle.
\end{equation*}
The recurrence relation imply the following explicit formula for $\psi_k(a,b)$:
\begin{equation}\label{eq:psi_polynomial_formula}
\psi_k(a,b)=\biggl(\frac{m}{2}+\frac{b}{2i\sqrt{a}}\biggr)\bigl(i\sqrt{a}\bigr)^{k}+\biggl(\frac{m}{2}-\frac{b}{2i\sqrt{a}}\biggr)\bigl(-i\sqrt{a}\bigr)^{k}.
\end{equation}

Now let us calculate explicitly the relation for $P_1$ and $c_2$ in $\mathscr{B}_{\mathfrak{sl}_2}(m\varpi_1)$. Consider the power series
\begin{equation}\label{eq:power_sum_z_series}
\exp\biggl\{-p_1z-\frac{1}{2}p_2z^2-\frac{1}{3}p_3z^3-\cdots\biggr\}=\exp\Biggl\{-\sum_{k=1}^{\infty}\frac{p_k}{k}z^k\Biggr\},
\end{equation}
where $p_k=p_k(u_1,\dots,u_m)$ is the $k$-th power sum in variables $u_1,\dots,u_m$. Observe that the $z^{m+1}$ term of this power series equals
\begin{equation*}
\sum_{\substack{i_1,\dots,i_{m+1}\ge 0\\i_1+2i_2+\dots+(m+1)i_{m+1}=m+1}}\frac{(-1)^{i_1+i_2+\dots+i_{m+1}}}{i_1!\cdots i_{m+1}!}\prod_{j=1}^{k}\biggl(\frac{1}{j}p_{j}\biggr)^{i_j}=\Theta_{m+1}^{(m)}(p_1,\dots,p_m)-p_{m+1},
\end{equation*}
where we used \eqref{eq:girard_waring_1} for $N=m$ and $k=m+1$. This identity together with \eqref{eq:medium_sl2_relation} implies that the relation for $P_1$ and $c_2$ in $\mathscr{B}_{\mathfrak{sl}_2}(m\varpi_1)$ is equal to the $z^{m+1}$ term of the power series \eqref{eq:power_sum_z_series} after the substitution $(p_1,p_2,\dots,p_m,\dots)=(\psi_1(c_2,P_1),\psi_2(c_2,P_1),\dots,\psi_m(c_2,P_1),\dots)$. To compute the resulting power series, note that by \eqref{eq:psi_polynomial_formula} we have
\begin{equation*}
\psi_k(c_2,P_1)=\alpha(i\sqrt{c_2})^{k}+\beta(-i\sqrt{c_2})^{k},\qquad \text{where}~\alpha=\frac{m}{2}+\frac{P_1}{2i\sqrt{c_2}},\ \beta=\frac{m}{2}-\frac{P_1}{2i\sqrt{c_2}}.
\end{equation*}
Then, the substitution gives
\begin{equation*}
\exp\Biggl\{-\sum_{k=1}^{\infty}\frac{p_k}{k}z^k\Biggr\}\Bigg|_{p_j=\psi_j(c_2,P_1) } =\exp\Biggl\{-\sum_{k=1}^{\infty}\frac{\alpha(i\sqrt{c_2})^k+\beta(-i\sqrt{c_2})^k}{k}z^k\Biggr\}.
\end{equation*}
The discussion above implies that the remaining relation between $P_1$ and $c_2$ is given by the $z^{m+1}$ term of the power series above. A direct calculation shows that
\begin{gather*}
\exp\Biggl\{-\sum_{k=1}^{\infty}\frac{\alpha(i\sqrt{c_2})^k+\beta(-i\sqrt{c_2})^k}{k}z^k\Biggr\}=(1-i\sqrt{c_2}z)^{\alpha}(1+i\sqrt{c_2}z)^{\beta}
\\
\qquad{}=\sum_{k,l\ge 0}(-1)^k\binom{\alpha}{k}\binom{\beta}{l}(i\sqrt{c_2})^{k+l}z^{k+l}=\sum_{k,l\ge 0}(-1)^{k+l}\binom{\alpha}{k}\binom{-\beta+l-1}{l}(i\sqrt{c_2})^{k+l}z^{k+l}
\\
\qquad{}=\sum_{k,l\ge 0}\frac{(-1)^{k+l}z^{k+l}}{2^{k+l}\cdot k!\cdot l!}\prod_{j=0}^{k-1}(P_1+(m-2j)i\sqrt{c_2})\prod_{j=0}^{l-1}(P_1-(m-2j)i\sqrt{c_2}).
\end{gather*}
The $z^{m+1}$ term of this expression equals
\begin{equation*}
\frac{(-1)^{m+1}}{2^{m+1}}\sum_{k+l=m+1}\frac{1}{k!\cdot l!}\prod_{j=0}^{k-1}(P_1+(m-2j)i\sqrt{c_2})\prod_{j=0}^{l-1}(P_1-(m-2j)i\sqrt{c_2}).
\end{equation*}
However, for $k+l=m+1$, we have
\begin{equation*}
\prod_{j=0}^{k-1}(P_1+(m-2j)i\sqrt{c_2})\prod_{j=0}^{l-1}(P_1-(m-2j)i\sqrt{c_2})=\prod_{j=0}^{m}(P_1+(m-2j) i\sqrt{c_2}),
\end{equation*}
so the sum above equals
\begin{equation*}
\frac{(-1)^{m+1}}{2^{m+1}}\biggl(\sum_{k+l=m+1}\frac{1}{k!\cdot l!}\biggr)\prod_{j=0}^{l-1}(P_1-(m-2j)i\sqrt{c_2})=\frac{(-1)^{m+1}}{(m+1)!}\prod_{j=0}^{m}(P_1+(m-2j)i\sqrt{c_2}).
\end{equation*}
Therefore, the relation between $P_1$ and $c_2$ in $\mathscr{B}_{\mathfrak{sl}_2}(m\varpi_1)$ is just $\prod_{j=0}^{m}(P_1+(m-2j)i\sqrt{c_2})=0$. Finally, recall that \smash{$\widetilde{M}_1=M_1-\frac{m}{2}c_1=\frac{m}{2}c_1-P_1=-P_1\pmod {c_1}$} in $\mathscr{B}_{\mathfrak{sl}_2}(m\varpi_1)$.

Our computations can be summarized in the following proposition.
\begin{Proposition}
We have the following presentations of the big algebra $\mathscr{B}_{\mathfrak{sl}_2}(m\varpi_1)$:
\begin{gather*}
\mathscr{B}_{\mathfrak{sl}_2}(m\varpi_1)\simeq\frac{\mathbb{C}[t_1,t_2]^{\mathfrak{S}_2}[u_1,\dots,u_m]^{\mathfrak{S}_m}}{\bigl\langle t_1+t_2,\bigl\{u_1^{r}(u_1-t_1)(u_1-t_2)+\dots+u_m^r(u_m-t_1)(u_m-t_2)\bigr\}_{r=0}^{m-1}\bigr\rangle},\quad \text{and}
\\
\mathscr{B}_{\mathfrak{sl}_2}(m\varpi_1)\simeq
\begin{cases}
\mathbb{C}[c_2]\bigl[\widetilde{M}_1\bigr]\big/\bigl\langle \bigl(\widetilde{M}_1^2+m^2c_2\bigr)\bigl(\widetilde{M}_1^2+(m-2)^2c_2\bigr)\cdots\bigl(\widetilde{M}_1^2+4c_2\bigr)\widetilde{M}_1\bigr\rangle,\\
\quad{}\text{if}~m~\text{is even},
\\
\mathbb{C}[c_2]\bigl[\widetilde{M}_1\bigr]\big/\bigl\langle\bigl(\widetilde{M}_1^2+m^2c_2\bigr)\bigl(\widetilde{M}_1^2+(m-2)^2c_2\bigr)\cdots\bigl(\widetilde{M}_1^2+c_2\bigr)\bigr\rangle,\\
\quad{}\text{if}~m~\text{is odd}.
\end{cases}
\end{gather*}
These two presentations are related by the following change of generators:
\begin{equation*}
\widetilde{M}_1=-(u_1+\dots+u_m),\qquad c_2=t_1t_2.
\end{equation*}
\end{Proposition}

A similar calculation gives an analogous presentation of $\mathscr{B}_{\mathfrak{gl}_2}(m\varpi_1)$.

\begin{Proposition}
We have the following presentations of the big algebra $\mathscr{B}_{\mathfrak{gl}_2}(m\varpi_1)$:
\begin{gather*}
\mathscr{B}_{\mathfrak{gl}_2}(m\varpi_1)\simeq\frac{\mathbb{C}[t_1,t_2]^{\mathfrak{S}_2}[u_1,\dots,u_m]^{\mathfrak{S}_m}}{\bigl\langle \{u_1^{r}(u_1-t_1)(u_1-t_2)+\dots+u_m^r(u_m-t_1)(u_m-t_2)\}_{r=0}^{m-1}\bigr\rangle}, \qquad \text{and}
\\
\mathscr{B}_{\mathfrak{gl}_2}(m\varpi_1)\simeq\mathbb{C}[c_1,c_2][M_1]/I_{m\varpi_1},
\end{gather*}
where
\begin{gather*}
I_{m\varpi_1}=
\begin{cases}
\biggl\langle\! \bigl(M_1-\frac{m}{2}c_1\bigr)\prod\limits_{j=0}^{m/2}\bigl(M_1^2-mc_1M_1+j(m-j)c_1^2+(m-2j)^2c_2\bigr)\!\biggr\rangle,&~\text{if}~m~\text{is even},
\\
\biggl\langle\prod\limits_{j=0}^{(m-1)/2}\bigl(M_1^2-mc_1M_1+j(m-j)c_1^2+(m-2j)^2c_2\bigr)\!\biggr\rangle,&~\text{if}~m~\text{is odd}.
\end{cases}
\end{gather*}
These two presentations are related by the following change of generators:
\begin{equation*}
M_1=\frac{m}{2}(t_1+t_2)-(u_1+\dots+u_m),\qquad c_1=t_1+t_2,\qquad c_2=t_1t_2.
\end{equation*}
\end{Proposition}
\end{Example}

\begin{Example}[$n=3$ and $m=2$]
This example was computed\footnote{Strictly speaking, she computed a presentation of the \emph{Kirillov algebra} $\mathscr{C}_{\mathfrak{sl_3}}(2\varpi_1)$, but in this case one can check that the Kirillov, medium and big algebras all coincide (cf.\ Corollary~\ref{cor:wmf_medium_equal_kirillov}).} by Rozhkovskaya, see \cite[Proposition 4.2]{Rozhkovskaya}. Corollary~\ref{cor:big_sym_power_gen_rel} gives the following presentation of the $\mathscr{B}_{\mathfrak{gl}_3}(2\varpi_1)$:
\begin{equation*}
\mathscr{B}_{\mathfrak{gl}_3}(2\varpi_1)\simeq\mathbb{C}[c_1,c_2,c_3][P_1,P_2]\bigg/\!\!
\left\langle
\begin{matrix}
P_3-c_1P_2+c_2P_1-c_3P_0,
\\
P_4-c_1P_3+c_2P_2-c_3P_1
\end{matrix}
\right\rangle,
\end{equation*}
where we substitute $P_0=2$, \smash{$P_3=-\frac{1}{2}P_1^3+\frac{3}{2}P_1P_2$} and \smash{$P_4=-\frac{1}{2}P_1^4+P_1^2P_2$}. This yields the following presentation:
\begin{gather*}
\mathscr{B}_{\mathfrak{gl}_3}(2\varpi_1)\nonumber\\
\qquad{}
\simeq\mathbb{C}[c_1,c_2,c_3][P_1,P_2]\bigg/\!\!\left\langle
\begin{matrix}
P_1^3-3P_1P_2-2c_2P_1+2c_1P_2+4c_3,
\\
P_1^4-c_1P_1^3-2P_1^2P_2+3c_1P_1P_2+2c_3P_1-P_2^2-2c_2P_2
\end{matrix}
\right\rangle.
\end{gather*}
Using the fact that $M_1=2c_1-P_1$ and $M_2=P_2-c_1P_1+2c_2$, one can obtain a presentation of~$\mathscr{B}_{\mathfrak{gl}_2}(m\varpi_1)$ in terms of the medium operators:
\begin{equation*}
\mathscr{B}_{\mathfrak{gl}_3}(2\varpi_1)
\simeq\mathbb{C}[c_1,c_2,c_3][M_1,M_2]/I_{2\varpi_1},
\end{equation*}
where the ideal $I_{2\varpi_1}$ is given by
\begin{equation*}
\begin{aligned}
I_{2\varpi_1}
&=
\left\langle
\begin{matrix}
M_1^3-3c_1M_1^2-3M_1M_2+2c_1^2M_1+4c_2M_1+4c_1M_2-4c_1c_2-4c_3,
\\
M_1^4-5c_1M_1^3-2M_1^2M_2+8c_1^2M_1^2+4c_2M_1^2+7c_1M_1M_2-4c_1^3M_1
\\
-\,12c_1c_2M_1-2c_3M_1-M_2^2-6c_1^2M_2+2c_2M_2+8c_1^2c_2+4c_1c_3
\end{matrix}
\right\rangle
\\
&=
\left\langle
\begin{matrix}
M_1^3-3c_1M_1^2-3M_1M_2+2c_1^2M_1+4c_2M_1+4c_1M_2-4c_1c_2-4c_3,
\\
M_1^2M_2-3c_1M_1M_2+2c_3M_1-M_2^2+2c_1^2M_2+2c_2M_2-4c_1c_3
\end{matrix}
\right\rangle.
\end{aligned}
\end{equation*}
In order to get a presentation of the $\mathfrak{sl}_2$-version, we use
\begin{gather*}
\widetilde{M}_1\equiv M_1-\frac{4}{3}c_1\equiv-P_1+\frac{2}{3}c_1\equiv-P_1\pmod{c_1},
\\
\widetilde{M}_2\equiv M_2-\frac{2}{3}c_2\equiv P_2-c_1P_1+\frac{4}{3}c_2\equiv P_2+\frac{4}{3}c_2\pmod{c_1}.
\end{gather*}
In particular, \smash{$P_1\equiv-\widetilde{M}_1\pmod{c_1}$} and \smash{$P_2\equiv\widetilde{M}_2-\frac{4}{3}c_2\pmod{c_1}$}. Plugging this into the presentation of $\mathscr{B}_{\mathfrak{gl}_3}(2\varpi_1)$ above and using the fact that $\mathscr{B}_{\mathfrak{sl}_3}(2\varpi_1)\simeq\mathscr{B}_{\mathfrak{sl}_3}(2\varpi_1)/\langle c_1\rangle$, we obtain
\begin{gather}\nonumber
\mathscr{B}_{\mathfrak{sl}_3}(2\varpi_1) \\ \nonumber
\quad{}\simeq\mathbb{C}[c_2,c_3]\bigl[\widetilde{M}_1,\widetilde{M}_2\bigr]\bigg/\!\!
\left\langle
\begin{matrix}
\widetilde{M}_1^3-3\widetilde{M}_1\widetilde{M}_2+2c_2\widetilde{M}_1-4c_3,
\\
9\widetilde{M}_1^4-18\widetilde{M}_1^2\widetilde{M}_2+24c_2\widetilde{M}_1^2-18c_3\widetilde{M}_1-9\widetilde{M}_2^2+6c_2\widetilde{M}_2+8c_2^2
\end{matrix}
\right\rangle
\\
\quad{}\simeq\mathbb{C}[c_2,c_3]\bigl[\widetilde{M}_1,\widetilde{M}_2\bigr]\bigg/\!\!
\left\langle
\begin{matrix}
\widetilde{M}_1^3-3\widetilde{M}_1\widetilde{M}_2+2c_2\widetilde{M}_1-4c_3,
\\
9\widetilde{M}_1^2\widetilde{M}_2+6c_2\widetilde{M}_1^2+18c_3\widetilde{M}_1-9\widetilde{M}_2^2+6c_2\widetilde{M}_2+8c_2^2
\label{eq:big_sl3_2omega_1_presentation}
\end{matrix}
\right\rangle.
\end{gather}
To obtain Rozhkovskaya's presentation, we introduce the following elements:
\begin{equation*}
M=-\widetilde{M}_1,
\qquad N=-\frac{1}{12}\widetilde{M}_1^2+\frac{1}{4}\widetilde{M}_2-\frac{1}{6}c_2,
\qquad C_2=-\frac{1}{3}c_2,
\qquad C_3=\frac{1}{3}c_3,
\end{equation*}
which is equivalent to
\begin{equation*}
\widetilde{M}_1=-M,
\qquad \widetilde{M}_2=\frac{1}{3}M^2+4N-2C_2,
\qquad c_2=-3C_2,
\qquad c_3=3C_3.
\end{equation*}
Substituting this into \eqref{eq:big_sl3_2omega_1_presentation} gives the following presentation:
\begin{gather*}
\mathscr{B}_{\mathfrak{sl}_3}(2\varpi_1)\\
\qquad\simeq\mathbb{C}[C_2,C_3][M,N]\bigg/\!\!
\left\langle
\begin{matrix}
MN-C_3,
\\
M^4+6M^2N-15C_2M^2-27C_3M-72N^2+36C_2N+36C_2^2
\end{matrix}
\right\rangle
\\
\qquad\simeq\mathbb{C}[C_2,C_3][M,N]\bigg/\!\!
\left\langle
\begin{matrix}
MN-C_3,
\\
M^4-15C_2M^2-21C_3M-72N^2+36C_2N+36C_2^2
\end{matrix}
\right\rangle.
\end{gather*}
Therefore, we obtain the generators and relations given in \cite[Proposition 4.2]{Rozhkovskaya}.
\end{Example}

\begin{Example}[$n=3$ and $m=3$] This example was computed by Hausel, see \cite[Section~4.2.3]{Hausel_24}.
In this case Corollary~\ref{cor:big_sym_power_gen_rel} gives the following of $\mathscr{B}_{\mathfrak{gl}_3}(3\varpi_1)$:
\begin{equation*}
\mathscr{B}_{\mathfrak{gl}_3}(3\varpi_1)\simeq\mathbb{C}[c_1,c_2,c_3][P_1,P_2,P_3]\Bigg/\!\!\left\langle
\begin{matrix}
P_3-c_1P_2+c_2P_1-c_3P_0,
\\
P_4-c_1P_3+c_2P_2-c_3P_1,
\\
P_5-c_1P_4+c_2P_3-c_3P_2
\end{matrix}
\right\rangle,
\end{equation*}
where we substitute $P_0=3$, $P_4=\frac{1}{6}P_1^4-P_1^2P_2+\frac{1}{2}P_2^2+\frac{4}{3}P_1P_3$ and $P_5=\frac{1}{6}P_1^5-\frac{5}{6}P_1^3P_2+\frac{5}{6}P_1^2P_3+\frac{5}{6}P_2P_3$. We obtain the following presentation:
\begin{equation*}
\mathscr{B}_{\mathfrak{gl}_3}(3\varpi_1)
\simeq\mathbb{C}[c_1,c_2,c_3][P_1,P_2,P_3]/J_{3\varpi_1},
\end{equation*}
where the ideal $J_{3\varpi_1}$ is given by
\begin{equation*}
\begin{aligned}
J_{3\varpi_1}
&=\scalebox{1.4}{\Bigg\langle}
\begin{matrix}
c_2P_1\!-\!c_1P_2\!+\!P_3\!-\!3c_3,
\\
P_1^4\!-\!6P_1^2P_2\!+\!8P_1P_3\!-\!6c_3P_1\!+\!3P_2^2\!+\!6c_2P_2\!-\!6c_1P_3,
\\
P_1^5\!-\!c_1P_1^4\!-\!5P_1^3P_2\!+\!6c_1P_1^2P_2\!+\!5P_1^2P_3\!-\!
8c_1P_1P_3\!-\!3c_1P_2^2\!+\!5P_2P_3\!-\!6c_3P_2\!+\!6c_2P_3
\end{matrix}
\scalebox{1.4}{\Bigg\rangle}
\\
&=\scalebox{1.4}{\Bigg\langle}
\begin{matrix}
P_1^4\!-\!6P_1^2P_2\!-\!8c_2P_1^2\!+\!8c_1P_1P_2\!+\!6c_1c_2P_1\!+\!18c_3P_1\!+\!3P_2^2\!-\!6c_1^2P_2\!+\!
6c_2P_2\!-\!18c_1c_3,
\\
P_1^5\!-\!c_1P_1^4\!-\!5P_1^3P_2\!-\!5c_2P_1^3\!+\!11c_1P_1^2P_2\!+\!8c_1c_2P_1^2\!+\!15c_3P_1^2\!-\!8c_1^2P_1P_2
\\
+\,2c_1P_2^2\!-\!5c_2P_1P_2\!-\!6c_2^2P_1\!-\!24c_1c_3P_1\!+\!6c_1c_2P_2\!+\!9c_3P_2\!+\!18c_2c_3
\end{matrix}
\scalebox{1.4}{\Bigg\rangle}.
\end{aligned}
\end{equation*}
As $M_1=3c_1-P_1$ and $M_2=P_2-c_1p_1+3c_2$, we have
\begin{equation*}
\mathscr{B}_{\mathfrak{gl}_3}(3\varpi_1)
\simeq\mathbb{C}[c_1,c_2,c_3][M_1,M_2]/I_{3\varpi_1},
\end{equation*}
where the ideal $I_{3\varpi_1}$ is given by
\begin{gather*}
I_{3\varpi_1}
\!=\!\scalebox{1.8}{\Bigg\langle}\!\!
\begin{matrix}
M_1^4\!-\!6c_1M_1^3\!-\!6M_1^2M_2\!+\!11c_1^2M_1^2\!+\!10c_2M_1^2\!+\!22c_1M_1M_2\!-\!6c_1^3M_1\!-\!30c_1c_2M_1
\\
 -\,18c_3M_1\!+\!3M_2^2\!-\!18c_1^2M_2\!-\!12c_2M_2\!+\!18c_1^2c_2\!+\!36c_1c_3\!+\!9c_2^2,
\\
M_1^5\!-\!9c_1M_1^4\!-\!5M_1^3M_2\!+\!29c_1^2M_1^3\!+\!10c_2M_1^3\!+\!34c_1M_1^2M_2\!-\!39c_1^3M_1^2\!-\!60c_1c_2M_1^2
\\
-\,15c_3M_1^2\!-\!73c_1^2M_1M_2\!-\!5c_2M_1M_2\!+\!18c_1^4M_1\!+\!108c_1^2c_2M_1\!+\!75c_1c_3M_1\!+\!9c_2^2M_1
\\
 -\,2c_1M_2^2\!+\!48c_1^3M_2\!+\!21c_1c_2M_2\!-\!9c_3M_2\!-\!54c_1^3c_2\!-\!90c_1^2c_3\!-\!27c_1c_2^2\!+\!9c_2c_3
\end{matrix}\!\!
\scalebox{1.8}{\Bigg\rangle}.\!\!
\end{gather*}
The medium operators for the $\mathscr{B}_{\mathfrak{sl}_3}(3\varpi_1)$ are given by the following formulas:
\begin{gather*}
\widetilde{M}_1\equiv M_1-2c_1\equiv -P_1+c_1=-P_1\pmod{c_1},
\\
\widetilde{M}_2\equiv M_2-c_2\equiv P_2-c_1P_1+2c_2\equiv P_2+2c_2\pmod{c_1}.
\end{gather*}
Hence, \smash{$P_1=-\widetilde{M}_1\pmod{c_1}$} and \smash{$P_2=\widetilde{M}_2-2c_2\pmod{c_1}$}. Plugging this into the presentation above and using the isomorphism $\mathscr{B}_{\mathfrak{sl}_3}(3\varpi_1)\simeq \mathscr{B}_{\mathfrak{gl}_3}(3\varpi_1)/\langle c_1\rangle$, we get
\begin{equation*}
\mathscr{B}_{\mathfrak{sl}_3}(3\varpi_1)=\frac{\mathbb{C}[c_2,c_3]\bigl[\widetilde{M}_1,\widetilde{M}_2\bigr]}{
\left\langle
\begin{matrix}
\widetilde{M}_1^4-6\widetilde{M}_1^2\widetilde{M}_2+4c_2\widetilde{M}_1^2-18c_3\widetilde{M}_1+3\widetilde{M}_2^2-6c_2\widetilde{M}_2,
\\
\widetilde{M}_1^5-5\widetilde{M}_1^3\widetilde{M}_2+5c_2\widetilde{M}_1^3-15c_3\widetilde{M}_1^2-5c_2\widetilde{M}_1\widetilde{M}_2+4c_2^2\widetilde{M}_1-9c_3\widetilde{M}_2
\end{matrix}
\right\rangle}.
\end{equation*}
Subtracting from the second relation the first one multiplied by $\widetilde{M}_1$, we obtain
\begin{equation*}
\mathscr{B}_{\mathfrak{sl}_3}(3\varpi_1)=\frac{\mathbb{C}[c_2,c_3]\bigl[\widetilde{M}_1,\widetilde{M}_2\bigr]}{
\left\langle
\begin{matrix}
\widetilde{M}_1^4-6\widetilde{M}_1^2\widetilde{M}_2+4c_2\widetilde{M}_1^2-18c_3\widetilde{M}_1+3\widetilde{M}_2^2-6c_2\widetilde{M}_2,
\\
\widetilde{M}_1^3\widetilde{M}_2+c_2\widetilde{M}_1^3+3c_3\widetilde{M}_1^2-3\widetilde{M}_1\widetilde{M}_2^2+c_2\widetilde{M}_1\widetilde{M}_2+4c_2^2\widetilde{M}_1-9c_3\widetilde{M}_2
\end{matrix}
\right\rangle},
\end{equation*}
which coincides with the presentation given in \cite[equation~(4.6)]{Hausel_24}.
\end{Example}

\subsection{Formula for the Kirillov--Wei operator}\label{subsect:kirillov_wei_sym_power}

The aim of this subsection is to derive a more explicit formula for the Kirillov--Wei operator on the Kirillov algebra of $V(m\varpi_1)\simeq S^m(\mathbb{C}^n)$.

\begin{Lemma}\label{lemma:leibniz_rule_sym_part}
For any $A\in S(\mathfrak{gl}_{n}^{*})^{\mathfrak{gl}_n}$ and any element $B$ of the Kirillov algebra $\mathscr{C}(V)$, we have the identity
\begin{equation*}
\DKir(A\cdot B)=\DKir(A)\cdot B+A\cdot \DKir(B).
\end{equation*}
\end{Lemma}
\begin{proof}
This is a direct consequence of the definition of $\DKir$ and the usual Leibniz rule, see formula~\eqref{eq:kirillov_wei_op_def}.
\end{proof}
\begin{Remark}
Note however that $\DKir$ does \textit{not} satisfy the Leibniz rule for arbitrary ${A,B\in\mathscr{C}(V)}$. It is important here that $A$ is a scalar matrix.
\end{Remark}

We use the map $\mathfrak{i}_m$ in order to obtain a different interpretation for the Kirillov--Wei operator on $\mathscr{B}(m\varpi_1)$.\footnote{If one accepts Conjecture~\ref{conj:kirillov_alg_center} for $\lambda=m\varpi_1$, then $\mathscr{B}(m\varpi_1)=\mathscr{C}(m\varpi_1)$ and hence $\DKir$ indeed acts on $\mathscr{B}(m\varpi_1)$. In our case, one can also check this directly, see Remark~\ref{rem:D_op_inv_big_alg}.}
Namely, define the operator $\widehat{\DKir}\colon\mathfrak{i}_{m}(\mathscr{B}(m\varpi_1))\to\mathfrak{i}_{m}(\mathscr{B}(m\varpi_1))$ so that the diagram%
\begin{equation*}
\xymatrix{
\mathscr{B}(m\varpi_1) \ar[r]^{\DKir} \ar[d]_{\mathfrak{i}_{m}} & \mathscr{B}(m\varpi_1) \ar[d]^{\mathfrak{i}_{m}}
\\
\mathfrak{i}_{m}(\mathscr{B}(m\varpi_1)) \ar[r]^{\widehat{\DKir}} & \mathfrak{i}_{m}(\mathscr{B}(m\varpi_1))
}
\end{equation*}
commutes. It turns out that $\widehat{\DKir}$ can be described using some constructions from symmetric function theory.

\begin{Proposition}\label{prop:kirillov_wei_operator_formula}
The Kirillov--Wei operator $\DKir$ on $\mathscr{B}(m\varpi_1)$ induces on $\mathfrak{i}_{m}(\mathscr{B}(m\varpi_1))\subset\Func(\wt(m\varpi_1),S(\mathfrak{h}^{*}))$ the operator $\widehat{\DKir}$ which acts as follows: for any symmetric functions $f$ and $g$ it maps the function $\mu\mapsto f[t_1+\dots+t_n]g[\mu_1t_1+\dots+\mu_nt_n]$ to the function
\begin{gather}
\mu\mapsto
\biggl(\mu_1\frac{\partial}{\partial t_1}+\dots+\mu_n\frac{\partial}{\partial t_n}\biggr)(f[t_1+\dots+t_n]g[\mu_1t_1+\dots+\mu_nt_n])
\nonumber \\ \hphantom{\mu\mapsto}{}
+f[t_1+\dots+t_n]
\nonumber \\ \hphantom{\mu\mapsto+}{}
\times\sum_{\substack{1\le i,j\le n\\i\neq j}}\frac{\mu_j(\mu_i+1)}{t_i-t_j}(g[\mu_1t_1+\dots+\mu_nt_n+(t_i-t_j)]-g[\mu_1t_1+\dots\mu_nt_n]).\label{eq:kirillov_wei_sym_func}
\end{gather}
Here, $F[-]$ for a symmetric function $F$ stands for the plethystic substitution $($see Remark~$\ref{remark:plethystic_subst}$ and {\rm \cite[{\it Section}~I.8]{Macdonald})}.
\end{Proposition}
\begin{Remark}
It is suggestive to regard the ``shifts'' $t_i-t_j$ as corresponding to the elements of the $A_{n-1}$ root system.
\end{Remark}

\begin{proof}[Proof of Proposition~\ref{prop:kirillov_wei_operator_formula}]
Throughout the proof, we perform computations in the big algebra $\mathscr{B}(\mathcal{P}(n,1))$. To obtain formulas for $\mathscr{B}(m\varpi_1)$, one simply restricts all operators to the subspace $V(m\varpi_1)\simeq\mathbb{C}[x_1,\dots,x_n]_{m}\subset\mathcal{P}(n,1)$.

Since the operators $\DKir$ and $\widehat{\DKir}$ are additive, it suffices to verify the statement from the proposition only for homogeneous $f$ and $g$. Applying Lemma~\ref{lemma:leibniz_rule_sym_part}, we can reduce\footnote{Note that if we denote the right-hand side of \eqref{eq:kirillov_wei_sym_func} by $D_{\mu}(f,g)$, then $D_{\mu}(f,g)$ satisfies an analogue of Lemma~\ref{lemma:leibniz_rule_sym_part}. Namely, $D_{\mu}(f,g)=f[t_1+\dots+t_n]D_{\mu}(1,g)+D_{\mu}(f,1)g[\mu_1t_1+\dots+\mu_nt_n]$.} this to the two cases: $f\equiv 1$ and $g\equiv 1$.

\textit{Case 1: assume that $g\equiv 1$.}
It is well known that the ring $S(\mathfrak{gl}_n^{*})^{\mathfrak{gl}_n}$ is a free polynomial ring in the elements $\theta_1,\dots,\theta_n\in S(\mathfrak{gl}_n^{*})$, where $\theta_k(Y)=\tr\bigl(Y^k\bigr)$, $Y\in\mathfrak{gl}_n$. Now let us show that $\DKir(\theta_{\alpha})=\alpha P_{\alpha-1}$ for all $\alpha=1,\dots,n$. Indeed, we have
\begin{align*}
\DKir(\theta_{\alpha})(Y)&{}=\sum_{i,j=1}^{n}\frac{\partial (\tr(Y^{\alpha}))}{\partial y_{ji}}\cdot L(E_{ij})=\sum_{i,j=1}^{n}\sum_{l=1}^{n}\sum_{\beta=1}^{\alpha}\bigl[Y^{\beta-1}\bigr]_{lj}\bigl[Y^{\alpha-\beta}\bigr]_{il}\cdot L(E_{ij})
\\
&{}=\sum_{i,j=1}^{n}\sum_{\beta=1}^{\alpha}\bigl[Y^{\beta-1}\cdot Y^{\alpha-\beta}\bigr]_{ij}\cdot L(E_{ij})=\alpha\sum_{i,j=1}^{n}\bigl[Y^{\alpha-1}\bigr]\cdot L(E_{ij})\\
&{}=\alpha L\bigl(Y^{\alpha-1}\bigr)=\alpha P_{\alpha-1}(Y).
\end{align*}
Now note that for $Y=\diag(t_1,\dots,t_n)\in\mathfrak{h}$ we have $\theta_{\alpha}(Y)=t_1^{\alpha}+\dots+t_n^{\alpha}$.
Therefore, $\widehat{\DKir}$ maps a~constant function \smash{$\mu\mapsto\sum_{i=1}^{n}t_i^{\alpha}$} to the function
\[
\mu\mapsto\sum_{i=1}^{n}\alpha\mu_i t_i^{\alpha-1}=\biggl(\sum_{i}\mu_i\frac{\partial}{\partial t_i}\biggr)\biggl(\sum_{i=1}^{n}t_i^{\alpha}\biggr).
\]
 This together with Lemma~\ref{lemma:leibniz_rule_sym_part} implies that for any \smash{$f\in\mathbb{C}[t_1,\dots,t_n]^{\mathfrak{S}_{n}}$} the operator $\widehat{\DKir}$ maps any constant function $\mu\mapsto f$ to \smash{$\mu\mapsto\bigl(\sum_{i}\mu_i\frac{\partial}{\partial t_i}\bigr)(f)$}.

\textit{Case 2: assume that $f\equiv 1$.}
As \smash{$\widehat{\DKir}$} is $\mathbb{C}$-linear we may assume without loss of generality that \smash{$g=\prod_{k=1}^{l}p_{\alpha_k}$} for some positive integers $\alpha_1,\dots,\alpha_l$, where $p_{\alpha}$ is the $\alpha$-th power sum. In~other words, we consider the function
\begin{equation*}
\mu\mapsto\prod_{k=1}^{l}\Biggl(\sum_{i=1}^{n}\mu_{i}t_{i}^{\alpha_{k}}\Biggr)=\prod_{k=1}^{l}p_{\alpha_k}[\mu_1t_1+\dots+\mu_nt_n].
\end{equation*}
Note that this element of $\mathfrak{i}_{m}(\mathscr{B}(m\varpi_1))$ corresponds to the element $P_{\alpha_1}P_{\alpha_2}\cdots P_{\alpha_l}$ in $\mathscr{B}(m\varpi_1)$. Then, we can compute the action of the Kirillov--Wei operator on this element:
\begin{gather*}
\DKir(P_{\alpha_1}\cdots P_{\alpha_l})(Y)
\\
\qquad{}=\sum_{i_1,j_1,\dots,i_l,j_l=1}^{n}\sum_{i,j=1}^{n}\frac{\partial}{\partial y_{ji}}([Y^{\alpha_1}]_{i_1j_1}\cdots[Y^{\alpha_l}]_{i_lj_l})\cdot L(E_{i_1j_1})\cdots L(E_{i_lj_l})L(E_{ij})
\\
\qquad{}=\sum_{i_1,j_1,\dots,i_l,j_l=1}^{n}\sum_{i,j=1}^{n} \Biggl(\sum_{k=1}^{l}\sum_{\alpha=0}^{\alpha_k-1}[Y^{\alpha_1}]_{i_1j_1}\cdots[Y^{\alpha}]_{i_kj}\bigl[Y^{\alpha_k-\alpha-1}\bigr]_{ij_k}\cdots[Y^{\alpha_l}]_{i_lj_l}\Biggr)
\\
\qquad{}\hphantom{=\sum_{i_1,j_1,\dots,i_l,j_l=1}^{n}\sum_{i,j=1}^{n}}{}\
\times L(E_{i_1j_1})\cdots L(E_{i_lj_l})L(E_{ij})
\\
\qquad{}=\sum_{k=1}^{l}L(Y^{\alpha_1})\cdots L(Y^{\alpha_{k-1}})\sum_{i_{k},j_{k},\dots,i_l,j_l=1}^{n}\sum_{i,j=1}^{n}\Biggl(\sum_{\alpha=0}^{\alpha_k-1}[Y^{\alpha}]_{i_kj}\bigl[Y^{\alpha_k-\alpha-1}\bigr]_{ij_k}\Biggr)\\
\qquad\hphantom{=\sum_{k=1}^{l}}{}\
\times [Y^{\alpha_{k+1}}]_{i_{k+1}j_{k+1}}\cdots[Y^{\alpha_l}]_{i_lj_l}L(E_{i_kj_k})L(E_{i_{k+1}j_{k+1}})\cdots L(E_{i_lj_l})L(E_{ij}).
\end{gather*}
For $Y=\diag(t_1,\dots,t_n)\in\mathfrak{h}$, this simplifies as follows:
\begin{gather*}
\DKir(P_{\alpha_1}\cdots P_{\alpha_l})(Y)
\\
\quad=\sum_{k=1}^{l}\Biggl(\prod_{s=1}^{k-1}P_{\alpha_{s}}(Y)\Biggr)\sum_{i_{k+1},\dots,i_l=1}^{n} \sum_{i,j=1}^{n}h_{\alpha_k-1}(t_{i},t_{j})t_{i_{k+1}}^{\alpha_{k+1}}\cdots t_{i_{l}}^{\alpha_{l}}\cdot x_{j}\partial_{i}x_{i_{k+1}}\partial_{i_{k+1}}\cdots x_{i_{l}}\partial_{i_l}x_i\partial_{j},
\end{gather*}
where $h_{\alpha}$ denotes the $\alpha$-th \emph{complete homogeneous symmetric polynomial}.
This formula can be simplified if we restrict the operator on the right to $V(m\varpi_1)\simeq\mathbb{C}[x_1,\dots,x_n]_{m}$. Note that the differential operator $x_{j}\partial_{i}x_{i_{k+1}}\partial_{i_{k+1}}\cdots x_{i_{l}}\partial_{i_l}x_i\partial_{j}$ acts on the monomial $x^{\mu}=x_1^{\mu_1}\cdots x_n^{\mu_n}$ as follows:
\begin{equation*}
\bigl(x_{j}\partial_{i}x_{i_{k+1}}\partial_{i_{k+1}}\cdots x_{i_{l}}\partial_{i_l}x_i\partial_{j}\bigr)(x^{\mu})=\mu_{j}(\mu_{i}+1-\delta_{ij})\prod_{s=1}^{l}(\mu_s+\delta_{i,i_s}-\delta_{j,i_s})\cdot x^{\mu}.
\end{equation*}
The last equality follows from the identity $x_i\partial_{i}(x^{\mu})=\mu_i\cdot x^{\mu}$. Therefore, the map $\mathfrak{i}_{m}$ sends the element $\DKir(P_{\alpha_1}\cdots P_{\alpha_k})$ of $\mathscr{B}(m\varpi_1)$ to the function in $\Func(\wt(m\varpi_1),S(\mathfrak{h}^{*}))$ whose value at~${\mu\in\wt(m\varpi_1)}$ is equal to
\begin{gather*}
\sum_{k=1}^{l}\prod_{s=1}^{k-1}\bigl(\mu_1t_1^{\alpha_s}+\dots+\mu_nt_n^{\alpha_s}\bigr)
\\
\hphantom{\sum_{k=1}^{l}\prod_{s=1}^{k-1}}{}
\times\sum_{i_{k+1},\dots,i_l=1}^{n}\sum_{i,j=1}^{n}h_{\alpha_k-1}(t_{i},t_{j})t_{i_{k+1}}^{\alpha_{k+1}}\cdots t_{i_{l}}^{\alpha_{l}}\cdot \mu_{j}(\mu_{i}+1-\delta_{ij})\prod_{s=k+1}^{l}(\mu_s+\delta_{i,i_s}-\delta_{j,i_s})
\\
\qquad{}=\sum_{k=1}^{l}\prod_{s=1}^{k-1}(\mu_1t_1^{\alpha_s}+\dots+\mu_nt_n^{\alpha_s})
\\
\qquad\quad{}
\times\Biggl(\sum_{i,j=1}^{n}\mu_{j}(\mu_{i}+1-\delta_{ij})h_{\alpha_k-1}(t_{i},t_{j}) \prod_{s=k+1}^{l}\bigl(\mu_1t_1^{\alpha_s}+\dots+\mu_nt_n^{\alpha_s}+\bigl(t_i^{\alpha_s}-t_j^{\alpha_s}\bigr)\bigr)\Biggr),
\end{gather*}
where we used that for any $s\in\{k+1,\dots,l\}$, we have
\begin{equation*}
\sum_{i_{s}=1}^{n}(\mu_{s}+\delta_{i,i_s}-\delta_{j,i_s})t_{i_s}^{\alpha_s}=\bigl(\mu_1t_1^{\alpha_s}+\mu_2t_2^{\alpha_s}+\dots+\mu_nt_n^{\alpha_s}\bigr) +\bigl(t_{i}^{\alpha_s}-t_j^{\alpha_s}\bigr).
\end{equation*}
Using the fact that $h_{\alpha_k-1}(t_i,t_j)$ equals \smash{$\alpha_kt_i^{\alpha_k-1}$} if $i=j$ and \smash{$\frac{1}{t_i-t_j}\bigl(t_i^{\alpha_k}-t_j^{\alpha_k}\bigr)$} otherwise, we obtain that the value of $\mathfrak{i}_{m}(\DKir(P_{\alpha_1}\dots P_{\alpha_k}))$ at $\mu$ equals
\begin{gather*}
\sum_{k=1}^{l}\Biggl(\sum_{i=1}^{n}\mu_i^2t_i^{\alpha_k-1}\Biggr)\prod_{\substack{1\le s\le l\\s\neq k}}\bigl(\mu_1t_1^{\alpha_s}+\dots+\mu_nt_n^{\alpha_s}\bigr)+\sum_{i\neq j}\frac{\mu_{j}(\mu_{i}+1)}{t_i-t_j}
\\
\quad\! \times\sum_{k=1}^{l}\Biggl(\prod_{s=1}^{k-1}\bigl(\mu_1t_1^{\alpha_s}+\dots+\mu_nt_n^{\alpha_s}\bigr)\! \Biggr)\!\bigl(t_i^{\alpha_k}-t_j^{\alpha_k}\bigr)\! \Biggl(\prod_{s=k+1}^{l}\bigl(\mu_1t_1^{\alpha_s}+\dots+\mu_nt_n^{\alpha_s}+\bigl(t_i^{\alpha_s}-t_j^{\alpha_s}\bigr)\bigr)\!\Biggr)
\\
=\Biggl(\sum_{i=1}^{n}\mu_i\frac{\partial}{\partial t_i}\Biggr)\Biggl(\prod_{s=1}^{l}\bigl(\mu_1t_1^{\alpha_s}+\dots+\mu_nt_n^{\alpha_s}\bigr)\Biggr)
\\
\quad\! +\sum_{i\neq j}\frac{\mu_j(\mu_i+1)}{t_i-t_j}\Biggl(\prod_{s=1}^{l}\bigl(\mu_1t_1^{\alpha_s}+\dots+\mu_nt_n^{\alpha_s}+\bigl(t_i^{\alpha_s}-t_j^{\alpha_s}\bigr)\bigr) -\prod_{s=1}^{l}\bigl(\mu_1t_1^{\alpha_s}+\dots+\mu_nt_n^{\alpha_s}\bigr)\Biggr).
\end{gather*}
Combining everything, we obtain that for \smash{$g=\prod_{k=1}^{l}p_{\alpha_k}$} the operator $\widehat{\DKir}$ maps the function $\mu\mapsto g[\mu_1t_1+\dots+\mu_nt_n]$ to the function
\begin{gather*}
\mu\mapsto
\biggl(\mu_1\frac{\partial}{\partial t_1}+\dots+\mu_n\frac{\partial}{\partial t_n}\biggr)(g[\mu_1t_1+\dots+\mu_nt_n])
\\ \hphantom{\mu\mapsto}{}
+\sum_{\substack{1\le i,j\le n\\i\neq j}}\frac{\mu_j(\mu_i+1)}{t_i-t_j}(g[\mu_1t_1+\dots+\mu_nt_n+(t_i-t_j)]-g[\mu_1t_1+\dots\mu_nt_n]),
\end{gather*}
as claimed.
\end{proof}

\begin{Remark}\label{rem:D_op_inv_big_alg}
In fact, one can modify the argument above in order to show directly that $\DKir$ leaves $\mathscr{B}(m\varpi_1)$ invariant. Namely, one simply needs to verify that for any symmetric functions~$f$~and~$g$ the expression in \eqref{eq:kirillov_wei_sym_func} defines an element of the subalgebra generated by \smash{$\mathcal{F}_{0}^{(m)}$} and~\smash{$\mathcal{F}_{1}^{(m)}$} inside $\Func(\wt(m\varpi_1),S(\mathfrak{h}^*))$. This is essentially done in the proof of Proposition~\ref{prop:D_xy_diag_inv}.
\end{Remark}

\subsubsection{Discussion of the formula} The expression in \eqref{eq:kirillov_wei_sym_func} can be put in a broader context. Consider the ring of \emph{diagonal invariants} $DI_{n}=\mathbb{C}[x_1,\dots,x_n,y_1,\dots,y_n]^{\mathfrak{S}_n}$, where $\mathfrak{S}_n$ acts diagonally on variables $x_i$, $y_i$. It is known that this ring is generated by the so-called \emph{polarized power sums} (see \cite[Chapter~II, Section~A.3]{Weyl}),
\begin{equation*}
p_{a,b}=\sum_{i=1}^{n}x_i^{a} y_i^{b}.
\end{equation*}
Note that the generators $\{p_{a,b}\}_{a,b\ge 0}$ are not algebraically independent.

Let $DI_{n}^{1}$ be the subring of $DI_n$ generated by $p_{a,0}$ and $p_{a,1}$ for all non-negative integers $a$. It~turns out that $DI_{n}^{1}$ is isomorphic to a polynomial ring in $2n$ variables.
\begin{Proposition}\label{prop:DI_n1_gens}
The subring $DI_{n}^{1}$ is a free polynomial ring in $2n$ generators $p_{1,0},\dots,p_{n,0}$ and $p_{0,1},\dots,p_{n-1,1}$.
\end{Proposition}
\begin{proof}
Indeed, these $2n$ polynomials are algebraically independent since the corresponding Jacobian
\begin{equation*}
\frac{\partial(p_{1,0},\dots,p_{n,0},p_{0,1},\dots,p_{n-1,1})}{\partial(x_1,\dots,x_n,y_1,\dots,y_n)}=n!\cdot\prod_{i<j}(x_i-x_j)^2
\end{equation*}
is non-zero. The elements $p_{1,0},\dots,p_{n,0}$ generate all symmetric polynomials in the variables $x_1,\dots,x_n$ and in particular all $p_{a,0}$ for $a>n$. Finally, for any $a\ge n$, we have the identity
\begin{equation*}
p_{a,1}-e_{1}p_{a-1,1}+e_{2}p_{a-1,1}-\dots+(-1)^{n}e_{n}p_{a-n,1}=0,
\end{equation*}
where $e_i=e_i(x_1,\dots,x_n)$ is the $i$-th elementary symmetric polynomial in $x_1,\dots,x_n$. It follows that all elements $p_{a,1}$ belong to the algebra generated by $p_{1,0},\dots,p_{n,0}$ and $p_{0,1},\dots,p_{n-1,1}$. Thus, these $2n$ polynomials freely generate $DI_{n}^{1}$.
\end{proof}

We can now restate the results of the previous subsection.

\begin{Proposition}
The big algebra $\mathscr{B}(S^m(\mathbb{C}^n))$ is a homomorphic image of the ring $DI_{n}^{1}$. Namely, the map
\begin{equation*}
\mathfrak{i}_{m}^{-1}\circ\biggl(\prod_{\mu\in\wt(m\varpi_1)}\mathrm{ev}_{\mu}\biggr)\colon\ DI_{n}^{1}\to\mathscr{B}(S^m(\mathbb{C}^n)),
\end{equation*}
where the ring homomorphism $\mathrm{ev}_{\mu}\colon DI_n^{1}\to \mathbb{C}[t_1,\dots,t_n]$ sends the generator $p_{a,b}$ to \smash{$\sum_{i=1}^{n}\mu_i^bt_i^a$}, is a surjective \smash{$\mathbb{C}[t_1,\dots,t_n]^{\mathfrak{S}_{n}}$}-algebra homomorphism.
\end{Proposition}

Recall that $\Lambda$ is the ring of symmetric functions in infinitely many variables.
For $f\in\Lambda$, we~use the notation $f[A]$ for the \emph{plethystic substitution}, where $A$ is a certain expression in~$x_i$,~$y_i$. In plethystic substitutions, we treat $x_i$ as variables and $y_i$ as constants. It turns out that the expression on the right-hand side of \eqref{eq:kirillov_wei_sym_func} can be lifted to an operator on~$DI_{n}^{1}$.

\begin{Proposition}\label{prop:D_xy_diag_inv}
There exists a $\mathbb{C}$-linear operator \smash{$\widehat{\DKir}_{x,y}\colon DI_{n}^{1}\!\to\!DI_{n}^{1}$} such that for any ${f,g\!\in\!\Lambda}$, we have
\begin{gather}\label{eq:D_xy_formula}
\widehat{\DKir}_{x,y}(f[x_1+\dots+x_n]g[x_1y_1+\dots+x_ny_n])
\\ \nonumber
\qquad{}=\biggl(y_1\frac{\partial}{\partial x_1}+\dots+y_n\frac{\partial}{\partial x_n}\biggr)(f[x_1+\dots+x_n]g[x_1y_1+\dots+x_ny_n])
\\ \nonumber
\qquad\quad{}+f[x_1+\dots+x_n]
\\ \nonumber
\qquad\qquad\times\sum_{\substack{1\le i,j\le n\\i\neq j}}y_j(y_i+1)\frac{g[x_1y_1+\dots+x_ny_n+(x_i-x_j)]-g[x_1y_1+\dots+x_ny_n]}{x_i-x_j}.
\end{gather}
The map $\widehat{\DKir}_{x,y}$ induces the Kirillov--Wei operator $\DKir$ on $\mathscr{B}(S^m(\mathbb{C}^n))$.
\end{Proposition}
\begin{proof}
Denote the expression on the right-hand side of \eqref{eq:D_xy_formula} by $D(f,g)$. It is clear that $D(f,g)\in DI_n$ for any $f,g\in\Lambda$.

Note that the subalgebras $\{f[x_1+\dots+x_n]\mid f\in\Lambda\}$ and $\{g[x_1y_1+\dots+x_ny_n]\mid g\in\Lambda\}$ of $DI_{n}^{1}$ coincide with $\mathbb{C}[p_{1,0},\dots,p_{n,0}]$ and $\mathbb{C}[p_{0,1},\dots,p_{n-1,1}]$, respectively.
It follows now from Proposition~\ref{prop:DI_n1_gens} that one can define (uniquely) a $\mathbb{C}$-linear map \smash{$\widehat{\DKir}_{x,y}\colon DI_{n}^{1}\to DI_{n}$} by the formula
\begin{equation*}
\widehat{\DKir}_{x,y}(f[x_1+\dots+x_n]g[x_1y_1+\dots+x_ny_n])=D(f,g),\qquad\text{where}~f,g\in\Lambda.
\end{equation*}
Observe that the right-hand side depends only on $f[x_1+\dots+x_n]$ and $g[x_1y_1+\dots+x_ny_n]$ rather than on the choice of $f$ and $g$. Indeed, it follows from the fact that if for $f_0,g_0\in\Lambda$ we have $f_0[x_1+\dots+x_n]\equiv 0$ or $g_0[x_1y_1+\dots+x_ny_n]\equiv 0$, then $D(f_0,g_0)\equiv 0$ as well. Therefore, we~obtain a~well-defined $\mathbb{C}$-linear map \smash{$\widehat{\DKir}_{x,y}\colon DI_{n}^{1}\to DI_{n}$}.

To show that \smash{$\widehat{\DKir}_{x,y}$} is an operator acting on \smash{$DI_{n}^{1}$} it remains to verify that $D(f,g)$ is indeed an element of $DI_{n}^{1}$ for any $f,g\in\Lambda$. Since{\samepage
\begin{equation*}
D(f,g)=D(f,1)\cdot g[x_1y_1+\dots+x_ny_n]+f[x_1+\dots+x_n]\cdot D(1,g),
\end{equation*}
it suffices to check that $D(f,1)$ and $D(1,g)$ belong to $DI_{n}^{1}$ for any $f,g\in\Lambda$.}

To compute $D(f,1)$, note that
\begin{equation*}
D(f,1)=\biggl(y_1\frac{\partial}{\partial x_{1}}+\dots+y_n\frac{\partial}{\partial x_{n}}\biggr)(f[x_1+\dots+x_n]),
\end{equation*}
Without loss of generality, we may assume that $f=p_{\alpha_1}\cdots p_{\alpha_l}$ for $\alpha_i\in\{1,\dots,n\}$, i.e., that $f[x_1+\dots+x_n]=p_{\alpha_1,0}\cdots p_{\alpha_l,0}$. Then,
\begin{equation*}
D(f,1)=\biggl(y_1\frac{\partial}{\partial x_{1}}+\dots+y_n\frac{\partial}{\partial x_{n}}\biggr)(p_{\alpha_1,0}\cdots p_{\alpha_l,0})=\sum_{k=1}^{l}\alpha_{k}p_{\alpha_{k}-1,1}\prod_{s\neq k}p_{\alpha_s,0}\in DI_{n}^{1}.
\end{equation*}

Similarly, to check that $D(1,g)\in DI_{n}^{1}$, we may assume that $g=p_{\alpha_1}\cdots p_{\alpha_l}$ for some $\alpha_i\in\{0,1,\dots,n-1\}$, i.e., that $g[x_1y_1+\dots+x_ny_n]=p_{\alpha_1,1}\cdots p_{\alpha_l,1}$. We obtain
\begin{align*}
D(1,g)={}&\Biggl(\sum_{i=1}^{n}y_{i}\frac{\partial}{\partial x_{i}}\Biggr)(p_{\alpha_1,1}\cdots p_{\alpha_l,1})\\
&{}{+}\,\sum_{\substack{1\le i,j\le n\\i\neq j}}\frac{y_j(y_i+1)}{x_i-x_j}\Biggl(\prod_{k=1}^{l}\bigl(p_{\alpha_k,1}+\bigl(x_i^{\alpha_k}-x_j^{\alpha_k}\bigr)\bigr)-\prod_{k=1}^{l}p_{\alpha_k,1}\Biggr)
\\
={}&\sum_{k=1}^{l}\alpha_k\Biggl(\sum_{i=1}^{n}y_i^2x_i^{\alpha_k-1}\Biggr)\prod_{s\neq k}p_{\alpha_s,1}+\sum_{\substack{1\le i,j\le n\\i\neq j}}\frac{y_j(y_i+1)}{x_i-x_j}\sum_{\substack{I\subset[l]\\I\neq\varnothing}}\prod_{s\in I}\bigl(x_{i}^{\alpha_s}-x_{j}^{\alpha_s}\bigr)\prod_{s\notin I}p_{\alpha_s,1}.
\end{align*}
Observe that for any non-empty subset $I$ of $[l]=\{1,\dots,l\}$ the expression \smash{$\frac{1}{x_i-x_j}\prod_{s\in I}(x_{i}^{\alpha_s}-x_{j}^{\alpha_s})$} is in fact a polynomial in $x_i$ and $x_j$. Moreover, if $I$ is not a singleton, then this expression vanishes when $x_i=x_j$. Thus, we can rewrite the formula above as follows:
\begin{gather*}
\widehat{\DKir}_{x,y}(p_{\alpha_1,1}\dots p_{\alpha_l,1})
\\
\quad{}=-\sum_{k=1}^{l}\alpha_k\Biggl(\sum_{i=1}^{n}y_ix_i^{\alpha_k-1}\Biggr)\prod_{s\neq k}p_{\alpha_s,1}+\sum_{i,j= 1}^{n}y_j(y_i+1)\sum_{\substack{I\subset[l]\\I\neq\varnothing}}\biggl(\frac{\prod_{s\in I}\bigl(x_{i}^{\alpha_s}-x_{j}^{\alpha_s}\bigr)}{x_i-x_j}\biggr)\prod_{s\notin I}p_{\alpha_s,1}
\\
\quad{}=-\sum_{k=1}^{l}\alpha_k p_{\alpha_k-1,1}\prod_{s\neq k}p_{\alpha_s,1}+\sum_{i,j= 1}^{n}y_j(y_i+1)\sum_{\substack{I\subset[l]\\I\neq\varnothing}}Q_{I}(x_i,x_j)\prod_{s\notin I}p_{\alpha_s,1},
\end{gather*}
where \smash{$Q_{I}(x_i,x_j)=\frac{1}{x_i-x_j}\prod_{s\in I}(x_{i}^{\alpha_s}-x_{j}^{\alpha_s})$} is a certain polynomial in $x_i$ and $x_j$ whose coefficients depend on $I$. Finally, note that for any polynomial \smash{$Q(u,v)=\sum_{a,b}Q_{a,b}u^av^b$}, we have
\begin{align*}
\sum_{i,j=1}^{n}y_{j}(y_{i}+1)Q(x_i,x_j)&{}=\sum_{a,b}Q_{a,b}\Biggl(\sum_{i=1}^{n}(y_i+1)x_i^{a}\Biggr)\Biggl(\sum_{j=1}^{n}y_jx_j^b\Biggr)\\
&{}=\sum_{a,b}Q_{a,b}(p_{a,0}+p_{a,1})p_{b,1}\in DI_{n}^{1}.
\end{align*}
This completes the proof of the fact that $D(1,g)\in DI_{n}^{1}$. Thus, \smash{$\widehat{\DKir}_{x,y}$} is an operator on $DI_{n}^{1}$ satisfying \eqref{eq:D_xy_formula}.
\end{proof}

\begin{Remark}
Note that a priori it is not clear that \smash{$\widehat{\DKir}_{x,y}$} descends to a well defined map on~$\mathscr{B}(S^m(\mathbb{C}^n))$. It is only a consequence of the computations from Proposition~\ref{prop:kirillov_wei_operator_formula}.
\end{Remark}

\appendix

\section{Proof of Lemma~\ref{lemma:vanishing_lemma}}\label{appendix}

In this appendix, we prove Lemma~\ref{lemma:vanishing_lemma}. In fact, it is a direct consequence of Corollary~\ref{cor:comb_nullstellensatz}. The latter fact is one of the variants of the celebrated combinatorial Nullstellensatz by Alon \cite[Theorems~1.1 and~1.2]{Alon} (see also \cite[Theorem~4]{Karasev_Petrov}). For the sake of completeness, we present here a short proof of this result following \cite{Karasev_Petrov}.

\begin{Lemma}\label{lagrange_interp_formula}
Let $A$ any finite subset of a field $\mathbb{F}$. Then, for any $k\in\{0,1,\dots,|A|-1\}$, we have
\begin{equation*}
\sum_{a\in A}\frac{a^k}{\prod_{b\in A\setminus\{a\}}(b-a)}=
\begin{cases}
0,& k\in\{0,1,\dots,|A|-2\},
\\
1,& k=|A|-1.
\end{cases}
\end{equation*}
\end{Lemma}
\begin{proof}
Applying the Lagrange interpolation formula to the polynomial $x^{k}$ and the elements of~$A$ (note that $k<|A|$) gives the identity
\begin{equation*}
\sum_{a\in A}a^k\cdot\frac{\prod_{b\in A\setminus\{a\}}(x-a)}{\prod_{b\in A\setminus\{a\}}(b-a)}=x^k.
\end{equation*}
Equating the coefficients of $x^{|A|-1}$ on both sides, we obtain the required equality.
\end{proof}

\begin{Proposition}\label{coeff_interp_formula}
Let $g\in\mathbb{F}[X_1,\dots,X_m]$ be a polynomial in $m$ variables over a field $\mathbb{F}$ of total degree $d$. Let $d_1,\dots,d_m$ be non-negative integers such that $d_1+\dots+d_m=d$. Then, for any subsets $A_1,\dots,A_m$ of $\mathbb{F}$ such that $|A_i|=d_i+1$ for all $i$, the $X_1^{d_1}\cdots X_m^{d_m}$ term of $g$ equals
\begin{equation*}
[g]_{X_1^{d_1}\cdots X_m^{d_m}}=\sum_{\substack{(a_1,\dots,a_m)\in\\A_1\times\dots\times A_m}}g(a_1,\dots,a_m)\prod_{i=1}^{m}\frac{1}{\prod_{a\in A_i\setminus\{a_i\}}(a_i-a)}.
\end{equation*}
\end{Proposition}
\begin{proof}
Because of linearity, it suffices to verify the identity in the case when \smash{$g=X_1^{i_1}\cdots X_{m}^{i_m}$} is a monomial such that $i_1+\dots+i_m\le d$. Then, either $(i_1,\dots,i_m)=(d_1,\dots,d_m)$, or $i_j<d_j$ for some $j$. Applying Lemma~\ref{lagrange_interp_formula} $m$ times now yields the required identity.
\end{proof}

\begin{Corollary}\label{cor:comb_nullstellensatz}
Let $g\in\mathbb{F}[X_1,\dots,X_m]$ be a non-zero polynomial in $m$ variables over a field~$\mathbb{F}$. Let $d_1,\dots,d_m$ be non-negative integers such that for each $i\in\{1,\dots,m\}$ the degree of $g$ in the variable $x_i$ is less than $d_i$. Let $A_1,\dots,A_m$ be subsets of $\mathbb{F}$ such that $|A_i|=d_i+1$ for all $i$. Then,~$g$~cannot vanish on $A_1\times\dots\times A_m$.
\end{Corollary}
\begin{proof}
The corollary is a direct consequence of the proposition above. Indeed, assume that $g$ vanishes on $A_1\times\dots\times A_m$ and choose a monomial \smash{$X_{1}^{i_1}\cdots X_m^{i_m}$} of the largest degree which occurs in $g$. Then, by the assumption, we have $0\le i_j\le d_j-1$ for all $j$. Choose arbitrary subsets $B_j\subset A_j$ such that $|B_j|=i_j+1$. Applying Proposition~\ref{coeff_interp_formula} to $g$, the subsets $B_1,\dots,B_m$ and the monomial $X_1^{i_1}\cdots X_m^{i_m}$ yields a contradiction.
\end{proof}

Now Lemma~\ref{lemma:vanishing_lemma} follows from 
Corollary~\ref{cor:comb_nullstellensatz} applied to the field $\mathbb{F}=\mathrm{Frac}(S(\mathfrak{gl}_n^{*}))$.

\subsection*{Acknowledgements}

I~would like to express my gratitude to Tam\'{a}s Hausel for introducing me to the subject and for his constant guidance throughout this work. I~would also like to thank Tam\'{a}s Hausel, Mischa Elkner, Jakub L\"{o}wit, Anton Mellit, Marino Romero, Leonid Rybnikov for many fruitful discussions and feedback on earlier drafts of this paper. We are grateful to the anonymous referees for many useful comments and suggestions that improved the manuscript.

This work was done during the author's PhD studies at the Institute of Science and Technology Austria (ISTA).
The author was supported by the Austrian Science Fund (FWF) grant ``Geometry of the tip of the global nilpotent cone'' no. 10.55776/P35847 and the DOC Fellowship of the Austrian
Academy of Sciences. The author also acknowledges the long-term program of support of the Ukrainian research teams at the Polish Academy of Sciences carried out in collaboration with the U.S. National Academy of Sciences with the financial support of external partners.
For open access purposes, the author has applied a CC BY public copyright license to any author-accepted manuscript version arising from this submission.


\begin{thebibliography}{99}
\footnotesize\itemsep=0pt

\bibitem{Alon}
Alon N., Combinatorial {N}ullstellensatz,
 \href{https://doi.org/10.1017/S0963548398003411}{\textit{Combin. Probab.
 Comput.}} \textbf{8} (1999), 7--29.

\bibitem{Bracken_Green}
Bracken A.J., Green H.S., Vector operators and a~polynomial identity for {${\rm
 SO}(n)$}, \href{https://doi.org/10.1063/1.1665506}{\textit{J.~Math. Phys.}}
 \textbf{12} (1971), 2099--2106.

\bibitem{Capelli}
Capelli A., \"Uber die {Z}ur\"uckf\"uhrung der {C}ayley'schen {O}peration
 {$\Omega$} auf gew\"ohnliche {P}olar-{O}perationen,
 \href{https://doi.org/10.1007/BF01447728}{\textit{Math. Ann.}} \textbf{29}
 (1887), 331--338.

\bibitem{CSS}
Caracciolo S., Sokal A.D., Sportiello A., Noncommutative determinants,
 {C}auchy--{B}inet formulae, and {C}apelli-type identities.~{I}.
 {G}eneralizations of the {C}apelli and {T}urnbull identities,
 \href{https://doi.org/10.37236/192}{\textit{Electron.~J.~Combin.}}
 \textbf{16} (2009), 103, 43~pages,
 \href{http://arxiv.org/abs/0809.3516}{arXiv:0809.3516}.

\bibitem{Feigin_Frenkel}
Feigin B., Frenkel E., Affine {K}ac--{M}oody algebras at the critical level and
 {G}elfand--{D}ikii algebras,
 \href{https://doi.org/10.1142/s0217751x92003781}{\textit{Internat.~J. Modern
 Phys.~A}} \textbf{7} (1992), suppl.~1, 97--215.

\bibitem{Gould}
Gould M.D., Characteristic identities for semisimple {L}ie algebras,
 \href{https://doi.org/10.1017/S0334270000004501}{\textit{J.~Austral. Math.
 Soc. Ser. B}} \textbf{26} (1985), 257--283.

\bibitem{Gould_Zhang_Bracken}
Gould M.D., Zhang R.B., Bracken A.J., Generalized {G}el'fand invariants and
 characteristic identities for quantum groups,
 \href{https://doi.org/10.1063/1.529152}{\textit{J.~Math. Phys.}} \textbf{32}
 (1991), 2298--2303.

\bibitem{Green}
Green H.S., Characteristic identities for generators of {${\rm GL}(n)$, ${\rm
 O}(n)$} and {${\rm SP}(n)$},
 \href{https://doi.org/10.1063/1.1665508}{\textit{J.~Math. Phys.}} \textbf{12}
 (1971), 2106--2113.

\bibitem{Gurevich_Pyatov_Saponov}
Gurevich D.I., Pyatov P.N., Saponov P.A., Hecke symmetries and characteristic
 relations on reflection equation algebras,
 \href{https://doi.org/10.1023/A:1007386006326}{\textit{Lett. Math. Phys.}}
 \textbf{41} (1997), 255--264,
 \href{http://arxiv.org/abs/q-alg/9605048}{arXiv:q-alg/9605048}.

\bibitem{Gurevich_Saponov_Zaitsev}
Gurevich D., Saponov P., Zaitsev M., Wick theorem and matrix {C}apelli identity
 for quantum differential operators on reflection equation algebras,
 \href{http://arxiv.org/abs/2412.13373}{arXiv:2412.13373}.

\bibitem{Hausel_24}
Hausel T., Commutative avatars of representations of semisimple {L}ie groups,
 \href{https://doi.org/10.1073/pnas.2319341121}{\textit{Proc. Natl. Acad. Sci.
 USA}} \textbf{121} (2024), e2319341121, 7~pages,
 \href{http://arxiv.org/abs/2311.02711}{arXiv:2311.02711}.

\bibitem{Hausel_25}
Hausel T., Center of {K}ostant algebra, in Algebraic Geometry and the Langlands
 Program~-- A~Volume in Honour of G\'erard Laumon,
 \url{https://www.laumonvolume.fr/en/articles},
 \href{http://arxiv.org/abs/2509.20159}{arXiv:2509.20159}.

\bibitem{Hausel_Rychlewicz}
Hausel T., Rychlewicz K., Spectrum of equivariant cohomology as a~fixed point
 scheme, \href{https://doi.org/10.46298/epiga.2025.12591}{\textit{\'Epijournal
 G\'eom. Alg\'ebrique}} \textbf{9} (2025), 1, 57~pages,
 \href{http://arxiv.org/abs/2212.11836}{arXiv:2212.11836}.

\bibitem{Howe_89}
Howe R., Remarks on classical invariant theory,
 \href{https://doi.org/10.2307/2001418}{\textit{Trans. Amer. Math. Soc.}}
 \textbf{313} (1989), 539--570.

\bibitem{Howe_95}
Howe R., Perspectives on invariant theory: {S}chur duality, multiplicity-free
 actions and beyond, in The {S}chur Lectures (1992) ({T}el {A}viv),
 \textit{Israel Math. Conf. Proc.}, Vol.~8, Bar-Ilan University, Ramat Gan,
 1995, 1--182.

\bibitem{Howe_Umeda}
Howe R., Umeda T., The {C}apelli identity, the double commutant theorem, and
 multiplicity-free actions,
 \href{https://doi.org/10.1007/BF01459261}{\textit{Math. Ann.}} \textbf{290}
 (1991), 565--619.

\bibitem{Humphreys}
Humphreys J.E., Introduction to {L}ie algebras and representation theory,
 \textit{Grad. Texts in Math.}, Vol.~9,
 \href{https://doi.org/10.1007/978-1-4612-6398-2}{Springer}, New York, 1972.

\bibitem{Ilin_Rybnikov}
Ilin A., Rybnikov L., On classical limits of {B}ethe subalgebras in {Y}angians,
 \href{https://doi.org/10.1007/s00031-021-09648-x}{\textit{Transform. Groups}}
 \textbf{26} (2021), 537--564,
 \href{http://arxiv.org/abs/2009.06934}{arXiv:2009.06934}.

\bibitem{Jing_Liu_Molev}
Jing N., Liu M., Molev A., The {$q$}-immanants and higher quantum {C}apelli
 identities, \href{https://doi.org/10.1007/s00220-025-05273-x}{\textit{Comm.
 Math. Phys.}} \textbf{406} (2025), 99, 16~pages,
 \href{http://arxiv.org/abs/2408.09855}{arXiv:2408.09855}.

\bibitem{Karasev_Petrov}
Karasev R.N., Petrov F.V., Partitions of nonzero elements of a~finite field
 into pairs,
 \href{https://doi.org/10.1007/s11856-012-0020-5}{\textit{Israel~J.~Math.}}
 \textbf{192} (2012), 143--156,
 \href{http://arxiv.org/abs/1005.1177}{arXiv:1005.1177}.

\bibitem{Kirillov_00}
Kirillov A.A., Family algebras,
 \href{https://doi.org/10.1090/S1079-6762-00-00075-5}{\textit{Electron. Res.
 Announc. Amer. Math. Soc.}} \textbf{6} (2000), 7--20.

\bibitem{Kirillov_01}
Kirillov A.A., Introduction to family algebras,
 \href{https://doi.org/10.17323/1609-4514-2001-1-1-49-63}{\textit{Mosc.
 Math.~J.}} \textbf{1} (2001), 49--63.

\bibitem{Kostant_63}
Kostant B., Lie group representations on polynomial rings,
 \href{https://doi.org/10.2307/2373130}{\textit{Amer.~J. Math.}} \textbf{85}
 (1963), 327--404.

\bibitem{Kostant_75}
Kostant B., On the tensor product of a~finite and an infinite dimensional
 representation,
 \href{https://doi.org/10.1016/0022-1236(75)90035-x}{\textit{J.~Funct. Anal.}}
 \textbf{20} (1975), 257--285.

\bibitem{Macdonald}
Macdonald I.G., Symmetric functions and {H}all polynomials, 2nd~ed., \textit{Oxford
 Math. Monogr.}, \href{https://doi.org/10.1093/oso/9780198534891.001.0001}{The
 Clarendon Press}, New York, 1995.

\bibitem{Molev_07}
Molev A., Yangians and classical {L}ie algebras, \textit{Math. Surveys
 Monogr.}, Vol.~143, \href{https://doi.org/10.1090/surv/143}{American
 Mathematical Society}, Providence, RI, 2007.

\bibitem{Molev_18}
Molev A., Sugawara operators for classical {L}ie algebras, \textit{Math.
 Surveys Monogr.}, Vol.~229, \href{https://doi.org/10.1090/surv/229}{American
 Mathematical Society}, Providence, RI, 2018.

\bibitem{Mukhin_Tarasov_Varchenko}
Mukhin E., Tarasov V., Varchenko A., A generalization of the {C}apelli
 identity, in Algebra, Arithmetic, and Geometry: in Honor of
 {Y}u.{I}.~{M}anin. {V}ol.~{II}, \textit{Progr. Math.}, Vol. 270,
 \href{https://doi.org/10.1007/978-0-8176-4747-6_12}{Birkh\"auser},
 Boston, MA, 2009, 383--398,
 \href{http://arxiv.org/abs/math.QA/0610799}{arXiv:math.QA/0610799}.

\bibitem{Nazarov_Olshanski}
Nazarov M., Olshanski G., Bethe subalgebras in twisted {Y}angians,
 \href{https://doi.org/10.1007/BF02099459}{\textit{Comm. Math. Phys.}}
 \textbf{178} (1996), 483--506,
 \href{http://arxiv.org/abs/q-alg/9507003}{arXiv:q-alg/9507003}.

\bibitem{Okounkov}
Okounkov A., Quantum immanants and higher {C}apelli identities,
 \href{https://doi.org/10.1007/BF02587738}{\textit{Transform. Groups}}
 \textbf{1} (1996), 99--126,
 \href{http://arxiv.org/abs/q-alg/9602028}{arXiv:q-alg/9602028}.

\bibitem{Panyushev}
Panyushev D.I., Weight multiplicity free representations,
 {$\mathfrak{g}$}-endomorphism algebras, and {D}ynkin polynomials,
 \href{https://doi.org/10.1112/S0024610703004873}{\textit{J.~London Math.
 Soc.}} \textbf{69} (2004), 273--290,
 \href{http://arxiv.org/abs/math.AG/0112314}{arXiv:math.AG/0112314}.

\bibitem{Rozhkovskaya}
Rozhkovskaya N., Commutativity of quantum family algebras,
 \href{https://doi.org/10.1023/A:1023037100013}{\textit{Lett. Math. Phys.}}
 \textbf{63} (2003), 87--103.

\bibitem{Rybnikov}
Rybnikov L.G., The argument shift method and the {G}audin model,
 \href{https://doi.org/10.1007/s10688-006-0030-3}{\textit{Funct. Anal. Appl.}}
 \textbf{40} (2006), 188--199,
 \href{http://arxiv.org/abs/math.RT/0606380}{arXiv:math.RT/0606380}.

\bibitem{Talalaev}
Talalaev D.V., The quantum {G}audin system,
 \href{https://doi.org/10.1007/s10688-006-0012-5}{\textit{Funct. Anal. Appl.}}
 \textbf{40} (2006), 73--77,
 \href{http://arxiv.org/abs/hep-th/0404153}{arXiv:hep-th/0404153}.

\bibitem{Umeda}
Umeda T., On the proof of the {C}apelli identities,
 \href{https://doi.org/10.1619/fesi.51.1}{\textit{Funkcial. Ekvac.}}
 \textbf{51} (2008), 1--15.

\bibitem{Wei}
Wei Z., The noncommutative {P}oisson bracket and the deformation of the family
 algebras, \href{https://doi.org/10.1063/1.4927337}{\textit{J.~Math. Phys.}}
 \textbf{56} (2015), 071703, 21~pages,
 \href{http://arxiv.org/abs/1211.5865}{arXiv:1211.5865}.

\bibitem{Weyl}
Weyl H., The classical groups. {T}heir invariants and pepresentations,
 Princeton University Press, Princeton, NJ, 1939.

\bibitem{Yakimova}
Yakimova O., Symmetrisation and the {F}eigin--{F}renkel centre,
 \href{https://doi.org/10.1112/S0010437X22007485}{\textit{Compos. Math.}}
 \textbf{158} (2022), 585--622,
 \href{http://arxiv.org/abs/1910.10204}{arXiv:1910.10204}.

\bibitem{Zaitsev}
Zaitsev M.R., Universal matrix {C}apelli identity,
 \href{https://doi.org/10.4213/rm10237e}{\textit{Russian Math. Surveys}}
 \textbf{80} (2025), 338--340,
 \href{http://arxiv.org/abs/2411.13178}{arXiv:2411.13178}.

\bibitem{Zhang_Gould_Bracken}
Zhang R.B., Gould M.D., Bracken A.J., Quantum group invariants and link
 polynomials, \href{https://doi.org/10.1007/BF02099115}{\textit{Comm. Math.
 Phys.}} \textbf{137} (1991), 13--27.

\end{thebibliography}

\pdfbookmark[1]{References}{ref}
\LastPageEnding

\end{document}